
\ifx\shlhetal\undefinedcontrolsequence\let\shlhetal\relax\fi
\def\fmtname{AmS-TeX}

\def\fmtversion{2.1}
\catcode`\@=11
\ifx\amstexloaded@\relax\catcode`\@=\active
  \endinput\else\let\amstexloaded@\relax\fi
\newlinechar=`\^^J
\def\W@{\immediate\write\sixt@@n}
\def\CR@{\W@{^^J\fmtname - Version \fmtversion^^J%
COPYRIGHT 1985, 1990, 1991 - AMERICAN MATHEMATICAL SOCIETY^^J%
Use of this macro package is not restricted provided^^J%
each use is acknowledged upon publication.^^J}}
\CR@ \everyjob{\CR@}
\message{Loading definitions for}
\message{misc utility macros,}
\toksdef\toks@@=2
\long\def\rightappend@#1\to#2{\toks@{\\{#1}}\toks@@
 =\expandafter{#2}\xdef#2{\the\toks@@\the\toks@}\toks@{}\toks@@{}}
\def\alloclist@{}
\newif\ifalloc@
\def\showallocations{{\def\\{\immediate\write\m@ne}\alloclist@}\alloc@true}
\def\alloc@#1#2#3#4#5{\global\advance\count1#1by\@ne
 \ch@ck#1#4#2\allocationnumber=\count1#1
 \global#3#5=\allocationnumber
 \edef\next@{\string#5=\string#2\the\allocationnumber}%
 \expandafter\rightappend@\next@\to\alloclist@}
\newcount\count@@
\newcount\count@@@
\def\FN@{\futurelet\next}
\def\DN@{\def\next@}
\def\DNii@{\def\nextii@}
\def\RIfM@{\relax\ifmmode}
\def\RIfMIfI@{\relax\ifmmode\ifinner}
\def\setboxz@h{\setbox\z@\hbox}
\def\wdz@{\wd\z@}
\def\boxz@{\box\z@}
\def\setbox@ne{\setbox\@ne}
\def\wd@ne{\wd\@ne}
\def\iterate{\body\expandafter\iterate\else\fi}
\def\err@#1{\errmessage{AmS-TeX error: #1}}
\newhelp\defaulthelp@{Sorry, I already gave what help I could...^^J
Maybe you should try asking a human?^^J
An error might have occurred before I noticed any problems.^^J
``If all else fails, read the instructions.''}
\def\Err@{\errhelp\defaulthelp@\err@}
\def\eat@#1{}
\def\in@#1#2{\def\in@@##1#1##2##3\in@@{\ifx\in@##2\in@false\else\in@true\fi}%
 \in@@#2#1\in@\in@@}
\newif\ifin@
\def\space@.{\futurelet\space@\relax}
\space@. %
\newhelp\athelp@
{Only certain combinations beginning with @ make sense to me.^^J
Perhaps you wanted \string\@\space for a printed @?^^J
I've ignored the character or group after @.}
{\catcode`\~=\active 
 \lccode`\~=`\@ \lowercase{\gdef~{\FN@\at@}}}
\def\at@{\let\next@\at@@
 \ifcat\noexpand\next a\else\ifcat\noexpand\next0\else
 \ifcat\noexpand\next\relax\else
   \let\next\at@@@\fi\fi\fi
 \next@}
\def\at@@#1{\expandafter
 \ifx\csname\space @\string#1\endcsname\relax
  \expandafter\at@@@ \else
  \csname\space @\string#1\expandafter\endcsname\fi}
\def\at@@@#1{\errhelp\athelp@ \err@{\Invalid@@ @}}
\def\atdef@#1{\expandafter\def\csname\space @\string#1\endcsname}
\newhelp\defahelp@{If you typed \string\define\space cs instead of
\string\define\string\cs\space^^J
I've substituted an inaccessible control sequence so that your^^J
definition will be completed without mixing me up too badly.^^J
If you typed \string\define{\string\cs} the inaccessible control sequence^^J
was defined to be \string\cs, and the rest of your^^J
definition appears as input.}
\newhelp\defbhelp@{I've ignored your definition, because it might^^J
conflict with other uses that are important to me.}
\def\define{\FN@\define@}
\def\define@{\ifcat\noexpand\next\relax
 \expandafter\define@@\else\errhelp\defahelp@                               
 \err@{\string\define\space must be followed by a control
 sequence}\expandafter\def\expandafter\nextii@\fi}                          
\def\undefined@@@@@@@@@@{}
\def\preloaded@@@@@@@@@@{}
\def\next@@@@@@@@@@{}
\def\define@@#1{\ifx#1\relax\errhelp\defbhelp@                              
 \err@{\string#1\space is already defined}\DN@{\DNii@}\else
 \expandafter\ifx\csname\expandafter\eat@\string                            
 #1@@@@@@@@@@\endcsname\undefined@@@@@@@@@@\errhelp\defbhelp@
 \err@{\string#1\space can't be defined}\DN@{\DNii@}\else
 \expandafter\ifx\csname\expandafter\eat@\string#1\endcsname\relax          
 \global\let#1\undefined\DN@{\def#1}\else\errhelp\defbhelp@
 \err@{\string#1\space is already defined}\DN@{\DNii@}\fi
 \fi\fi\next@}

\def\predefine#1#2{\let#1#2}
\def\undefine#1{\let#1\undefined}
\message{page layout,}
\newdimen\captionwidth@
\captionwidth@\hsize
\advance\captionwidth@-1.5in
\def\pagewidth#1{\hsize#1\relax
 \captionwidth@\hsize\advance\captionwidth@-1.5in}
\def\pageheight#1{\vsize#1\relax}
\def\hcorrection#1{\advance\hoffset#1\relax}
\def\vcorrection#1{\advance\voffset#1\relax}
\message{accents/punctuation,}

\let\graveaccent\`
\let\acuteaccent\'
\let\tildeaccent\~
\let\hataccent\^
\let\underscore\_
\let\B\=
\let\D\.
\let\ic@\/
\def\/{\unskip\ic@}
\def\textfonti{\the\textfont\@ne}
\def\t#1#2{{\edef\next@{\the\font}\textfonti\accent"7F \next@#1#2}}
\def~{\unskip\nobreak\ \ignorespaces}
\def\.{.\spacefactor\@m}
\atdef@;{\leavevmode\null;}
\atdef@:{\leavevmode\null:}
\atdef@?{\leavevmode\null?}
\edef\@{\string @}
\def\relaxnext@{\let\next\relax}
\atdef@-{\relaxnext@\leavevmode
 \DN@{\ifx\next-\DN@-{\FN@\nextii@}\else
  \DN@{\leavevmode\hbox{-}}\fi\next@}%
 \DNii@{\ifx\next-\DN@-{\leavevmode\hbox{---}}\else
  \DN@{\leavevmode\hbox{--}}\fi\next@}%
 \FN@\next@}
\def\srdr@{\kern.16667em}
\def\drsr@{\kern.02778em}
\def\sldl@{\drsr@}
\def\dlsl@{\srdr@}
\atdef@"{\unskip\relaxnext@
 \DN@{\ifx\next\space@\DN@. {\FN@\nextii@}\else
  \DN@.{\FN@\nextii@}\fi\next@.}%
 \DNii@{\ifx\next`\DN@`{\FN@\nextiii@}\else
  \ifx\next\lq\DN@\lq{\FN@\nextiii@}\else
  \DN@####1{\FN@\nextiv@}\fi\fi\next@}%
 \def\nextiii@{\ifx\next`\DN@`{\sldl@``}\else\ifx\next\lq
  \DN@\lq{\sldl@``}\else\DN@{\dlsl@`}\fi\fi\next@}%
 \def\nextiv@{\ifx\next'\DN@'{\srdr@''}\else
  \ifx\next\rq\DN@\rq{\srdr@''}\else\DN@{\drsr@'}\fi\fi\next@}%
 \FN@\next@}

\def\textfontii{\the\textfont\tw@}
\def\lbrace@{\delimiter"4266308 }
\def\rbrace@{\delimiter"5267309 }
\def\{{\RIfM@\lbrace@\else{\textfontii f}\spacefactor\@m\fi}
\def\}{\RIfM@\rbrace@\else
 \let\@sf\empty\ifhmode\edef\@sf{\spacefactor\the\spacefactor}\fi
 {\textfontii g}\@sf\relax\fi}
\let\lbrace\{
\let\rbrace\}
\def\AmSTeX{{\textfontii A\kern-.1667em%
  \lower.5ex\hbox{M}\kern-.125emS}-\TeX}
\message{line and page breaks,}
\def\vmodeerr@#1{\Err@{\string#1\space not allowed between paragraphs}}
\def\mathmodeerr@#1{\Err@{\string#1\space not allowed in math mode}}
\def\linebreak{\RIfM@\mathmodeerr@\linebreak\else
 \ifhmode\unskip\unkern\break\else\vmodeerr@\linebreak\fi\fi}

\newskip\saveskip@
\def\allowlinebreak{\RIfM@\mathmodeerr@\allowlinebreak\else
 \ifhmode\saveskip@\lastskip\unskip
 \allowbreak\ifdim\saveskip@>\z@\hskip\saveskip@\fi
 \else\vmodeerr@\allowlinebreak\fi\fi}
\def\nolinebreak{\RIfM@\mathmodeerr@\nolinebreak\else
 \ifhmode\saveskip@\lastskip\unskip
 \nobreak\ifdim\saveskip@>\z@\hskip\saveskip@\fi
 \else\vmodeerr@\nolinebreak\fi\fi}
\def\newline{\relaxnext@
 \DN@{\RIfM@\expandafter\mathmodeerr@\expandafter\newline\else
  \ifhmode\ifx\next\par\else
  \expandafter\unskip\expandafter\null\expandafter\hfill\expandafter\break\fi
  \else
  \expandafter\vmodeerr@\expandafter\newline\fi\fi}%
 \FN@\next@}
\def\dmatherr@#1{\Err@{\string#1\space not allowed in display math mode}}
\def\nondmatherr@#1{\Err@{\string#1\space not allowed in non-display math
 mode}}
\def\onlydmatherr@#1{\Err@{\string#1\space allowed only in display math mode}}
\def\nonmatherr@#1{\Err@{\string#1\space allowed only in math mode}}
\def\mathbreak{\RIfMIfI@\break\else
 \dmatherr@\mathbreak\fi\else\nonmatherr@\mathbreak\fi}
\def\nomathbreak{\RIfMIfI@\nobreak\else
 \dmatherr@\nomathbreak\fi\else\nonmatherr@\nomathbreak\fi}
\def\allowmathbreak{\RIfMIfI@\allowbreak\else
 \dmatherr@\allowmathbreak\fi\else\nonmatherr@\allowmathbreak\fi}
\def\pagebreak{\RIfM@
 \ifinner\nondmatherr@\pagebreak\else\postdisplaypenalty-\@M\fi
 \else\ifvmode\removelastskip\break\else\vadjust{\break}\fi\fi}
\def\nopagebreak{\RIfM@
 \ifinner\nondmatherr@\nopagebreak\else\postdisplaypenalty\@M\fi
 \else\ifvmode\nobreak\else\vadjust{\nobreak}\fi\fi}
\def\nonvmodeerr@#1{\Err@{\string#1\space not allowed within a paragraph
 or in math}}
\def\vnonvmode@#1#2{\relaxnext@\DNii@{\ifx\next\par\DN@{#1}\else
 \DN@{#2}\fi\next@}%
 \ifvmode\DN@{#1}\else
 \DN@{\FN@\nextii@}\fi\next@}
\def\newpage{\vnonvmode@{\vfill\break}{\nonvmodeerr@\newpage}}
\def\smallpagebreak{\vnonvmode@\smallbreak{\nonvmodeerr@\smallpagebreak}}
\def\medpagebreak{\vnonvmode@\medbreak{\nonvmodeerr@\medpagebreak}}
\def\bigpagebreak{\vnonvmode@\bigbreak{\nonvmodeerr@\bigpagebreak}}
\def\NoBlackBoxes{\global\overfullrule\z@}
\def\BlackBoxes{\global\overfullrule5\p@}
\def\Invalid@#1{\def#1{\Err@{\Invalid@@\string#1}}}
\def\Invalid@@{Invalid use of }
\message{figures,}
\Invalid@\caption
\Invalid@\captionwidth
\newdimen\smallcaptionwidth@
\def\topspace{\mid@false\ins@}
\def\midspace{\mid@true\ins@}
\newif\ifmid@
\def\captionfont@{}
\def\ins@#1{\relaxnext@\allowbreak
 \smallcaptionwidth@\captionwidth@\gdef\thespace@{#1}%
 \DN@{\ifx\next\space@\DN@. {\FN@\nextii@}\else
  \DN@.{\FN@\nextii@}\fi\next@.}%
 \DNii@{\ifx\next\caption\DN@\caption{\FN@\nextiii@}%
  \else\let\next@\nextiv@\fi\next@}%
 \def\nextiv@{\vnonvmode@
  {\ifmid@\expandafter\midinsert\else\expandafter\topinsert\fi
   \vbox to\thespace@{}\endinsert}
  {\ifmid@\nonvmodeerr@\midspace\else\nonvmodeerr@\topspace\fi}}%
 \def\nextiii@{\ifx\next\captionwidth\expandafter\nextv@
  \else\expandafter\nextvi@\fi}%
 \def\nextv@\captionwidth##1##2{\smallcaptionwidth@##1\relax\nextvi@{##2}}%
 \def\nextvi@##1{\def\thecaption@{\captionfont@##1}%
  \DN@{\ifx\next\space@\DN@. {\FN@\nextvii@}\else
   \DN@.{\FN@\nextvii@}\fi\next@.}%
  \FN@\next@}%
 \def\nextvii@{\vnonvmode@
  {\ifmid@\expandafter\midinsert\else
  \expandafter\topinsert\fi\vbox to\thespace@{}\nobreak\smallskip
  \setboxz@h{\noindent\ignorespaces\thecaption@\unskip}%
  \ifdim\wdz@>\smallcaptionwidth@\centerline{\vbox{\hsize\smallcaptionwidth@
   \noindent\ignorespaces\thecaption@\unskip}}%
  \else\centerline{\boxz@}\fi\endinsert}
  {\ifmid@\nonvmodeerr@\midspace
  \else\nonvmodeerr@\topspace\fi}}%
 \FN@\next@}
\message{comments,}
\def\newcodes@{\catcode`\\12\catcode`\{12\catcode`\}12\catcode`\#12%
 \catcode`\%12\relax}
\def\oldcodes@{\catcode`\\0\catcode`\{1\catcode`\}2\catcode`\#6%
 \catcode`\%14\relax}
\def\comment{\newcodes@\endlinechar=10 \comment@}
{\lccode`\0=`\\
\lowercase{\gdef\comment@#1^^J{\comment@@#10endcomment\comment@@@}%
\gdef\comment@@#10endcomment{\FN@\comment@@@}%
\gdef\comment@@@#1\comment@@@{\ifx\next\comment@@@\let\next\comment@
 \else\def\next{\oldcodes@\endlinechar=`\^^M\relax}%
 \fi\next}}}
\def\pr@m@s{\ifx'\next\DN@##1{\prim@s}\else\let\next@\egroup\fi\next@}
\def\prime{{\null\prime@\null}}
\mathchardef\prime@="0230
\let\dsize\displaystyle

\let\ssize\scriptstyle

\message{math spacing,}
\def\,{\RIfM@\mskip\thinmuskip\relax\else\kern.16667em\fi}
\def\!{\RIfM@\mskip-\thinmuskip\relax\else\kern-.16667em\fi}
\let\thinspace\,
\let\negthinspace\!
\def\medspace{\RIfM@\mskip\medmuskip\relax\else\kern.222222em\fi}
\def\negmedspace{\RIfM@\mskip-\medmuskip\relax\else\kern-.222222em\fi}
\def\thickspace{\RIfM@\mskip\thickmuskip\relax\else\kern.27777em\fi}
\let\;\thickspace
\def\negthickspace{\RIfM@\mskip-\thickmuskip\relax\else
 \kern-.27777em\fi}
\atdef@,{\RIfM@\mskip.1\thinmuskip\else\leavevmode\null,\fi}
\atdef@!{\RIfM@\mskip-.1\thinmuskip\else\leavevmode\null!\fi}
\atdef@.{\RIfM@&&\else\leavevmode.\spacefactor3000 \fi}
\def\and{\DOTSB\;\mathchar"3026 \;}

\message{fractions,}
\def\frac#1#2{{#1\over#2}}

\newdimen\ex@
\ex@.2326ex
\Invalid@\thickness
\def\thickfrac{\relaxnext@
 \DN@{\ifx\next\thickness\let\next@\nextii@\else
 \DN@{\nextii@\thickness1}\fi\next@}%
 \DNii@\thickness##1##2##3{{##2\above##1\ex@##3}}%
 \FN@\next@}

\def\thickfracwithdelims#1#2{\relaxnext@\def\ldelim@{#1}\def\rdelim@{#2}%
 \DN@{\ifx\next\thickness\let\next@\nextii@\else
 \DN@{\nextii@\thickness1}\fi\next@}%
 \DNii@\thickness##1##2##3{{##2\abovewithdelims
 \ldelim@\rdelim@##1\ex@##3}}%
 \FN@\next@}

\def\:{\nobreak\hskip.1111em\mathpunct{}\nonscript\mkern-\thinmuskip{:}\hskip
 .3333emplus.0555em\relax}
\def\snug{\unskip\kern-\mathsurround}
\message{smash commands,}
\def\topsmash{\top@true\bot@false\smash@}
\def\botsmash{\top@false\bot@true\smash@}
\newif\iftop@
\newif\ifbot@
\def\smash{\top@true\bot@true\smash@}
\def\smash@{\RIfM@\expandafter\mathpalette\expandafter\mathsm@sh\else
 \expandafter\makesm@sh\fi}
\def\finsm@sh{\iftop@\ht\z@\z@\fi\ifbot@\dp\z@\z@\fi\leavevmode\boxz@}
\message{large operator symbols,}
\def\LimitsOnSums{\global\let\slimits@\displaylimits}
\def\NoLimitsOnSums{\global\let\slimits@\nolimits}
\LimitsOnSums
\mathchardef\coprod@="1360       \def\coprod{\DOTSB\coprod@\slimits@}
\mathchardef\bigvee@="1357       \def\bigvee{\DOTSB\bigvee@\slimits@}
\mathchardef\bigwedge@="1356     \def\bigwedge{\DOTSB\bigwedge@\slimits@}
\mathchardef\biguplus@="1355     \def\biguplus{\DOTSB\biguplus@\slimits@}
\mathchardef\bigcap@="1354       \def\bigcap{\DOTSB\bigcap@\slimits@}
\mathchardef\bigcup@="1353       \def\bigcup{\DOTSB\bigcup@\slimits@}
\mathchardef\prod@="1351         \def\prod{\DOTSB\prod@\slimits@}
\mathchardef\sum@="1350          \def\sum{\DOTSB\sum@\slimits@}
\mathchardef\bigotimes@="134E    \def\bigotimes{\DOTSB\bigotimes@\slimits@}
\mathchardef\bigoplus@="134C     \def\bigoplus{\DOTSB\bigoplus@\slimits@}
\mathchardef\bigodot@="134A      \def\bigodot{\DOTSB\bigodot@\slimits@}
\mathchardef\bigsqcup@="1346     \def\bigsqcup{\DOTSB\bigsqcup@\slimits@}
\message{integrals,}
\def\LimitsOnInts{\global\let\ilimits@\displaylimits}
\def\NoLimitsOnInts{\global\let\ilimits@\nolimits}
\NoLimitsOnInts
\def\int{\DOTSI\intop\ilimits@}
\def\oint{\DOTSI\ointop\ilimits@}
\def\intic@{\mathchoice{\hskip.5em}{\hskip.4em}{\hskip.4em}{\hskip.4em}}
\def\negintic@{\mathchoice
 {\hskip-.5em}{\hskip-.4em}{\hskip-.4em}{\hskip-.4em}}
\def\intkern@{\mathchoice{\!\!\!}{\!\!}{\!\!}{\!\!}}
\def\intdots@{\mathchoice{\plaincdots@}
 {{\cdotp}\mkern1.5mu{\cdotp}\mkern1.5mu{\cdotp}}
 {{\cdotp}\mkern1mu{\cdotp}\mkern1mu{\cdotp}}
 {{\cdotp}\mkern1mu{\cdotp}\mkern1mu{\cdotp}}}
\newcount\intno@
\def\iint{\DOTSI\intno@\tw@\FN@\ints@}
\def\iiint{\DOTSI\intno@\thr@@\FN@\ints@}
\def\iiiint{\DOTSI\intno@4 \FN@\ints@}
\def\idotsint{\DOTSI\intno@\z@\FN@\ints@}
\def\ints@{\findlimits@\ints@@}
\newif\iflimtoken@
\newif\iflimits@
\def\findlimits@{\limtoken@true\ifx\next\limits\limits@true
 \else\ifx\next\nolimits\limits@false\else
 \limtoken@false\ifx\ilimits@\nolimits\limits@false\else
 \ifinner\limits@false\else\limits@true\fi\fi\fi\fi}
\def\multint@{\int\ifnum\intno@=\z@\intdots@                                
 \else\intkern@\fi                                                          
 \ifnum\intno@>\tw@\int\intkern@\fi                                         
 \ifnum\intno@>\thr@@\int\intkern@\fi                                       
 \int}                                                                      
\def\multintlimits@{\intop\ifnum\intno@=\z@\intdots@\else\intkern@\fi
 \ifnum\intno@>\tw@\intop\intkern@\fi
 \ifnum\intno@>\thr@@\intop\intkern@\fi\intop}
\def\ints@@{\iflimtoken@                                                    
 \def\ints@@@{\iflimits@\negintic@\mathop{\intic@\multintlimits@}\limits    
  \else\multint@\nolimits\fi                                                
  \eat@}                                                                    
 \else                                                                      
 \def\ints@@@{\iflimits@\negintic@
  \mathop{\intic@\multintlimits@}\limits\else
  \multint@\nolimits\fi}\fi\ints@@@}
\def\LimitsOnNames{\global\let\nlimits@\displaylimits}
\def\NoLimitsOnNames{\global\let\nlimits@\nolimits@}
\LimitsOnNames
\def\nolimits@{\relaxnext@
 \DN@{\ifx\next\limits\DN@\limits{\nolimits}\else
  \let\next@\nolimits\fi\next@}%
 \FN@\next@}
\message{operator names,}
\def\newmcodes@{\mathcode`\'"27\mathcode`\*"2A\mathcode`\."613A%
 \mathcode`\-"2D\mathcode`\/"2F\mathcode`\:"603A }
\def\operatorname#1{\mathop{\newmcodes@\kern\z@\fam\z@#1}\nolimits@}
\def\operatornamewithlimits#1{\mathop{\newmcodes@\kern\z@\fam\z@#1}\nlimits@}
\def\qopname@#1{\mathop{\fam\z@#1}\nolimits@}
\def\qopnamewl@#1{\mathop{\fam\z@#1}\nlimits@}
\def\arccos{\qopname@{arccos}}
\def\arcsin{\qopname@{arcsin}}
\def\arctan{\qopname@{arctan}}
\def\arg{\qopname@{arg}}
\def\cos{\qopname@{cos}}
\def\cosh{\qopname@{cosh}}
\def\cot{\qopname@{cot}}
\def\coth{\qopname@{coth}}
\def\csc{\qopname@{csc}}
\def\deg{\qopname@{deg}}
\def\det{\qopnamewl@{det}}
\def\dim{\qopname@{dim}}
\def\exp{\qopname@{exp}}
\def\gcd{\qopnamewl@{gcd}}
\def\hom{\qopname@{hom}}
\def\inf{\qopnamewl@{inf}}
\def\injlim{\qopnamewl@{inj\,lim}}
\def\ker{\qopname@{ker}}
\def\lg{\qopname@{lg}}
\def\lim{\qopnamewl@{lim}}
\def\liminf{\qopnamewl@{lim\,inf}}
\def\limsup{\qopnamewl@{lim\,sup}}
\def\ln{\qopname@{ln}}
\def\log{\qopname@{log}}
\def\max{\qopnamewl@{max}}
\def\min{\qopnamewl@{min}}
\def\Pr{\qopnamewl@{Pr}}
\def\projlim{\qopnamewl@{proj\,lim}}
\def\sec{\qopname@{sec}}
\def\sin{\qopname@{sin}}
\def\sinh{\qopname@{sinh}}
\def\sup{\qopnamewl@{sup}}
\def\tan{\qopname@{tan}}
\def\tanh{\qopname@{tanh}}
\def\varinjlim{\mathop{\vtop{\ialign{##\crcr
 \hfil\rm lim\hfil\crcr\noalign{\nointerlineskip}\rightarrowfill\crcr
 \noalign{\nointerlineskip\kern-\ex@}\crcr}}}}
\def\varprojlim{\mathop{\vtop{\ialign{##\crcr
 \hfil\rm lim\hfil\crcr\noalign{\nointerlineskip}\leftarrowfill\crcr
 \noalign{\nointerlineskip\kern-\ex@}\crcr}}}}
\def\varliminf{\mathop{\underline{\vrule height\z@ depth.2exwidth\z@
 \hbox{\rm lim}}}}

\newdimen\buffer@
\buffer@\fontdimen13 \tenex
\newdimen\buffer
\buffer\buffer@

\def\ResetBuffer{\fontdimen13 \tenex\buffer@\global\buffer\buffer@}
\def\shave#1{\mathop{\hbox{$\m@th\fontdimen13 \tenex\z@                     
 \displaystyle{#1}$}}\fontdimen13 \tenex\buffer}

\message{multilevel sub/superscripts,}
\Invalid@\\
\def\Let@{\relax\iffalse{\fi\let\\=\cr\iffalse}\fi}
\Invalid@\vspace
\def\vspace@{\def\vspace##1{\crcr\noalign{\vskip##1\relax}}}
\def\multilimits@{\bgroup\vspace@\Let@
 \baselineskip\fontdimen10 \scriptfont\tw@
 \advance\baselineskip\fontdimen12 \scriptfont\tw@
 \lineskip\thr@@\fontdimen8 \scriptfont\thr@@
 \lineskiplimit\lineskip
 \vbox\bgroup\ialign\bgroup\hfil$\m@th\scriptstyle{##}$\hfil\crcr}
\def\Sb{_\multilimits@}
\def\endSb{\crcr\egroup\egroup\egroup}
\def\Sp{^\multilimits@}

\def\spreadlines#1{\RIfMIfI@\onlydmatherr@\spreadlines\else
 \openup#1\relax\fi\else\onlydmatherr@\spreadlines\fi}
\def\Mathstrut@{\copy\Mathstrutbox@}
\newbox\Mathstrutbox@
\setbox\Mathstrutbox@\null
\setboxz@h{$\m@th($}
\ht\Mathstrutbox@\ht\z@
\dp\Mathstrutbox@\dp\z@
\message{matrices,}
\newdimen\spreadmlines@
\def\spreadmatrixlines#1{\RIfMIfI@
 \onlydmatherr@\spreadmatrixlines\else
 \spreadmlines@#1\relax\fi\else\onlydmatherr@\spreadmatrixlines\fi}
\def\matrix{\null\,\vcenter\bgroup\Let@\vspace@
 \normalbaselines\openup\spreadmlines@\ialign
 \bgroup\hfil$\m@th##$\hfil&&\quad\hfil$\m@th##$\hfil\crcr
 \Mathstrut@\crcr\noalign{\kern-\baselineskip}}
\def\endmatrix{\crcr\Mathstrut@\crcr\noalign{\kern-\baselineskip}\egroup
 \egroup\,}
\def\format{\crcr\egroup\iffalse{\fi\ifnum`}=0 \fi\format@}
\newtoks\hashtoks@
\hashtoks@{#}
\def\format@#1\\{\def\preamble@{#1}%
 \def\l{$\m@th\the\hashtoks@$\hfil}%
 \def\c{\hfil$\m@th\the\hashtoks@$\hfil}%
 \def\r{\hfil$\m@th\the\hashtoks@$}%
 \edef\preamble@@{\preamble@}\ifnum`{=0 \fi\iffalse}\fi
 \ialign\bgroup\span\preamble@@\crcr}
\def\smallmatrix{\null\,\vcenter\bgroup\vspace@\Let@
 \baselineskip9\ex@\lineskip\ex@
 \ialign\bgroup\hfil$\m@th\scriptstyle{##}$\hfil&&\thickspace\hfil
 $\m@th\scriptstyle{##}$\hfil\crcr}
\def\endsmallmatrix{\crcr\egroup\egroup\,}

\newmuskip\dotsspace@
\dotsspace@1.5mu
\def\strip@#1 {#1}
\def\spacehdots#1\for#2{\multispan{#2}\xleaders
 \hbox{$\m@th\mkern\strip@#1 \dotsspace@.\mkern\strip@#1 \dotsspace@$}\hfill}
\def\hdotsfor#1{\spacehdots\@ne\for{#1}}
\def\multispan@#1{\omit\mscount#1\unskip\loop\ifnum\mscount>\@ne\sp@n\repeat}
\def\spaceinnerhdots#1\for#2\after#3{\multispan@{\strip@#2 }#3\xleaders
 \hbox{$\m@th\mkern\strip@#1 \dotsspace@.\mkern\strip@#1 \dotsspace@$}\hfill}
\def\innerhdotsfor#1\after#2{\spaceinnerhdots\@ne\for#1\after{#2}}
\def\cases{\bgroup\spreadmlines@\jot\left\{\,\matrix\format\l&\quad\l\\}
\def\endcases{\endmatrix\right.\egroup}
\message{multiline displays,}
\newif\ifinany@
\newif\ifinalign@
\newif\ifingather@
\def\strut@{\copy\strutbox@}
\newbox\strutbox@
\setbox\strutbox@\hbox{\vrule height8\p@ depth3\p@ width\z@}
\def\topaligned{\null\,\vtop\aligned@}
\def\botaligned{\null\,\vbox\aligned@}
\def\aligned{\null\,\vcenter\aligned@}
\def\aligned@{\bgroup\vspace@\Let@
 \ifinany@\else\openup\jot\fi\ialign
 \bgroup\hfil\strut@$\m@th\displaystyle{##}$&
 $\m@th\displaystyle{{}##}$\hfil\crcr}
\def\endaligned{\crcr\egroup\egroup}

\def\alignedat#1{\null\,\vcenter\bgroup\doat@{#1}\vspace@\Let@
 \ifinany@\else\openup\jot\fi\ialign\bgroup\span\preamble@@\crcr}
\newcount\atcount@
\def\doat@#1{\toks@{\hfil\strut@$\m@th
 \displaystyle{\the\hashtoks@}$&$\m@th\displaystyle
 {{}\the\hashtoks@}$\hfil}
 \atcount@#1\relax\advance\atcount@\m@ne                                    
 \loop\ifnum\atcount@>\z@\toks@=\expandafter{\the\toks@&\hfil$\m@th
 \displaystyle{\the\hashtoks@}$&$\m@th
 \displaystyle{{}\the\hashtoks@}$\hfil}\advance
  \atcount@\m@ne\repeat                                                     
 \xdef\preamble@{\the\toks@}\xdef\preamble@@{\preamble@}}

\def\gathered{\null\,\vcenter\bgroup\vspace@\Let@
 \ifinany@\else\openup\jot\fi\ialign
 \bgroup\hfil\strut@$\m@th\displaystyle{##}$\hfil\crcr}
\def\endgathered{\crcr\egroup\egroup}
\newif\iftagsleft@
\def\TagsOnLeft{\global\tagsleft@true}
\def\TagsOnRight{\global\tagsleft@false}
\TagsOnLeft
\newif\ifmathtags@
\def\TagsAsMath{\global\mathtags@true}
\def\TagsAsText{\global\mathtags@false}
\TagsAsText
\def\tagform@#1{\hbox{\rm(\ignorespaces#1\unskip)}}
\def\thetag{\leavevmode\tagform@}
\def\tag#1$${\iftagsleft@\leqno\else\eqno\fi                                
 \maketag@#1\maketag@                                                       
 $$}                                                                        
\def\maketag@{\FN@\maketag@@}
\def\maketag@@{\ifx\next"\expandafter\maketag@@@\else\expandafter\maketag@@@@
 \fi}
\def\maketag@@@"#1"#2\maketag@{\hbox{\rm#1}}                                
\def\maketag@@@@#1\maketag@{\ifmathtags@\tagform@{$\m@th#1$}\else
 \tagform@{#1}\fi}
\interdisplaylinepenalty\@M
\def\allowdisplaybreaks{\RIfMIfI@
 \onlydmatherr@\allowdisplaybreaks\else
 \interdisplaylinepenalty\z@\fi\else\onlydmatherr@\allowdisplaybreaks\fi}
\Invalid@\allowdisplaybreak
\Invalid@\displaybreak
\Invalid@\intertext
\def\allowdisplaybreak@{\def\allowdisplaybreak{\crcr\noalign{\allowbreak}}}
\def\displaybreak@{\def\displaybreak{\crcr\noalign{\break}}}
\def\intertext@{\def\intertext##1{\crcr\noalign{%
 \penalty\postdisplaypenalty \vskip\belowdisplayskip
 \vbox{\normalbaselines\noindent##1}%
 \penalty\predisplaypenalty \vskip\abovedisplayskip}}}
\newskip\centering@
\centering@\z@ plus\@m\p@
\def\align{\relax\ifingather@\DN@{\csname align (in
  \string\gather)\endcsname}\else
 \ifmmode\ifinner\DN@{\onlydmatherr@\align}\else
  \let\next@\align@\fi
 \else\DN@{\onlydmatherr@\align}\fi\fi\next@}
\newhelp\andhelp@
{An extra & here is so disastrous that you should probably exit^^J
and fix things up.}
\newif\iftag@
\newcount\and@
\def\align@{\inalign@true\inany@true
 \vspace@\allowdisplaybreak@\displaybreak@\intertext@
 \def\tag{\global\tag@true\ifnum\and@=\z@\DN@{&&}\else
          \DN@{&}\fi\next@}%
 \iftagsleft@\DN@{\csname align \endcsname}\else
  \DN@{\csname align \space\endcsname}\fi\next@}
\def\Tag@{\iftag@\else\errhelp\andhelp@\err@{Extra & on this line}\fi}
\newdimen\lwidth@
\newdimen\rwidth@
\newdimen\maxlwidth@
\newdimen\maxrwidth@
\newdimen\totwidth@
\def\measure@#1\endalign{\lwidth@\z@\rwidth@\z@\maxlwidth@\z@\maxrwidth@\z@
 \global\and@\z@                                                            
 \setbox@ne\vbox                                                            
  {\everycr{\noalign{\global\tag@false\global\and@\z@}}\Let@                
  \halign{\setboxz@h{$\m@th\displaystyle{\@lign##}$}
   \global\lwidth@\wdz@                                                     
   \ifdim\lwidth@>\maxlwidth@\global\maxlwidth@\lwidth@\fi                  
   \global\advance\and@\@ne                                                 
   &\setboxz@h{$\m@th\displaystyle{{}\@lign##}$}\global\rwidth@\wdz@        
   \ifdim\rwidth@>\maxrwidth@\global\maxrwidth@\rwidth@\fi                  
   \global\advance\and@\@ne                                                
   &\Tag@
   \eat@{##}\crcr#1\crcr}}
 \totwidth@\maxlwidth@\advance\totwidth@\maxrwidth@}                       
\def\displ@y@{\global\dt@ptrue\openup\jot
 \everycr{\noalign{\global\tag@false\global\and@\z@\ifdt@p\global\dt@pfalse
 \vskip-\lineskiplimit\vskip\normallineskiplimit\else
 \penalty\interdisplaylinepenalty\fi}}}
\def\black@#1{\noalign{\ifdim#1>\displaywidth
 \dimen@\prevdepth\nointerlineskip                                          
 \vskip-\ht\strutbox@\vskip-\dp\strutbox@                                   
 \vbox{\noindent\hbox to#1{\strut@\hfill}}
 \prevdepth\dimen@                                                          
 \fi}}
\expandafter\def\csname align \space\endcsname#1\endalign
 {\measure@#1\endalign\global\and@\z@                                       
 \ifingather@\everycr{\noalign{\global\and@\z@}}\else\displ@y@\fi           
 \Let@\tabskip\centering@                                                   
 \halign to\displaywidth
  {\hfil\strut@\setboxz@h{$\m@th\displaystyle{\@lign##}$}
  \global\lwidth@\wdz@\boxz@\global\advance\and@\@ne                        
  \tabskip\z@skip                                                           
  &\setboxz@h{$\m@th\displaystyle{{}\@lign##}$}
  \global\rwidth@\wdz@\boxz@\hfill\global\advance\and@\@ne                  
  \tabskip\centering@                                                       
  &\setboxz@h{\@lign\strut@\maketag@##\maketag@}
  \dimen@\displaywidth\advance\dimen@-\totwidth@
  \divide\dimen@\tw@\advance\dimen@\maxrwidth@\advance\dimen@-\rwidth@     
  \ifdim\dimen@<\tw@\wdz@\llap{\vtop{\normalbaselines\null\boxz@}}
  \else\llap{\boxz@}\fi                                                    
  \tabskip\z@skip                                                          
  \crcr#1\crcr                                                             
  \black@\totwidth@}}                                                      
\newdimen\lineht@
\expandafter\def\csname align \endcsname#1\endalign{\measure@#1\endalign
 \global\and@\z@
 \ifdim\totwidth@>\displaywidth\let\displaywidth@\totwidth@\else
  \let\displaywidth@\displaywidth\fi                                        
 \ifingather@\everycr{\noalign{\global\and@\z@}}\else\displ@y@\fi
 \Let@\tabskip\centering@\halign to\displaywidth
  {\hfil\strut@\setboxz@h{$\m@th\displaystyle{\@lign##}$}%
  \global\lwidth@\wdz@\global\lineht@\ht\z@                                 
  \boxz@\global\advance\and@\@ne
  \tabskip\z@skip&\setboxz@h{$\m@th\displaystyle{{}\@lign##}$}%
  \global\rwidth@\wdz@\ifdim\ht\z@>\lineht@\global\lineht@\ht\z@\fi         
  \boxz@\hfil\global\advance\and@\@ne
  \tabskip\centering@&\kern-\displaywidth@                                  
  \setboxz@h{\@lign\strut@\maketag@##\maketag@}%
  \dimen@\displaywidth\advance\dimen@-\totwidth@
  \divide\dimen@\tw@\advance\dimen@\maxlwidth@\advance\dimen@-\lwidth@
  \ifdim\dimen@<\tw@\wdz@
   \rlap{\vbox{\normalbaselines\boxz@\vbox to\lineht@{}}}\else
   \rlap{\boxz@}\fi
  \tabskip\displaywidth@\crcr#1\crcr\black@\totwidth@}}
\expandafter\def\csname align (in \string\gather)\endcsname
  #1\endalign{\vcenter{\align@#1\endalign}}
\Invalid@\endalign
\newif\ifxat@
\def\alignat{\RIfMIfI@\DN@{\onlydmatherr@\alignat}\else
 \DN@{\csname alignat \endcsname}\fi\else
 \DN@{\onlydmatherr@\alignat}\fi\next@}
\newif\ifmeasuring@
\newbox\savealignat@
\expandafter\def\csname alignat \endcsname#1#2\endalignat                   
 {\inany@true\xat@false
 \def\tag{\global\tag@true\count@#1\relax\multiply\count@\tw@
  \xdef\tag@{}\loop\ifnum\count@>\and@\xdef\tag@{&\tag@}\advance\count@\m@ne
  \repeat\tag@}%
 \vspace@\allowdisplaybreak@\displaybreak@\intertext@
 \displ@y@\measuring@true                                                   
 \setbox\savealignat@\hbox{$\m@th\displaystyle\Let@
  \attag@{#1}
  \vbox{\halign{\span\preamble@@\crcr#2\crcr}}$}%
 \measuring@false                                                           
 \Let@\attag@{#1}
 \tabskip\centering@\halign to\displaywidth
  {\span\preamble@@\crcr#2\crcr                                             
  \black@{\wd\savealignat@}}}                                               
\Invalid@\endalignat
\def\xalignat{\RIfMIfI@
 \DN@{\onlydmatherr@\xalignat}\else
 \DN@{\csname xalignat \endcsname}\fi\else
 \DN@{\onlydmatherr@\xalignat}\fi\next@}
\expandafter\def\csname xalignat \endcsname#1#2\endxalignat
 {\inany@true\xat@true
 \def\tag{\global\tag@true\def\tag@{}\count@#1\relax\multiply\count@\tw@
  \loop\ifnum\count@>\and@\xdef\tag@{&\tag@}\advance\count@\m@ne\repeat\tag@}%
 \vspace@\allowdisplaybreak@\displaybreak@\intertext@
 \displ@y@\measuring@true\setbox\savealignat@\hbox{$\m@th\displaystyle\Let@
 \attag@{#1}\vbox{\halign{\span\preamble@@\crcr#2\crcr}}$}%
 \measuring@false\Let@
 \attag@{#1}\tabskip\centering@\halign to\displaywidth
 {\span\preamble@@\crcr#2\crcr\black@{\wd\savealignat@}}}
\def\attag@#1{\let\Maketag@\maketag@\let\TAG@\Tag@                          
 \let\Tag@=0\let\maketag@=0
 \ifmeasuring@\def\llap@##1{\setboxz@h{##1}\hbox to\tw@\wdz@{}}%
  \def\rlap@##1{\setboxz@h{##1}\hbox to\tw@\wdz@{}}\else
  \let\llap@\llap\let\rlap@\rlap\fi                                         
 \toks@{\hfil\strut@$\m@th\displaystyle{\@lign\the\hashtoks@}$\tabskip\z@skip
  \global\advance\and@\@ne&$\m@th\displaystyle{{}\@lign\the\hashtoks@}$\hfil
  \ifxat@\tabskip\centering@\fi\global\advance\and@\@ne}
 \iftagsleft@
  \toks@@{\tabskip\centering@&\Tag@\kern-\displaywidth
   \rlap@{\@lign\maketag@\the\hashtoks@\maketag@}%
   \global\advance\and@\@ne\tabskip\displaywidth}\else
  \toks@@{\tabskip\centering@&\Tag@\llap@{\@lign\maketag@
   \the\hashtoks@\maketag@}\global\advance\and@\@ne\tabskip\z@skip}\fi      
 \atcount@#1\relax\advance\atcount@\m@ne
 \loop\ifnum\atcount@>\z@
 \toks@=\expandafter{\the\toks@&\hfil$\m@th\displaystyle{\@lign
  \the\hashtoks@}$\global\advance\and@\@ne
  \tabskip\z@skip&$\m@th\displaystyle{{}\@lign\the\hashtoks@}$\hfil\ifxat@
  \tabskip\centering@\fi\global\advance\and@\@ne}\advance\atcount@\m@ne
 \repeat                                                                    
 \xdef\preamble@{\the\toks@\the\toks@@}
 \xdef\preamble@@{\preamble@}
 \let\maketag@\Maketag@\let\Tag@\TAG@}                                      
\Invalid@\endxalignat
\def\xxalignat{\RIfMIfI@
 \DN@{\onlydmatherr@\xxalignat}\else\DN@{\csname xxalignat
  \endcsname}\fi\else
 \DN@{\onlydmatherr@\xxalignat}\fi\next@}
\expandafter\def\csname xxalignat \endcsname#1#2\endxxalignat{\inany@true
 \vspace@\allowdisplaybreak@\displaybreak@\intertext@
 \displ@y\setbox\savealignat@\hbox{$\m@th\displaystyle\Let@
 \xxattag@{#1}\vbox{\halign{\span\preamble@@\crcr#2\crcr}}$}%
 \Let@\xxattag@{#1}\tabskip\z@skip\halign to\displaywidth
 {\span\preamble@@\crcr#2\crcr\black@{\wd\savealignat@}}}
\def\xxattag@#1{\toks@{\tabskip\z@skip\hfil\strut@
 $\m@th\displaystyle{\the\hashtoks@}$&%
 $\m@th\displaystyle{{}\the\hashtoks@}$\hfil\tabskip\centering@&}%
 \atcount@#1\relax\advance\atcount@\m@ne\loop\ifnum\atcount@>\z@
 \toks@=\expandafter{\the\toks@&\hfil$\m@th\displaystyle{\the\hashtoks@}$%
  \tabskip\z@skip&$\m@th\displaystyle{{}\the\hashtoks@}$\hfil
  \tabskip\centering@}\advance\atcount@\m@ne\repeat
 \xdef\preamble@{\the\toks@\tabskip\z@skip}\xdef\preamble@@{\preamble@}}
\Invalid@\endxxalignat
\newdimen\gwidth@
\newdimen\gmaxwidth@
\def\gmeasure@#1\endgather{\gwidth@\z@\gmaxwidth@\z@\setbox@ne\vbox{\Let@
 \halign{\setboxz@h{$\m@th\displaystyle{##}$}\global\gwidth@\wdz@
 \ifdim\gwidth@>\gmaxwidth@\global\gmaxwidth@\gwidth@\fi
 &\eat@{##}\crcr#1\crcr}}}
\def\gather{\RIfMIfI@\DN@{\onlydmatherr@\gather}\else
 \ingather@true\inany@true\def\tag{&}%
 \vspace@\allowdisplaybreak@\displaybreak@\intertext@
 \displ@y\Let@
 \iftagsleft@\DN@{\csname gather \endcsname}\else
  \DN@{\csname gather \space\endcsname}\fi\fi
 \else\DN@{\onlydmatherr@\gather}\fi\next@}
\expandafter\def\csname gather \space\endcsname#1\endgather
 {\gmeasure@#1\endgather\tabskip\centering@
 \halign to\displaywidth{\hfil\strut@\setboxz@h{$\m@th\displaystyle{##}$}%
 \global\gwidth@\wdz@\boxz@\hfil&
 \setboxz@h{\strut@{\maketag@##\maketag@}}%
 \dimen@\displaywidth\advance\dimen@-\gwidth@
 \ifdim\dimen@>\tw@\wdz@\llap{\boxz@}\else
 \llap{\vtop{\normalbaselines\null\boxz@}}\fi
 \tabskip\z@skip\crcr#1\crcr\black@\gmaxwidth@}}
\newdimen\glineht@
\expandafter\def\csname gather \endcsname#1\endgather{\gmeasure@#1\endgather
 \ifdim\gmaxwidth@>\displaywidth\let\gdisplaywidth@\gmaxwidth@\else
 \let\gdisplaywidth@\displaywidth\fi\tabskip\centering@\halign to\displaywidth
 {\hfil\strut@\setboxz@h{$\m@th\displaystyle{##}$}%
 \global\gwidth@\wdz@\global\glineht@\ht\z@\boxz@\hfil&\kern-\gdisplaywidth@
 \setboxz@h{\strut@{\maketag@##\maketag@}}%
 \dimen@\displaywidth\advance\dimen@-\gwidth@
 \ifdim\dimen@>\tw@\wdz@\rlap{\boxz@}\else
 \rlap{\vbox{\normalbaselines\boxz@\vbox to\glineht@{}}}\fi
 \tabskip\gdisplaywidth@\crcr#1\crcr\black@\gmaxwidth@}}
\newif\ifctagsplit@
\def\CenteredTagsOnSplits{\global\ctagsplit@true}
\def\TopOrBottomTagsOnSplits{\global\ctagsplit@false}
\TopOrBottomTagsOnSplits
\def\split{\relax\ifinany@\let\next@\insplit@\else
 \ifmmode\ifinner\def\next@{\onlydmatherr@\split}\else
 \let\next@\outsplit@\fi\else
 \def\next@{\onlydmatherr@\split}\fi\fi\next@}
\def\insplit@{\global\setbox\z@\vbox\bgroup\vspace@\Let@\ialign\bgroup
 \hfil\strut@$\m@th\displaystyle{##}$&$\m@th\displaystyle{{}##}$\hfill\crcr}
\def\endsplit{\crcr\egroup\egroup\iftagsleft@\expandafter\lendsplit@\else
 \expandafter\rendsplit@\fi}
\def\rendsplit@{\global\setbox9 \vbox
 {\unvcopy\z@\global\setbox8 \lastbox\unskip}
 \setbox@ne\hbox{\unhcopy8 \unskip\global\setbox\tw@\lastbox
 \unskip\global\setbox\thr@@\lastbox}
 \global\setbox7 \hbox{\unhbox\tw@\unskip}
 \ifinalign@\ifctagsplit@                                                   
  \gdef\split@{\hbox to\wd\thr@@{}&
   \vcenter{\vbox{\moveleft\wd\thr@@\boxz@}}}
 \else\gdef\split@{&\vbox{\moveleft\wd\thr@@\box9}\crcr
  \box\thr@@&\box7}\fi                                                      
 \else                                                                      
  \ifctagsplit@\gdef\split@{\vcenter{\boxz@}}\else
  \gdef\split@{\box9\crcr\hbox{\box\thr@@\box7}}\fi
 \fi
 \split@}                                                                   
\def\lendsplit@{\global\setbox9\vtop{\unvcopy\z@}
 \setbox@ne\vbox{\unvcopy\z@\global\setbox8\lastbox}
 \setbox@ne\hbox{\unhcopy8\unskip\setbox\tw@\lastbox
  \unskip\global\setbox\thr@@\lastbox}
 \ifinalign@\ifctagsplit@                                                   
  \gdef\split@{\hbox to\wd\thr@@{}&
  \vcenter{\vbox{\moveleft\wd\thr@@\box9}}}
  \else                                                                     
  \gdef\split@{\hbox to\wd\thr@@{}&\vbox{\moveleft\wd\thr@@\box9}}\fi
 \else
  \ifctagsplit@\gdef\split@{\vcenter{\box9}}\else
  \gdef\split@{\box9}\fi
 \fi\split@}
\def\outsplit@#1$${\align\insplit@#1\endalign$$}
\newdimen\multlinegap@
\multlinegap@1em
\newdimen\multlinetaggap@
\multlinetaggap@1em
\def\MultlineGap#1{\global\multlinegap@#1\relax}
\def\multlinegap#1{\RIfMIfI@\onlydmatherr@\multlinegap\else
 \multlinegap@#1\relax\fi\else\onlydmatherr@\multlinegap\fi}
\def\nomultlinegap{\multlinegap{\z@}}
\def\multline{\RIfMIfI@
 \DN@{\onlydmatherr@\multline}\else
 \DN@{\multline@}\fi\else
 \DN@{\onlydmatherr@\multline}\fi\next@}
\newif\iftagin@
\def\tagin@#1{\tagin@false\in@\tag{#1}\ifin@\tagin@true\fi}
\def\multline@#1$${\inany@true\vspace@\allowdisplaybreak@\displaybreak@
 \tagin@{#1}\iftagsleft@\DN@{\multline@l#1$$}\else
 \DN@{\multline@r#1$$}\fi\next@}
\newdimen\mwidth@
\def\rmmeasure@#1\endmultline{%
 \def\shoveleft##1{##1}\def\shoveright##1{##1}
 \setbox@ne\vbox{\Let@\halign{\setboxz@h
  {$\m@th\@lign\displaystyle{}##$}\global\mwidth@\wdz@
  \crcr#1\crcr}}}
\newdimen\mlineht@
\newif\ifzerocr@
\newif\ifonecr@
\def\lmmeasure@#1\endmultline{\global\zerocr@true\global\onecr@false
 \everycr{\noalign{\ifonecr@\global\onecr@false\fi
  \ifzerocr@\global\zerocr@false\global\onecr@true\fi}}
  \def\shoveleft##1{##1}\def\shoveright##1{##1}%
 \setbox@ne\vbox{\Let@\halign{\setboxz@h
  {$\m@th\@lign\displaystyle{}##$}\ifonecr@\global\mwidth@\wdz@
  \global\mlineht@\ht\z@\fi\crcr#1\crcr}}}
\newbox\mtagbox@
\newdimen\ltwidth@
\newdimen\rtwidth@
\def\multline@l#1$${\iftagin@\DN@{\lmultline@@#1$$}\else
 \DN@{\setbox\mtagbox@\null\ltwidth@\z@\rtwidth@\z@
  \lmultline@@@#1$$}\fi\next@}
\def\lmultline@@#1\endmultline\tag#2$${%
 \setbox\mtagbox@\hbox{\maketag@#2\maketag@}
 \lmmeasure@#1\endmultline\dimen@\mwidth@\advance\dimen@\wd\mtagbox@
 \advance\dimen@\multlinetaggap@                                            
 \ifdim\dimen@>\displaywidth\ltwidth@\z@\else\ltwidth@\wd\mtagbox@\fi       
 \lmultline@@@#1\endmultline$$}
\def\lmultline@@@{\displ@y
 \def\shoveright##1{##1\hfilneg\hskip\multlinegap@}%
 \def\shoveleft##1{\setboxz@h{$\m@th\displaystyle{}##1$}%
  \setbox@ne\hbox{$\m@th\displaystyle##1$}%
  \hfilneg
  \iftagin@
   \ifdim\ltwidth@>\z@\hskip\ltwidth@\hskip\multlinetaggap@\fi
  \else\hskip\multlinegap@\fi\hskip.5\wd@ne\hskip-.5\wdz@##1}
  \halign\bgroup\Let@\hbox to\displaywidth
   {\strut@$\m@th\displaystyle\hfil{}##\hfil$}\crcr
   \hfilneg                                                                 
   \iftagin@                                                                
    \ifdim\ltwidth@>\z@                                                     
     \box\mtagbox@\hskip\multlinetaggap@                                    
    \else
     \rlap{\vbox{\normalbaselines\hbox{\strut@\box\mtagbox@}%
     \vbox to\mlineht@{}}}\fi                                               
   \else\hskip\multlinegap@\fi}                                             
\def\multline@r#1$${\iftagin@\DN@{\rmultline@@#1$$}\else
 \DN@{\setbox\mtagbox@\null\ltwidth@\z@\rtwidth@\z@
  \rmultline@@@#1$$}\fi\next@}
\def\rmultline@@#1\endmultline\tag#2$${\ltwidth@\z@
 \setbox\mtagbox@\hbox{\maketag@#2\maketag@}%
 \rmmeasure@#1\endmultline\dimen@\mwidth@\advance\dimen@\wd\mtagbox@
 \advance\dimen@\multlinetaggap@
 \ifdim\dimen@>\displaywidth\rtwidth@\z@\else\rtwidth@\wd\mtagbox@\fi
 \rmultline@@@#1\endmultline$$}
\def\rmultline@@@{\displ@y
 \def\shoveright##1{##1\hfilneg\iftagin@\ifdim\rtwidth@>\z@
  \hskip\rtwidth@\hskip\multlinetaggap@\fi\else\hskip\multlinegap@\fi}%
 \def\shoveleft##1{\setboxz@h{$\m@th\displaystyle{}##1$}%
  \setbox@ne\hbox{$\m@th\displaystyle##1$}%
  \hfilneg\hskip\multlinegap@\hskip.5\wd@ne\hskip-.5\wdz@##1}%
 \halign\bgroup\Let@\hbox to\displaywidth
  {\strut@$\m@th\displaystyle\hfil{}##\hfil$}\crcr
 \hfilneg\hskip\multlinegap@}
\def\endmultline{\iftagsleft@\expandafter\lendmultline@\else
 \expandafter\rendmultline@\fi}
\def\lendmultline@{\hfilneg\hskip\multlinegap@\crcr\egroup}
\def\rendmultline@{\iftagin@                                                
 \ifdim\rtwidth@>\z@                                                        
  \hskip\multlinetaggap@\box\mtagbox@                                       
 \else\llap{\vtop{\normalbaselines\null\hbox{\strut@\box\mtagbox@}}}\fi     
 \else\hskip\multlinegap@\fi                                                
 \hfilneg\crcr\egroup}
\def\bmod{\mskip-\medmuskip\mkern5mu\mathbin{\fam\z@ mod}\penalty900
 \mkern5mu\mskip-\medmuskip}
\def\pmod#1{\allowbreak\ifinner\mkern8mu\else\mkern18mu\fi
 ({\fam\z@ mod}\,\,#1)}
\def\pod#1{\allowbreak\ifinner\mkern8mu\else\mkern18mu\fi(#1)}
\def\mod#1{\allowbreak\ifinner\mkern12mu\else\mkern18mu\fi{\fam\z@ mod}\,\,#1}
\message{continued fractions,}
\newcount\cfraccount@
\def\cfrac{\bgroup\bgroup\advance\cfraccount@\@ne\strut
 \iffalse{\fi\def\\{\over\displaystyle}\iffalse}\fi}
\def\lcfrac{\bgroup\bgroup\advance\cfraccount@\@ne\strut
 \iffalse{\fi\def\\{\hfill\over\displaystyle}\iffalse}\fi}
\def\rcfrac{\bgroup\bgroup\advance\cfraccount@\@ne\strut\hfill
 \iffalse{\fi\def\\{\over\displaystyle}\iffalse}\fi}
\def\gloop@#1\repeat{\gdef\body{#1}\iterate}
\def\endcfrac{\gloop@\ifnum\cfraccount@>\z@\global\advance\cfraccount@\m@ne
 \egroup\hskip-\nulldelimiterspace\egroup\repeat}
\message{compound symbols,}
\def\binrel@#1{\setboxz@h{\thinmuskip0mu
  \medmuskip\m@ne mu\thickmuskip\@ne mu$#1\m@th$}%
 \setbox@ne\hbox{\thinmuskip0mu\medmuskip\m@ne mu\thickmuskip
  \@ne mu${}#1{}\m@th$}%
 \setbox\tw@\hbox{\hskip\wd@ne\hskip-\wdz@}}
\def\overset#1\to#2{\binrel@{#2}\ifdim\wd\tw@<\z@
 \mathbin{\mathop{\kern\z@#2}\limits^{#1}}\else\ifdim\wd\tw@>\z@
 \mathrel{\mathop{\kern\z@#2}\limits^{#1}}\else
 {\mathop{\kern\z@#2}\limits^{#1}}{}\fi\fi}
\def\underset#1\to#2{\binrel@{#2}\ifdim\wd\tw@<\z@
 \mathbin{\mathop{\kern\z@#2}\limits_{#1}}\else\ifdim\wd\tw@>\z@
 \mathrel{\mathop{\kern\z@#2}\limits_{#1}}\else
 {\mathop{\kern\z@#2}\limits_{#1}}{}\fi\fi}
\def\oversetbrace#1\to#2{\overbrace{#2}^{#1}}
\def\undersetbrace#1\to#2{\underbrace{#2}_{#1}}
\def\sideset#1\and#2\to#3{%
 \setbox@ne\hbox{$\dsize{\vphantom{#3}}#1{#3}\m@th$}%
 \setbox\tw@\hbox{$\dsize{#3}#2\m@th$}%
 \hskip\wd@ne\hskip-\wd\tw@\mathop{\hskip\wd\tw@\hskip-\wd@ne
  {\vphantom{#3}}#1{#3}#2}}
\def\rightarrowfill@#1{\setboxz@h{$#1-\m@th$}\ht\z@\z@
  $#1\m@th\copy\z@\mkern-6mu\cleaders
  \hbox{$#1\mkern-2mu\box\z@\mkern-2mu$}\hfill
  \mkern-6mu\mathord\rightarrow$}
\def\leftarrowfill@#1{\setboxz@h{$#1-\m@th$}\ht\z@\z@
  $#1\m@th\mathord\leftarrow\mkern-6mu\cleaders
  \hbox{$#1\mkern-2mu\copy\z@\mkern-2mu$}\hfill
  \mkern-6mu\box\z@$}
\def\leftrightarrowfill@#1{\setboxz@h{$#1-\m@th$}\ht\z@\z@
  $#1\m@th\mathord\leftarrow\mkern-6mu\cleaders
  \hbox{$#1\mkern-2mu\box\z@\mkern-2mu$}\hfill
  \mkern-6mu\mathord\rightarrow$}
\def\overrightarrow{\mathpalette\overrightarrow@}
\def\overrightarrow@#1#2{\vbox{\ialign{##\crcr\rightarrowfill@#1\crcr
 \noalign{\kern-\ex@\nointerlineskip}$\m@th\hfil#1#2\hfil$\crcr}}}

\def\overleftarrow{\mathpalette\overleftarrow@}
\def\overleftarrow@#1#2{\vbox{\ialign{##\crcr\leftarrowfill@#1\crcr
 \noalign{\kern-\ex@\nointerlineskip}$\m@th\hfil#1#2\hfil$\crcr}}}
\def\overleftrightarrow{\mathpalette\overleftrightarrow@}
\def\overleftrightarrow@#1#2{\vbox{\ialign{##\crcr\leftrightarrowfill@#1\crcr
 \noalign{\kern-\ex@\nointerlineskip}$\m@th\hfil#1#2\hfil$\crcr}}}
\def\underrightarrow{\mathpalette\underrightarrow@}
\def\underrightarrow@#1#2{\vtop{\ialign{##\crcr$\m@th\hfil#1#2\hfil$\crcr
 \noalign{\nointerlineskip}\rightarrowfill@#1\crcr}}}

\def\underleftarrow{\mathpalette\underleftarrow@}
\def\underleftarrow@#1#2{\vtop{\ialign{##\crcr$\m@th\hfil#1#2\hfil$\crcr
 \noalign{\nointerlineskip}\leftarrowfill@#1\crcr}}}
\def\underleftrightarrow{\mathpalette\underleftrightarrow@}
\def\underleftrightarrow@#1#2{\vtop{\ialign{##\crcr$\m@th\hfil#1#2\hfil$\crcr
 \noalign{\nointerlineskip}\leftrightarrowfill@#1\crcr}}}
\message{various kinds of dots,}
\let\DOTSI\relax
\let\DOTSB\relax

\newif\ifmath@
{\uccode`7=`\\ \uccode`8=`m \uccode`9=`a \uccode`0=`t \uccode`!=`h
 \uppercase{\gdef\math@#1#2#3#4#5#6\math@{\global\math@false\ifx 7#1\ifx 8#2%
 \ifx 9#3\ifx 0#4\ifx !#5\xdef\meaning@{#6}\global\math@true\fi\fi\fi\fi\fi}}}
\newif\ifmathch@
{\uccode`7=`c \uccode`8=`h \uccode`9=`\"
 \uppercase{\gdef\mathch@#1#2#3#4#5#6\mathch@{\global\mathch@false
  \ifx 7#1\ifx 8#2\ifx 9#5\global\mathch@true\xdef\meaning@{9#6}\fi\fi\fi}}}
\newcount\classnum@
\def\getmathch@#1.#2\getmathch@{\classnum@#1 \divide\classnum@4096
 \ifcase\number\classnum@\or\or\gdef\thedots@{\dotsb@}\or
 \gdef\thedots@{\dotsb@}\fi}
\newif\ifmathbin@
{\uccode`4=`b \uccode`5=`i \uccode`6=`n
 \uppercase{\gdef\mathbin@#1#2#3{\relaxnext@
  \DNii@##1\mathbin@{\ifx\space@\next\global\mathbin@true\fi}%
 \global\mathbin@false\DN@##1\mathbin@{}%
 \ifx 4#1\ifx 5#2\ifx 6#3\DN@{\FN@\nextii@}\fi\fi\fi\next@}}}
\newif\ifmathrel@
{\uccode`4=`r \uccode`5=`e \uccode`6=`l
 \uppercase{\gdef\mathrel@#1#2#3{\relaxnext@
  \DNii@##1\mathrel@{\ifx\space@\next\global\mathrel@true\fi}%
 \global\mathrel@false\DN@##1\mathrel@{}%
 \ifx 4#1\ifx 5#2\ifx 6#3\DN@{\FN@\nextii@}\fi\fi\fi\next@}}}
\newif\ifmacro@
{\uccode`5=`m \uccode`6=`a \uccode`7=`c
 \uppercase{\gdef\macro@#1#2#3#4\macro@{\global\macro@false
  \ifx 5#1\ifx 6#2\ifx 7#3\global\macro@true
  \xdef\meaning@{\macro@@#4\macro@@}\fi\fi\fi}}}
\def\macro@@#1->#2\macro@@{#2}
\newif\ifDOTS@
\newcount\DOTSCASE@
{\uccode`6=`\\ \uccode`7=`D \uccode`8=`O \uccode`9=`T \uccode`0=`S
 \uppercase{\gdef\DOTS@#1#2#3#4#5{\global\DOTS@false\DN@##1\DOTS@{}%
  \ifx 6#1\ifx 7#2\ifx 8#3\ifx 9#4\ifx 0#5\let\next@\DOTS@@\fi\fi\fi\fi\fi
  \next@}}}
{\uccode`3=`B \uccode`4=`I \uccode`5=`X
 \uppercase{\gdef\DOTS@@#1{\relaxnext@
  \DNii@##1\DOTS@{\ifx\space@\next\global\DOTS@true\fi}%
  \DN@{\FN@\nextii@}%
  \ifx 3#1\global\DOTSCASE@\z@\else
  \ifx 4#1\global\DOTSCASE@\@ne\else
  \ifx 5#1\global\DOTSCASE@\tw@\else\DN@##1\DOTS@{}%
  \fi\fi\fi\next@}}}
\newif\ifnot@
{\uccode`5=`\\ \uccode`6=`n \uccode`7=`o \uccode`8=`t
 \uppercase{\gdef\not@#1#2#3#4{\relaxnext@
  \DNii@##1\not@{\ifx\space@\next\global\not@true\fi}%
 \global\not@false\DN@##1\not@{}%
 \ifx 5#1\ifx 6#2\ifx 7#3\ifx 8#4\DN@{\FN@\nextii@}\fi\fi\fi
 \fi\next@}}}
\newif\ifkeybin@
\def\keybin@{\keybin@true
 \ifx\next+\else\ifx\next=\else\ifx\next<\else\ifx\next>\else\ifx\next-\else
 \ifx\next*\else\ifx\next:\else\keybin@false\fi\fi\fi\fi\fi\fi\fi}
\def\dots{\RIfM@\expandafter\mdots@\else\expandafter\tdots@\fi}
\def\tdots@{\unskip\relaxnext@
 \DN@{$\m@th\mathinner{\ldotp\ldotp\ldotp}\,
   \ifx\next,\,$\else\ifx\next.\,$\else\ifx\next;\,$\else\ifx\next:\,$\else
   \ifx\next?\,$\else\ifx\next!\,$\else$ \fi\fi\fi\fi\fi\fi}%
 \ \FN@\next@}
\def\mdots@{\FN@\mdots@@}
\def\mdots@@{\gdef\thedots@{\dotso@}
 \ifx\next\boldkey\gdef\thedots@\boldkey{\boldkeydots@}\else                
 \ifx\next\boldsymbol\gdef\thedots@\boldsymbol{\boldsymboldots@}\else       
 \ifx,\next\gdef\thedots@{\dotsc}
 \else\ifx\not\next\gdef\thedots@{\dotsb@}
 \else\keybin@
 \ifkeybin@\gdef\thedots@{\dotsb@}
 \else\xdef\meaning@{\meaning\next..........}\xdef\meaning@@{\meaning@}
  \expandafter\math@\meaning@\math@
  \ifmath@
   \expandafter\mathch@\meaning@\mathch@
   \ifmathch@\expandafter\getmathch@\meaning@\getmathch@\fi                 
  \else\expandafter\macro@\meaning@@\macro@                                 
  \ifmacro@                                                                
   \expandafter\not@\meaning@\not@\ifnot@\gdef\thedots@{\dotsb@}
  \else\expandafter\DOTS@\meaning@\DOTS@
  \ifDOTS@
   \ifcase\number\DOTSCASE@\gdef\thedots@{\dotsb@}%
    \or\gdef\thedots@{\dotsi}\else\fi                                      
  \else\expandafter\math@\meaning@\math@                                   
  \ifmath@\expandafter\mathbin@\meaning@\mathbin@
  \ifmathbin@\gdef\thedots@{\dotsb@}
  \else\expandafter\mathrel@\meaning@\mathrel@
  \ifmathrel@\gdef\thedots@{\dotsb@}
  \fi\fi\fi\fi\fi\fi\fi\fi\fi\fi\fi\fi
 \thedots@}
\def\plainldots@{\mathinner{\ldotp\ldotp\ldotp}}
\def\plaincdots@{\mathinner{\cdotp\cdotp\cdotp}}
\def\dotsi{\!\plaincdots@}
\let\dotsb@\plaincdots@
\newif\ifextra@
\newif\ifrightdelim@
\def\rightdelim@{\global\rightdelim@true                                    
 \ifx\next)\else                                                            
 \ifx\next]\else
 \ifx\next\rbrack\else
 \ifx\next\}\else
 \ifx\next\rbrace\else
 \ifx\next\rangle\else
 \ifx\next\rceil\else
 \ifx\next\rfloor\else
 \ifx\next\rgroup\else
 \ifx\next\rmoustache\else
 \ifx\next\right\else
 \ifx\next\bigr\else
 \ifx\next\biggr\else
 \ifx\next\Bigr\else                                                        
 \ifx\next\Biggr\else\global\rightdelim@false
 \fi\fi\fi\fi\fi\fi\fi\fi\fi\fi\fi\fi\fi\fi\fi}
\def\extra@{%
 \global\extra@false\rightdelim@\ifrightdelim@\global\extra@true            
 \else\ifx\next$\global\extra@true                                          
 \else\xdef\meaning@{\meaning\next..........}
 \expandafter\macro@\meaning@\macro@\ifmacro@                               
 \expandafter\DOTS@\meaning@\DOTS@
 \ifDOTS@
 \ifnum\DOTSCASE@=\tw@\global\extra@true                                    
 \fi\fi\fi\fi\fi}
\newif\ifbold@
\def\dotso@{\relaxnext@
 \ifbold@
  \let\next\delayed@
  \DNii@{\extra@\plainldots@\ifextra@\,\fi}%
 \else
  \DNii@{\DN@{\extra@\plainldots@\ifextra@\,\fi}\FN@\next@}%
 \fi
 \nextii@}
\def\extrap@#1{%
 \ifx\next,\DN@{#1\,}\else
 \ifx\next;\DN@{#1\,}\else
 \ifx\next.\DN@{#1\,}\else\extra@
 \ifextra@\DN@{#1\,}\else
 \let\next@#1\fi\fi\fi\fi\next@}
\def\ldots{\DN@{\extrap@\plainldots@}%
 \FN@\next@}
\def\cdots{\DN@{\extrap@\plaincdots@}%
 \FN@\next@}

\def\dotsc{\relaxnext@
 \DN@{\ifx\next;\plainldots@\,\else
  \ifx\next.\plainldots@\,\else\extra@\plainldots@
  \ifextra@\,\fi\fi\fi}%
 \FN@\next@}
\def\cdot{\mathchar"2201 }

\def\mapsto{\DOTSB\mapstochar\rightarrow}

\message{special superscripts,}
\def\dddot#1{{\mathop{#1}\limits^{\vbox to-1.4\ex@{\kern-\tw@\ex@
 \hbox{\rm...}\vss}}}}
\def\ddddot#1{{\mathop{#1}\limits^{\vbox to-1.4\ex@{\kern-\tw@\ex@
 \hbox{\rm....}\vss}}}}
\def\sphat{^{\mathchoice{}{}%
 {\,\,\botsmash{\hbox{\lower4\ex@\hbox{$\m@th\widehat{\null}$}}}}%
 {\,\botsmash{\hbox{\lower3\ex@\hbox{$\m@th\hat{\null}$}}}}}}

\def\spacute{^{\!\botsmash{\hbox{\lower\@ne ex\hbox{\'{}}}}}}
\def\spgrave{^{\mathchoice{}{}{}{\!}%
 \botsmash{\hbox{\lower\@ne ex\hbox{\`{}}}}}}
\def\spdot{^{\hbox{\raise\ex@\hbox{\rm.}}}}
\def\spddot{^{\hbox{\raise\ex@\hbox{\rm..}}}}
\def\spdddot{^{\hbox{\raise\ex@\hbox{\rm...}}}}
\def\spddddot{^{\hbox{\raise\ex@\hbox{\rm....}}}}
\def\spbreve{^{\!\botsmash{\hbox{\lower4\ex@\hbox{\u{}}}}}}

\message{\string\text,}
\def\textonlyfont@#1#2{\def#1{\RIfM@
 \Err@{Use \string#1\space only in text}\else#2\fi}}
\textonlyfont@\rm\tenrm
\textonlyfont@\it\tenit
\textonlyfont@\sl\tensl
\textonlyfont@\bf\tenbf
\def\oldnos#1{\RIfM@{\mathcode`\,="013B \fam\@ne#1}\else
 \leavevmode\hbox{$\m@th\mathcode`\,="013B \fam\@ne#1$}\fi}
\def\text{\RIfM@\expandafter\text@\else\expandafter\text@@\fi}
\def\text@@#1{\leavevmode\hbox{#1}}
\def\mathhexbox@#1#2#3{\text{$\m@th\mathchar"#1#2#3$}}
\def\dag{{\mathhexbox@279}}
\def\ddag{{\mathhexbox@27A}}
\def\S{{\mathhexbox@278}}
\def\P{{\mathhexbox@27B}}
\newif\iffirstchoice@
\firstchoice@true
\def\text@#1{\mathchoice
 {\hbox{\everymath{\displaystyle}\def\textfonti{\the\textfont\@ne}%
  \def\textfontii{\the\textfont\tw@}\textdef@@ T#1}}
 {\hbox{\firstchoice@false
  \everymath{\textstyle}\def\textfonti{\the\textfont\@ne}%
  \def\textfontii{\the\textfont\tw@}\textdef@@ T#1}}
 {\hbox{\firstchoice@false
  \everymath{\scriptstyle}\def\textfonti{\the\scriptfont\@ne}%
  \def\textfontii{\the\scriptfont\tw@}\textdef@@ S\rm#1}}
 {\hbox{\firstchoice@false
  \everymath{\scriptscriptstyle}\def\textfonti
  {\the\scriptscriptfont\@ne}%
  \def\textfontii{\the\scriptscriptfont\tw@}\textdef@@ s\rm#1}}}
\def\textdef@@#1{\textdef@#1\rm\textdef@#1\bf\textdef@#1\sl\textdef@#1\it}
\def\rmfam{0}
\def\textdef@#1#2{%
 \DN@{\csname\expandafter\eat@\string#2fam\endcsname}%
 \if S#1\edef#2{\the\scriptfont\next@\relax}%
 \else\if s#1\edef#2{\the\scriptscriptfont\next@\relax}%
 \else\edef#2{\the\textfont\next@\relax}\fi\fi}
\scriptfont\itfam\tenit \scriptscriptfont\itfam\tenit
\scriptfont\slfam\tensl \scriptscriptfont\slfam\tensl
\newif\iftopfolded@
\newif\ifbotfolded@
\def\topfoldedtext{\topfolded@true\botfolded@false\foldedtext@}
\def\botfoldedtext{\botfolded@true\topfolded@false\foldedtext@}
\def\foldedtext{\topfolded@false\botfolded@false\foldedtext@}
\Invalid@\foldedwidth
\def\foldedtext@{\relaxnext@
 \DN@{\ifx\next\foldedwidth\let\next@\nextii@\else
  \DN@{\nextii@\foldedwidth{.3\hsize}}\fi\next@}%
 \DNii@\foldedwidth##1##2{\setbox\z@\vbox
  {\normalbaselines\hsize##1\relax
  \tolerance1600 \noindent\ignorespaces##2}\ifbotfolded@\boxz@\else
  \iftopfolded@\vtop{\unvbox\z@}\else\vcenter{\boxz@}\fi\fi}%
 \FN@\next@}
\message{math font commands,}
\def\bold{\RIfM@\expandafter\bold@\else
 \expandafter\nonmatherr@\expandafter\bold\fi}
\def\bold@#1{{\bold@@{#1}}}
\def\bold@@#1{\fam\bffam\relax#1}
\def\slanted{\RIfM@\expandafter\slanted@\else
 \expandafter\nonmatherr@\expandafter\slanted\fi}
\def\slanted@#1{{\slanted@@{#1}}}
\def\slanted@@#1{\fam\slfam\relax#1}
\def\roman{\RIfM@\expandafter\roman@\else
 \expandafter\nonmatherr@\expandafter\roman\fi}
\def\roman@#1{{\roman@@{#1}}}
\def\roman@@#1{\fam\rmfam\relax#1}
\def\italic{\RIfM@\expandafter\italic@\else
 \expandafter\nonmatherr@\expandafter\italic\fi}
\def\italic@#1{{\italic@@{#1}}}
\def\italic@@#1{\fam\itfam\relax#1}
\def\Cal{\RIfM@\expandafter\Cal@\else
 \expandafter\nonmatherr@\expandafter\Cal\fi}
\def\Cal@#1{{\Cal@@{#1}}}
\def\Cal@@#1{\noaccents@\fam\tw@#1}
\mathchardef\Gamma="0000
\mathchardef\Delta="0001
\mathchardef\Theta="0002
\mathchardef\Lambda="0003
\mathchardef\Xi="0004
\mathchardef\Pi="0005
\mathchardef\Sigma="0006
\mathchardef\Upsilon="0007
\mathchardef\Phi="0008
\mathchardef\Psi="0009
\mathchardef\Omega="000A
\mathchardef\varGamma="0100
\mathchardef\varDelta="0101
\mathchardef\varTheta="0102
\mathchardef\varLambda="0103
\mathchardef\varXi="0104
\mathchardef\varPi="0105
\mathchardef\varSigma="0106
\mathchardef\varUpsilon="0107
\mathchardef\varPhi="0108
\mathchardef\varPsi="0109
\mathchardef\varOmega="010A
\let\alloc@@\alloc@
\def\hexnumber@#1{\ifcase#1 0\or 1\or 2\or 3\or 4\or 5\or 6\or 7\or 8\or
 9\or A\or B\or C\or D\or E\or F\fi}
\def\loadmsam{%
 \font@\tenmsa=msam10
 \font@\sevenmsa=msam7
 \font@\fivemsa=msam5
 \alloc@@8\fam\chardef\sixt@@n\msafam
 \textfont\msafam=\tenmsa
 \scriptfont\msafam=\sevenmsa
 \scriptscriptfont\msafam=\fivemsa
 \edef\next{\hexnumber@\msafam}%
 \mathchardef\dabar@"0\next39
 \edef\dashrightarrow{\mathrel{\dabar@\dabar@\mathchar"0\next4B}}%
 \edef\dashleftarrow{\mathrel{\mathchar"0\next4C\dabar@\dabar@}}%
 \let\dasharrow\dashrightarrow
 \edef\ulcorner{\delimiter"4\next70\next70 }%
 \edef\urcorner{\delimiter"5\next71\next71 }%
 \edef\llcorner{\delimiter"4\next78\next78 }%
 \edef\lrcorner{\delimiter"5\next79\next79 }%
 \edef\yen{{\noexpand\mathhexbox@\next55}}%
 \edef\checkmark{{\noexpand\mathhexbox@\next58}}%
 \edef\circledR{{\noexpand\mathhexbox@\next72}}%
 \edef\maltese{{\noexpand\mathhexbox@\next7A}}%
 \global\let\loadmsam\empty}%
\def\loadmsbm{%
 \font@\tenmsb=msbm10 \font@\sevenmsb=msbm7 \font@\fivemsb=msbm5
 \alloc@@8\fam\chardef\sixt@@n\msbfam
 \textfont\msbfam=\tenmsb
 \scriptfont\msbfam=\sevenmsb \scriptscriptfont\msbfam=\fivemsb
 \global\let\loadmsbm\empty
 }
\def\widehat#1{\ifx\undefined\msbfam \DN@{362}%
  \else \setboxz@h{$\m@th#1$}%
    \edef\next@{\ifdim\wdz@>\tw@ em%
        \hexnumber@\msbfam 5B%
      \else 362\fi}\fi
  \mathaccent"0\next@{#1}}
\def\widetilde#1{\ifx\undefined\msbfam \DN@{365}%
  \else \setboxz@h{$\m@th#1$}%
    \edef\next@{\ifdim\wdz@>\tw@ em%
        \hexnumber@\msbfam 5D%
      \else 365\fi}\fi
  \mathaccent"0\next@{#1}}
\message{\string\newsymbol,}
\def\newsymbol#1#2#3#4#5{\define#1{}%
  \count@#2\relax \advance\count@\m@ne 
 \ifcase\count@
   \ifx\undefined\msafam\loadmsam\fi \let\next@\msafam
 \or \ifx\undefined\msbfam\loadmsbm\fi \let\next@\msbfam
 \else  \Err@{\Invalid@@\string\newsymbol}\let\next@\tw@\fi
 \mathchardef#1="#3\hexnumber@\next@#4#5\space}

\def\Bbb{\RIfM@\expandafter\Bbb@\else
 \expandafter\nonmatherr@\expandafter\Bbb\fi}
\def\Bbb@#1{{\Bbb@@{#1}}}
\def\Bbb@@#1{\noaccents@\fam\msbfam\relax#1}
\message{bold Greek and bold symbols,}
\def\loadbold{%
 \font@\tencmmib=cmmib10 \font@\sevencmmib=cmmib7 \font@\fivecmmib=cmmib5
 \skewchar\tencmmib'177 \skewchar\sevencmmib'177 \skewchar\fivecmmib'177
 \alloc@@8\fam\chardef\sixt@@n\cmmibfam
 \textfont\cmmibfam\tencmmib
 \scriptfont\cmmibfam\sevencmmib \scriptscriptfont\cmmibfam\fivecmmib
 \font@\tencmbsy=cmbsy10 \font@\sevencmbsy=cmbsy7 \font@\fivecmbsy=cmbsy5
 \skewchar\tencmbsy'60 \skewchar\sevencmbsy'60 \skewchar\fivecmbsy'60
 \alloc@@8\fam\chardef\sixt@@n\cmbsyfam
 \textfont\cmbsyfam\tencmbsy
 \scriptfont\cmbsyfam\sevencmbsy \scriptscriptfont\cmbsyfam\fivecmbsy
 \let\loadbold\empty
}
\def\boldnotloaded#1{\Err@{\ifcase#1\or First\else Second\fi
       bold symbol font not loaded}}
\def\mathchari@#1#2#3{\ifx\undefined\cmmibfam
    \boldnotloaded@\@ne
  \else\mathchar"#1\hexnumber@\cmmibfam#2#3\space \fi}
\def\mathcharii@#1#2#3{\ifx\undefined\cmbsyfam
    \boldnotloaded\tw@
  \else \mathchar"#1\hexnumber@\cmbsyfam#2#3\space\fi}
\edef\bffam@{\hexnumber@\bffam}
\def\boldkey#1{\ifcat\noexpand#1A%
  \ifx\undefined\cmmibfam \boldnotloaded\@ne
  \else {\fam\cmmibfam#1}\fi
 \else
 \ifx#1!\mathchar"5\bffam@21 \else
 \ifx#1(\mathchar"4\bffam@28 \else\ifx#1)\mathchar"5\bffam@29 \else
 \ifx#1+\mathchar"2\bffam@2B \else\ifx#1:\mathchar"3\bffam@3A \else
 \ifx#1;\mathchar"6\bffam@3B \else\ifx#1=\mathchar"3\bffam@3D \else
 \ifx#1?\mathchar"5\bffam@3F \else\ifx#1[\mathchar"4\bffam@5B \else
 \ifx#1]\mathchar"5\bffam@5D \else
 \ifx#1,\mathchari@63B \else
 \ifx#1-\mathcharii@200 \else
 \ifx#1.\mathchari@03A \else
 \ifx#1/\mathchari@03D \else
 \ifx#1<\mathchari@33C \else
 \ifx#1>\mathchari@33E \else
 \ifx#1*\mathcharii@203 \else
 \ifx#1|\mathcharii@06A \else
 \ifx#10\bold0\else\ifx#11\bold1\else\ifx#12\bold2\else\ifx#13\bold3\else
 \ifx#14\bold4\else\ifx#15\bold5\else\ifx#16\bold6\else\ifx#17\bold7\else
 \ifx#18\bold8\else\ifx#19\bold9\else
  \Err@{\string\boldkey\space can't be used with #1}%
 \fi\fi\fi\fi\fi\fi\fi\fi\fi\fi\fi\fi\fi\fi\fi
 \fi\fi\fi\fi\fi\fi\fi\fi\fi\fi\fi\fi\fi\fi}
\def\boldsymbol#1{%
 \DN@{\Err@{You can't use \string\boldsymbol\space with \string#1}#1}%
 \ifcat\noexpand#1A%
   \let\next@\relax
   \ifx\undefined\cmmibfam \boldnotloaded\@ne
   \else {\fam\cmmibfam#1}\fi
 \else
  \xdef\meaning@{\meaning#1.........}%
  \expandafter\math@\meaning@\math@
  \ifmath@
   \expandafter\mathch@\meaning@\mathch@
   \ifmathch@
    \expandafter\boldsymbol@@\meaning@\boldsymbol@@
   \fi
  \else
   \expandafter\macro@\meaning@\macro@
   \expandafter\delim@\meaning@\delim@
   \ifdelim@
    \expandafter\delim@@\meaning@\delim@@
   \else
    \boldsymbol@{#1}%
   \fi
  \fi
 \fi
 \next@}
\def\mathhexboxii@#1#2{\ifx\undefined\cmbsyfam
    \boldnotloaded\tw@
  \else \mathhexbox@{\hexnumber@\cmbsyfam}{#1}{#2}\fi}
\def\boldsymbol@#1{\let\next@\relax\let\next#1%
 \ifx\next\cdot\mathcharii@201 \else
 \ifx\next\prime{{\null\mathcharii@030 \null}}\else
 \ifx\next\lbrack\mathchar"4\bffam@5B \else
 \ifx\next\rbrack\mathchar"5\bffam@5D \else
 \ifx\next\{\mathcharii@466 \else
 \ifx\next\lbrace\mathcharii@466 \else
 \ifx\next\}\mathcharii@567 \else
 \ifx\next\rbrace\mathcharii@567 \else
 \ifx\next\surd{{\mathcharii@170}}\else
 \ifx\next\S{{\mathhexboxii@78}}\else
 \ifx\next\P{{\mathhexboxii@7B}}\else
 \ifx\next\dag{{\mathhexboxii@79}}\else
 \ifx\next\ddag{{\mathhexboxii@7A}}\else
 \DN@{\Err@{You can't use \string\boldsymbol\space with \string#1}#1}%
 \fi\fi\fi\fi\fi\fi\fi\fi\fi\fi\fi\fi\fi}
\def\boldsymbol@@#1.#2\boldsymbol@@{\classnum@#1 \count@@@\classnum@        
 \divide\classnum@4096 \count@\classnum@                                    
 \multiply\count@4096 \advance\count@@@-\count@ \count@@\count@@@           
 \divide\count@@@\@cclvi \count@\count@@                                    
 \multiply\count@@@\@cclvi \advance\count@@-\count@@@                       
 \divide\count@@@\@cclvi                                                    
 \multiply\classnum@4096 \advance\classnum@\count@@                         
 \ifnum\count@@@=\z@                                                        
  \count@"\bffam@ \multiply\count@\@cclvi
  \advance\classnum@\count@
  \DN@{\mathchar\number\classnum@}%
 \else
  \ifnum\count@@@=\@ne                                                      
   \ifx\undefined\cmmibfam \DN@{\boldnotloaded\@ne}%
   \else \count@\cmmibfam \multiply\count@\@cclvi
     \advance\classnum@\count@
     \DN@{\mathchar\number\classnum@}\fi
  \else
   \ifnum\count@@@=\tw@                                                    
     \ifx\undefined\cmbsyfam
       \DN@{\boldnotloaded\tw@}%
     \else
       \count@\cmbsyfam \multiply\count@\@cclvi
       \advance\classnum@\count@
       \DN@{\mathchar\number\classnum@}%
     \fi
  \fi
 \fi
\fi}
\newif\ifdelim@
\newcount\delimcount@
{\uccode`6=`\\ \uccode`7=`d \uccode`8=`e \uccode`9=`l
 \uppercase{\gdef\delim@#1#2#3#4#5\delim@
  {\delim@false\ifx 6#1\ifx 7#2\ifx 8#3\ifx 9#4\delim@true
   \xdef\meaning@{#5}\fi\fi\fi\fi}}}
\def\delim@@#1"#2#3#4#5#6\delim@@{\if#32%
\let\next@\relax
 \ifx\undefined\cmbsyfam \boldnotloaded\@ne
 \else \mathcharii@#2#4#5\space \fi\fi}
\def\vert{\delimiter"026A30C }
\def\Vert{\delimiter"026B30D }
\let\|\Vert
\def\backslash{\delimiter"026E30F }
\def\boldkeydots@#1{\bold@true\let\next=#1\let\delayed@=#1\mdots@@
 \boldkey#1\bold@false}  
\def\boldsymboldots@#1{\bold@true\let\next#1\let\delayed@#1\mdots@@
 \boldsymbol#1\bold@false}
\message{Euler fonts,}

\def\frak{\mathfont@\frak}

\def\loadmathfont#1{%
   \expandafter\font@\csname ten#1\endcsname=#110
   \expandafter\font@\csname seven#1\endcsname=#17
   \expandafter\font@\csname five#1\endcsname=#15
   \edef\next{\noexpand\alloc@@8\fam\chardef\sixt@@n
     \expandafter\noexpand\csname#1fam\endcsname}%
   \next
   \textfont\csname#1fam\endcsname \csname ten#1\endcsname
   \scriptfont\csname#1fam\endcsname \csname seven#1\endcsname
   \scriptscriptfont\csname#1fam\endcsname \csname five#1\endcsname
   \expandafter\def\csname #1\expandafter\endcsname\expandafter{%
      \expandafter\mathfont@\csname#1\endcsname}%
 \expandafter\gdef\csname load#1\endcsname{}%
}
\def\mathfont@#1{\RIfM@\expandafter\mathfont@@\expandafter#1\else
  \expandafter\nonmatherr@\expandafter#1\fi}
\def\mathfont@@#1#2{{\mathfont@@@#1{#2}}}
\def\mathfont@@@#1#2{\noaccents@
   \fam\csname\expandafter\eat@\string#1fam\endcsname
   \relax#2}
\message{math accents,}
\def\accentclass@{7}
\def\noaccents@{\def\accentclass@{0}}
\def\makeacc@#1#2{\def#1{\mathaccent"\accentclass@#2 }}
\makeacc@\hat{05E}
\makeacc@\check{014}
\makeacc@\tilde{07E}
\makeacc@\acute{013}
\makeacc@\grave{012}
\makeacc@\dot{05F}
\makeacc@\ddot{07F}
\makeacc@\breve{015}
\makeacc@\bar{016}

\newcount\skewcharcount@
\newcount\familycount@
\def\theskewchar@{\familycount@\@ne
 \global\skewcharcount@\the\skewchar\textfont\@ne                           
 \ifnum\fam>\m@ne\ifnum\fam<16
  \global\familycount@\the\fam\relax
  \global\skewcharcount@\the\skewchar\textfont\the\fam\relax\fi\fi          
 \ifnum\skewcharcount@>\m@ne
  \ifnum\skewcharcount@<128
  \multiply\familycount@256
  \global\advance\skewcharcount@\familycount@
  \global\advance\skewcharcount@28672
  \mathchar\skewcharcount@\else
  \global\skewcharcount@\m@ne\fi\else
 \global\skewcharcount@\m@ne\fi}                                            
\newcount\pointcount@
\def\getpoints@#1.#2\getpoints@{\pointcount@#1 }
\newdimen\accentdimen@
\newcount\accentmu@
\def\dimentomu@{\multiply\accentdimen@ 100
 \expandafter\getpoints@\the\accentdimen@\getpoints@
 \multiply\pointcount@18
 \divide\pointcount@\@m
 \global\accentmu@\pointcount@}
\def\Makeacc@#1#2{\def#1{\RIfM@\DN@{\mathaccent@
 {"\accentclass@#2 }}\else\DN@{\nonmatherr@{#1}}\fi\next@}}
\def\unbracefonts@{\let\Cal@\Cal@@\let\roman@\roman@@\let\bold@\bold@@
 \let\slanted@\slanted@@}
\def\mathaccent@#1#2{\ifnum\fam=\m@ne\xdef\thefam@{1}\else
 \xdef\thefam@{\the\fam}\fi                                                 
 \accentdimen@\z@                                                           
 \setboxz@h{\unbracefonts@$\m@th\fam\thefam@\relax#2$}
 \ifdim\accentdimen@=\z@\DN@{\mathaccent#1{#2}}
  \setbox@ne\hbox{\unbracefonts@$\m@th\fam\thefam@\relax#2\theskewchar@$}
  \setbox\tw@\hbox{$\m@th\ifnum\skewcharcount@=\m@ne\else
   \mathchar\skewcharcount@\fi$}
  \global\accentdimen@\wd@ne\global\advance\accentdimen@-\wdz@
  \global\advance\accentdimen@-\wd\tw@                                     
  \global\multiply\accentdimen@\tw@
  \dimentomu@\global\advance\accentmu@\@ne                                 
 \else\DN@{{\mathaccent#1{#2\mkern\accentmu@ mu}%
    \mkern-\accentmu@ mu}{}}\fi                                             
 \next@}\Makeacc@\Hat{05E}
\Makeacc@\Check{014}
\Makeacc@\Tilde{07E}
\Makeacc@\Acute{013}
\Makeacc@\Grave{012}
\Makeacc@\Dot{05F}
\Makeacc@\Ddot{07F}
\Makeacc@\Breve{015}
\Makeacc@\Bar{016}
\def\Vec{\RIfM@\DN@{\mathaccent@{"017E }}\else
 \DN@{\nonmatherr@\Vec}\fi\next@}
\def\accentedsymbol#1#2{\csname newbox\expandafter\endcsname
  \csname\expandafter\eat@\string#1@box\endcsname
 \expandafter\setbox\csname\expandafter\eat@
  \string#1@box\endcsname\hbox{$\m@th#2$}\define
  #1{\copy\csname\expandafter\eat@\string#1@box\endcsname{}}}
\message{roots,}
\def\sqrt#1{\radical"270370 {#1}}
\let\underline@\underline
\let\overline@\overline
\def\underline#1{\underline@{#1}}
\def\overline#1{\overline@{#1}}
\Invalid@\leftroot
\Invalid@\uproot
\newcount\uproot@
\newcount\leftroot@
\def\root{\relaxnext@
  \DN@{\ifx\next\uproot\let\next@\nextii@\else
   \ifx\next\leftroot\let\next@\nextiii@\else
   \let\next@\plainroot@\fi\fi\next@}%
  \DNii@\uproot##1{\uproot@##1\relax\FN@\nextiv@}%
  \def\nextiv@{\ifx\next\space@\DN@. {\FN@\nextv@}\else
   \DN@.{\FN@\nextv@}\fi\next@.}%
  \def\nextv@{\ifx\next\leftroot\let\next@\nextvi@\else
   \let\next@\plainroot@\fi\next@}%
  \def\nextvi@\leftroot##1{\leftroot@##1\relax\plainroot@}%
   \def\nextiii@\leftroot##1{\leftroot@##1\relax\FN@\nextvii@}%
  \def\nextvii@{\ifx\next\space@
   \DN@. {\FN@\nextviii@}\else
   \DN@.{\FN@\nextviii@}\fi\next@.}%
  \def\nextviii@{\ifx\next\uproot\let\next@\nextix@\else
   \let\next@\plainroot@\fi\next@}%
  \def\nextix@\uproot##1{\uproot@##1\relax\plainroot@}%
  \bgroup\uproot@\z@\leftroot@\z@\FN@\next@}
\def\plainroot@#1\of#2{\setbox\rootbox\hbox{$\m@th\scriptscriptstyle{#1}$}%
 \mathchoice{\r@@t\displaystyle{#2}}{\r@@t\textstyle{#2}}
 {\r@@t\scriptstyle{#2}}{\r@@t\scriptscriptstyle{#2}}\egroup}
\def\r@@t#1#2{\setboxz@h{$\m@th#1\sqrt{#2}$}%
 \dimen@\ht\z@\advance\dimen@-\dp\z@
 \setbox@ne\hbox{$\m@th#1\mskip\uproot@ mu$}\advance\dimen@ 1.667\wd@ne
 \mkern-\leftroot@ mu\mkern5mu\raise.6\dimen@\copy\rootbox
 \mkern-10mu\mkern\leftroot@ mu\boxz@}
\def\boxed#1{\setboxz@h{$\m@th\displaystyle{#1}$}\dimen@.4\ex@
 \advance\dimen@3\ex@\advance\dimen@\dp\z@
 \hbox{\lower\dimen@\hbox{%
 \vbox{\hrule height.4\ex@
 \hbox{\vrule width.4\ex@\hskip3\ex@\vbox{\vskip3\ex@\boxz@\vskip3\ex@}%
 \hskip3\ex@\vrule width.4\ex@}\hrule height.4\ex@}%
 }}}
\message{commutative diagrams,}
\let\ampersand@\relax
\newdimen\minaw@
\minaw@11.11128\ex@
\newdimen\minCDaw@
\minCDaw@2.5pc
\def\minCDarrowwidth#1{\RIfMIfI@\onlydmatherr@\minCDarrowwidth
 \else\minCDaw@#1\relax\fi\else\onlydmatherr@\minCDarrowwidth\fi}
\newif\ifCD@
\def\CD{\bgroup\vspace@\relax\let\ampersand@&\iffalse}\fi
 \CD@true\vcenter\bgroup\Let@\tabskip\z@skip\baselineskip20\ex@
 \lineskip3\ex@\lineskiplimit3\ex@\halign\bgroup
 &\hfill$\m@th##$\hfill\crcr}
\def\endCD{\crcr\egroup\egroup\egroup}
\newdimen\bigaw@
\atdef@>#1>#2>{\ampersand@                                                  
 \setboxz@h{$\m@th\ssize\;{#1}\;\;$}
 \setbox@ne\hbox{$\m@th\ssize\;{#2}\;\;$}
 \setbox\tw@\hbox{$\m@th#2$}
 \ifCD@\global\bigaw@\minCDaw@\else\global\bigaw@\minaw@\fi                 
 \ifdim\wdz@>\bigaw@\global\bigaw@\wdz@\fi
 \ifdim\wd@ne>\bigaw@\global\bigaw@\wd@ne\fi                                
 \ifCD@\enskip\fi                                                           
 \ifdim\wd\tw@>\z@
  \mathrel{\mathop{\hbox to\bigaw@{\rightarrowfill@\displaystyle}}%
    \limits^{#1}_{#2}}
 \else\mathrel{\mathop{\hbox to\bigaw@{\rightarrowfill@\displaystyle}}%
    \limits^{#1}}\fi                                                        
 \ifCD@\enskip\fi                                                          
 \ampersand@}                                                              
\atdef@<#1<#2<{\ampersand@\setboxz@h{$\m@th\ssize\;\;{#1}\;$}%
 \setbox@ne\hbox{$\m@th\ssize\;\;{#2}\;$}\setbox\tw@\hbox{$\m@th#2$}%
 \ifCD@\global\bigaw@\minCDaw@\else\global\bigaw@\minaw@\fi
 \ifdim\wdz@>\bigaw@\global\bigaw@\wdz@\fi
 \ifdim\wd@ne>\bigaw@\global\bigaw@\wd@ne\fi
 \ifCD@\enskip\fi
 \ifdim\wd\tw@>\z@
  \mathrel{\mathop{\hbox to\bigaw@{\leftarrowfill@\displaystyle}}%
       \limits^{#1}_{#2}}\else
  \mathrel{\mathop{\hbox to\bigaw@{\leftarrowfill@\displaystyle}}%
       \limits^{#1}}\fi
 \ifCD@\enskip\fi\ampersand@}
\begingroup
 \catcode`\~=\active \lccode`\~=`\@
 \lowercase{%
  \global\atdef@)#1)#2){~>#1>#2>}
  \global\atdef@(#1(#2({~<#1<#2<}}
\endgroup
\atdef@ A#1A#2A{\llap{$\m@th\vcenter{\hbox
 {$\ssize#1$}}$}\Big\uparrow\rlap{$\m@th\vcenter{\hbox{$\ssize#2$}}$}&&}
\atdef@ V#1V#2V{\llap{$\m@th\vcenter{\hbox
 {$\ssize#1$}}$}\Big\downarrow\rlap{$\m@th\vcenter{\hbox{$\ssize#2$}}$}&&}
\atdef@={&\enskip\mathrel
 {\vbox{\hrule width\minCDaw@\vskip3\ex@\hrule width
 \minCDaw@}}\enskip&}
\atdef@|{\Big\Vert&&}
\atdef@\vert{\Big\Vert&&}
\def\pretend#1\haswidth#2{\setboxz@h{$\m@th\scriptstyle{#2}$}\hbox
 to\wdz@{\hfill$\m@th\scriptstyle{#1}$\hfill}}
\message{poor man's bold,}
\def\pmb{\RIfM@\expandafter\mathpalette\expandafter\pmb@\else
 \expandafter\pmb@@\fi}
\def\pmb@@#1{\leavevmode\setboxz@h{#1}%
   \dimen@-\wdz@
   \kern-.5\ex@\copy\z@
   \kern\dimen@\kern.25\ex@\raise.4\ex@\copy\z@
   \kern\dimen@\kern.25\ex@\box\z@
}
\def\binrel@@#1{\ifdim\wd2<\z@\mathbin{#1}\else\ifdim\wd\tw@>\z@
 \mathrel{#1}\else{#1}\fi\fi}
\newdimen\pmbraise@
\def\pmb@#1#2{\setbox\thr@@\hbox{$\m@th#1{#2}$}%
 \setbox4\hbox{$\m@th#1\mkern.5mu$}\pmbraise@\wd4\relax
 \binrel@{#2}%
 \dimen@-\wd\thr@@
   \binrel@@{%
   \mkern-.8mu\copy\thr@@
   \kern\dimen@\mkern.4mu\raise\pmbraise@\copy\thr@@
   \kern\dimen@\mkern.4mu\box\thr@@
}}
\def\documentstyle#1{\W@{}\input #1.sty\relax}
\message{syntax check,}
\font\dummyft@=dummy
\fontdimen1 \dummyft@=\z@
\fontdimen2 \dummyft@=\z@
\fontdimen3 \dummyft@=\z@
\fontdimen4 \dummyft@=\z@
\fontdimen5 \dummyft@=\z@
\fontdimen6 \dummyft@=\z@
\fontdimen7 \dummyft@=\z@
\fontdimen8 \dummyft@=\z@
\fontdimen9 \dummyft@=\z@
\fontdimen10 \dummyft@=\z@
\fontdimen11 \dummyft@=\z@
\fontdimen12 \dummyft@=\z@
\fontdimen13 \dummyft@=\z@
\fontdimen14 \dummyft@=\z@
\fontdimen15 \dummyft@=\z@
\fontdimen16 \dummyft@=\z@
\fontdimen17 \dummyft@=\z@
\fontdimen18 \dummyft@=\z@
\fontdimen19 \dummyft@=\z@
\fontdimen20 \dummyft@=\z@
\fontdimen21 \dummyft@=\z@
\fontdimen22 \dummyft@=\z@
\def\fontlist@{\\{\tenrm}\\{\sevenrm}\\{\fiverm}\\{\teni}\\{\seveni}%
 \\{\fivei}\\{\tensy}\\{\sevensy}\\{\fivesy}\\{\tenex}\\{\tenbf}\\{\sevenbf}%
 \\{\fivebf}\\{\tensl}\\{\tenit}}
\def\font@#1=#2 {\rightappend@#1\to\fontlist@\font#1=#2 }
\def\dodummy@{{\def\\##1{\global\let##1\dummyft@}\fontlist@}}
\def\nopages@{\output{\setbox\z@\box\@cclv \deadcycles\z@}%
 \alloc@5\toks\toksdef\@cclvi\output}
\let\galleys\nopages@
\newif\ifsyntax@
\newcount\countxviii@
\def\syntax{\syntax@true\dodummy@\countxviii@\count18
 \loop\ifnum\countxviii@>\m@ne\textfont\countxviii@=\dummyft@
 \scriptfont\countxviii@=\dummyft@\scriptscriptfont\countxviii@=\dummyft@
 \advance\countxviii@\m@ne\repeat                                           
 \dummyft@\tracinglostchars\z@\nopages@\frenchspacing\hbadness\@M}
\def\first@#1#2\end{#1}
\def\printoptions{\W@{Do you want S(yntax check),
  G(alleys) or P(ages)?}%
 \message{Type S, G or P, followed by <return>: }%
 \begingroup 
 \endlinechar\m@ne 
 \read\m@ne to\ans@
 \edef\ans@{\uppercase{\def\noexpand\ans@{%
   \expandafter\first@\ans@ P\end}}}%
 \expandafter\endgroup\ans@
 \if\ans@ P
 \else \if\ans@ S\syntax
 \else \if\ans@ G\galleys
 \else\message{? Unknown option: \ans@; using the `pages' option.}%
 \fi\fi\fi}
\def\alloc@#1#2#3#4#5{\global\advance\count1#1by\@ne
 \ch@ck#1#4#2\allocationnumber=\count1#1
 \global#3#5=\allocationnumber
 \ifalloc@\wlog{\string#5=\string#2\the\allocationnumber}\fi}
\def\document{\def\alloclist@{}\def\fontlist@{}}
\let\enddocument\bye

\let\proclaim\undefined
\let\footnote\undefined
\let\=\undefined
\let\>\undefined

\catcode`\@=\active
\message{... finished}

\expandafter\ifx\csname mathdefs.tex\endcsname\relax
  \expandafter\gdef\csname mathdefs.tex\endcsname{}
\else \message{Hey!  Apparently you were trying to
  \string\input{mathdefs.tex} twice.   This does not make sense.} 
\errmessage{Please edit your file (probably \jobname.tex) and remove
any duplicate ``\string\input'' lines}\endinput\fi




\catcode`\X=12\catcode`\@=11

\def\n@wcount{\alloc@0\count\countdef\insc@unt}
\def\n@wwrite{\alloc@7\write\chardef\sixt@@n}
\def\n@wread{\alloc@6\read\chardef\sixt@@n}
\def\r@s@t{\relax}\def\v@idline{\par}\def\@mputate#1/{#1}
\def\l@c@l#1X{\firstpart.#1}\def\gl@b@l#1X{#1}\def\t@d@l#1X{{}}

\def\crossrefs#1{\ifx\all#1\let\tr@ce=\all\else\def\tr@ce{#1,}\fi
   \n@wwrite\cit@tionsout\openout\cit@tionsout=\jobname.cit 
   \write\cit@tionsout{\tr@ce}\expandafter\setfl@gs\tr@ce,}
\def\setfl@gs#1,{\def\@{#1}\ifx\@\empty\let\next=\relax
   \else\let\next=\setfl@gs\expandafter\xdef
   \csname#1tr@cetrue\endcsname{}\fi\next}
\def\m@ketag#1#2{\expandafter\n@wcount\csname#2tagno\endcsname
     \csname#2tagno\endcsname=0\let\tail=\all\xdef\all{\tail#2,}
   \ifx#1\l@c@l\let\tail=\r@s@t\xdef\r@s@t{\csname#2tagno\endcsname=0\tail}\fi
   \expandafter\gdef\csname#2cite\endcsname##1{\expandafter
     \ifx\csname#2tag##1\endcsname\relax?\else\csname#2tag##1\endcsname\fi
     \expandafter\ifx\csname#2tr@cetrue\endcsname\relax\else
     \write\cit@tionsout{#2tag ##1 cited on page \folio.}\fi}
   \expandafter\gdef\csname#2page\endcsname##1{\expandafter
     \ifx\csname#2page##1\endcsname\relax?\else\csname#2page##1\endcsname\fi
     \expandafter\ifx\csname#2tr@cetrue\endcsname\relax\else
     \write\cit@tionsout{#2tag ##1 cited on page \folio.}\fi}
   \expandafter\gdef\csname#2tag\endcsname##1{\expandafter
      \ifx\csname#2check##1\endcsname\relax
      \expandafter\xdef\csname#2check##1\endcsname{}%
      \else\immediate\write16{Warning: #2tag ##1 used more than once.}\fi
      \multit@g{#1}{#2}##1/X%
      \write\t@gsout{#2tag ##1 assigned number \csname#2tag##1\endcsname\space
      on page \number\count0.}%
   \csname#2tag##1\endcsname}}

\def\multit@g#1#2#3/#4X{\def\t@mp{#4}\ifx\t@mp\empty%
      \global\advance\csname#2tagno\endcsname by 1 
      \expandafter\xdef\csname#2tag#3\endcsname
      {#1\number\csname#2tagno\endcsnameX}%
   \else\expandafter\ifx\csname#2last#3\endcsname\relax
      \expandafter\n@wcount\csname#2last#3\endcsname
      \global\advance\csname#2tagno\endcsname by 1 
      \expandafter\xdef\csname#2tag#3\endcsname
      {#1\number\csname#2tagno\endcsnameX}
      \write\t@gsout{#2tag #3 assigned number \csname#2tag#3\endcsname\space
      on page \number\count0.}\fi
   \global\advance\csname#2last#3\endcsname by 1
   \def\t@mp{\expandafter\xdef\csname#2tag#3/}%
   \expandafter\t@mp\@mputate#4\endcsname
   {\csname#2tag#3\endcsname\lastpart{\csname#2last#3\endcsname}}\fi}
\def\t@gs#1{\def\all{}\m@ketag#1e\m@ketag#1s\m@ketag\t@d@l p
\let\realscite\scite
\let\realstag\stag
   \m@ketag\gl@b@l r \n@wread\t@gsin
   \openin\t@gsin=\jobname.tgs \re@der \closein\t@gsin
   \n@wwrite\t@gsout\openout\t@gsout=\jobname.tgs }
\outer\def\localtags{\t@gs\l@c@l}
\outer\def\globaltags{\t@gs\gl@b@l}
\outer\def\newlocaltag#1{\m@ketag\l@c@l{#1}}
\outer\def\newglobaltag#1{\m@ketag\gl@b@l{#1}}

\newif\ifpr@ 
\def\m@kecs #1tag #2 assigned number #3 on page #4.%
   {\expandafter\gdef\csname#1tag#2\endcsname{#3}
   \expandafter\gdef\csname#1page#2\endcsname{#4}
   \ifpr@\expandafter\xdef\csname#1check#2\endcsname{}\fi}
\def\re@der{\ifeof\t@gsin\let\next=\relax\else
   \read\t@gsin to\t@gline\ifx\t@gline\v@idline\else
   \expandafter\m@kecs \t@gline\fi\let \next=\re@der\fi\next}
\def\pretags#1{\pr@true\pret@gs#1,,}
\def\pret@gs#1,{\def\@{#1}\ifx\@\empty\let\n@xtfile=\relax
   \else\let\n@xtfile=\pret@gs \openin\t@gsin=#1.tgs \message{#1} \re@der 
   \closein\t@gsin\fi \n@xtfile}

\newcount\sectno\sectno=0\newcount\subsectno\subsectno=0
\newif\ifultr@local \def\ultralocal{\ultr@localtrue}
\def\firstpart{\number\sectno}
\def\lastpart#1{\ifcase#1 \or a\or b\or c\or d\or e\or f\or g\or h\or 
   i\or k\or l\or m\or n\or o\or p\or q\or r\or s\or t\or u\or v\or w\or 
   x\or y\or z \fi}

\def\resetall{\global\advance\sectno by 1\subsectno=0
   \gdef\firstpart{\number\sectno}\r@s@t}
\def\resetsub{\global\advance\subsectno by 1
   \gdef\firstpart{\number\sectno.\number\subsectno}\r@s@t}
\def\newsection#1\par{\resetall\vskip0pt plus.3\vsize\penalty-250
   \vskip0pt plus-.3\vsize\bigskip\bigskip
   \message{#1}\leftline{\bf#1}\nobreak\bigskip}
\def\subsection#1\par{\ifultr@local\resetsub\fi
   \vskip0pt plus.2\vsize\penalty-250\vskip0pt plus-.2\vsize
   \bigskip\smallskip\message{#1}\leftline{\bf#1}\nobreak\medskip}


\newdimen\marginshift

\newdimen\margindelta
\newdimen\marginmax
\newdimen\marginmin

\def\margininit{       
\marginmax=3 true cm                  
				      
\margindelta=0.1 true cm              
\marginmin=0.1true cm                 
\marginshift=\marginmin
}    

\def\t@gsjj#1,{\def\@{#1}\ifx\@\empty\let\next=\relax\else\let\next=\t@gsjj
   \def\@@{p}\ifx\@\@@\else
   \expandafter\gdef\csname#1cite\endcsname##1{\citejj{##1}}
   \expandafter\gdef\csname#1page\endcsname##1{?}
   \expandafter\gdef\csname#1tag\endcsname##1{\tagjj{##1}}\fi\fi\next}
\newif\ifshowstuffinmargin
\showstuffinmarginfalse
\def\jjtags{\ifx\shlhetal\relax 
  \else
\ifx\shlhetal\undefinedcontrolseq
\else
\showstuffinmargintrue
\ifx\all\relax\else\expandafter\t@gsjj\all,\fi\fi \fi
}

\def\tagjj#1{\realstag{#1}\mginpar{\zeigen{#1}}}
\def\citejj#1{\rechnen{#1}\mginpar{\zeigen{#1}}}     

\def\rechnen#1{\expandafter\ifx\csname stag#1\endcsname\relax ??\else
                           \csname stag#1\endcsname\fi}

\newdimen\theight

\def\marginfont{\sevenrm}

\def\trymarginbox#1{\setbox0=\hbox{\marginfont\hskip\marginshift #1}%
		\global\marginshift\wd0 
		\global\advance\marginshift\margindelta}

\def \mginpar#1{%
\ifvmode\setbox0\hbox to \hsize{\hfill\rlap{\marginfont\quad#1}}%
\ht0 0cm
\dp0 0cm
\box0\vskip-\baselineskip
\else 
             \vadjust{\trymarginbox{#1}%
		\ifdim\marginshift>\marginmax \global\marginshift\marginmin
			\trymarginbox{#1}%
                \fi
             \theight=\ht0
             \advance\theight by \dp0    \advance\theight by \lineskip
             \kern -\theight \vbox to \theight{\rightline{\rlap{\box0}}%
\vss}}\fi}


\def\t@gsoff#1,{\def\@{#1}\ifx\@\empty\let\next=\relax\else\let\next=\t@gsoff
   \def\@@{p}\ifx\@\@@\else
   \expandafter\gdef\csname#1cite\endcsname##1{\zeigen{##1}}
   \expandafter\gdef\csname#1page\endcsname##1{?}
   \expandafter\gdef\csname#1tag\endcsname##1{\zeigen{##1}}\fi\fi\next}
\def\verbatimtags{\showstuffinmarginfalse
\ifx\all\relax\else\expandafter\t@gsoff\all,\fi}
\def\zeigen#1{\hbox{$\langle$}#1\hbox{$\rangle$}}
\def\margincite#1{\ifshowstuffinmargin\mginpar{\rechnen{#1}}\fi}
\def\margintag#1{\ifshowstuffinmargin\mginpar{\zeigen{#1}}\fi}

\def\(#1){\edef\dot@g{\ifmmode\ifinner(\hbox{\noexpand\etag{#1}})
   \else\noexpand\eqno(\hbox{\noexpand\etag{#1}})\fi
   \else(\noexpand\ecite{#1})\fi}\dot@g}

\newif\ifbr@ck
\def\eat#1{}
\def\[#1]{\br@cktrue[\br@cket#1'X]}
\def\br@cket#1'#2X{\def\temp{#2}\ifx\temp\empty\let\next\eat
   \else\let\next\br@cket\fi
   \ifbr@ck\br@ckfalse\br@ck@t#1,X\else\br@cktrue#1\fi\next#2X}
\def\br@ck@t#1,#2X{\def\temp{#2}\ifx\temp\empty\let\neext\eat
   \else\let\neext\br@ck@t\def\temp{,}\fi
   \def\teemp{#1}\ifx\teemp\empty\else\rcite{#1}\fi\temp\neext#2X}
\def\resetbr@cket{\gdef\[##1]{[\rtag{##1}]}}
\def\references{\resetbr@cket\newsection References\par}

\newtoks\symb@ls\newtoks\s@mb@ls\newtoks\p@gelist\n@wcount\ftn@mber
    \ftn@mber=1\newif\ifftn@mbers\ftn@mbersfalse\newif\ifbyp@ge\byp@gefalse
\def\defm@rk{\ifftn@mbers\n@mberm@rk\else\symb@lm@rk\fi}
\def\n@mberm@rk{\xdef\m@rk{{\the\ftn@mber}}%
    \global\advance\ftn@mber by 1 }
\def\rot@te#1{\let\temp=#1\global#1=\expandafter\r@t@te\the\temp,X}
\def\r@t@te#1,#2X{{#2#1}\xdef\m@rk{{#1}}}
\def\b@@st#1{{$^{#1}$}}\def\str@p#1{#1}
\def\symb@lm@rk{\ifbyp@ge\rot@te\p@gelist\ifnum\expandafter\str@p\m@rk=1 
    \s@mb@ls=\symb@ls\fi\write\f@nsout{\number\count0}\fi \rot@te\s@mb@ls}
\def\byp@ge{\byp@getrue\n@wwrite\f@nsin\openin\f@nsin=\jobname.fns 
    \n@wcount\currentp@ge\currentp@ge=0\p@gelist={0}
    \re@dfns\closein\f@nsin\rot@te\p@gelist
    \n@wread\f@nsout\openout\f@nsout=\jobname.fns }
\def\m@kelist#1X#2{{#1,#2}}
\def\re@dfns{\ifeof\f@nsin\let\next=\relax\else\read\f@nsin to \f@nline
    \ifx\f@nline\v@idline\else\let\t@mplist=\p@gelist
    \ifnum\currentp@ge=\f@nline
    \global\p@gelist=\expandafter\m@kelist\the\t@mplistX0
    \else\currentp@ge=\f@nline
    \global\p@gelist=\expandafter\m@kelist\the\t@mplistX1\fi\fi
    \let\next=\re@dfns\fi\next}
\def\symbols#1{\symb@ls={#1}\s@mb@ls=\symb@ls} 
\def\bigsymbol{\textstyle}
\symbols{\bigsymbol\ast,\dagger,\ddagger,\sharp,\flat,\natural,\star}
\def\ftnumbers{\ftn@mberstrue} \def\ftsymbols{\ftn@mbersfalse}
\def\paginal{\byp@ge} \def\resetftnumbers{\ftn@mber=1}
\def\ftnote#1{\defm@rk\expandafter\expandafter\expandafter\footnote
    \expandafter\b@@st\m@rk{#1}}

\long\def\jump#1\endjump{}
\def\ssum{\mathop{\lower .1em\hbox{$\textstyle\Sigma$}}\nolimits}

\def\qed{\nobreak\kern 1em \vrule height .5em width .5em depth 0em}
\def\newneq{\hbox{\rlap{\hbox to 1\wd9{\hss$=$\hss}}\raise .1em 
   \hbox to 1\wd9{\hss$\scriptscriptstyle/$\hss}}}
\def\subsetne{\setbox9 = \hbox{$\subset$}\mathrel{\hbox{\rlap
   {\lower .4em \newneq}\raise .13em \hbox{$\subset$}}}}
\def\supsetne{\setbox9 = \hbox{$\subset$}\mathrel{\hbox{\rlap
   {\lower .4em \newneq}\raise .13em \hbox{$\supset$}}}}

\def\vbar{\mathchoice{\vrule height6.3ptdepth-.5ptwidth.8pt\kern-.8pt}
   {\vrule height6.3ptdepth-.5ptwidth.8pt\kern-.8pt}
   {\vrule height4.1ptdepth-.35ptwidth.6pt\kern-.6pt}
   {\vrule height3.1ptdepth-.25ptwidth.5pt\kern-.5pt}}
\def\f@dge{\mathchoice{}{}{\mkern.5mu}{\mkern.8mu}}
\def\b@c#1#2{{\rm \mkern#2mu\vbar\mkern-#2mu#1}}
\def\b@b#1{{\rm I\mkern-3.5mu #1}}
\def\b@a#1#2{{\rm #1\mkern-#2mu\f@dge #1}}
\def\bb#1{{\count4=`#1 \advance\count4by-64 \ifcase\count4\or\b@a A{11.5}\or
   \b@b B\or\b@c C{5}\or\b@b D\or\b@b E\or\b@b F \or\b@c G{5}\or\b@b H\or
   \b@b I\or\b@c J{3}\or\b@b K\or\b@b L \or\b@b M\or\b@b N\or\b@c O{5} \or
   \b@b P\or\b@c Q{5}\or\b@b R\or\b@a S{8}\or\b@a T{10.5}\or\b@c U{5}\or
   \b@a V{12}\or\b@a W{16.5}\or\b@a X{11}\or\b@a Y{11.7}\or\b@a Z{7.5}\fi}}

\catcode`\X=11 \catcode`\@=12


\expandafter\ifx\csname citeadd.tex\endcsname\relax
\expandafter\gdef\csname citeadd.tex\endcsname{}
\else \message{Hey!  Apparently you were trying to
\string\input{citeadd.tex} twice.   This does not make sense.} 
\errmessage{Please edit your file (probably \jobname.tex) and remove
any duplicate ``\string\input'' lines}\endinput\fi

\sectno=-1   
\localtags
\jjtags
\NoBlackBoxes
\define\mr{\medskip\roster}
\define\sn{\smallskip\noindent}
\define\mn{\medskip\noindent}
\define\bn{\bigskip\noindent}
\define\ub{\underbar}
\define\wilog{\text{without loss of generality}}
\define\ermn{\endroster\medskip\noindent}

\define\dbcu{\dsize\bigcup}
\define \nl{\newline}
\magnification=\magstep 1
\documentstyle{amsppt}

{    
\catcode`@11

\ifx\alicetwothousandloaded@\relax
  \endinput\else\global\let\alicetwothousandloaded@\relax\fi

\gdef\subjclass{\let\savedef@\subjclass
 \def\subjclass##1\endsubjclass{\let\subjclass\savedef@
   \toks@{\def\usualspace{{\rm\enspace}}\eightpoint}%
   \toks@@{##1\unskip.}%
   \edef\thesubjclass@{\the\toks@
     \frills@{{\noexpand\rm2000 {\noexpand\it Mathematics Subject
       Classification}.\noexpand\enspace}}%
     \the\toks@@}}%
  \nofrillscheck\subjclass}
} 


\expandafter\ifx\csname alice2jlem.tex\endcsname\relax
  \expandafter\xdef\csname alice2jlem.tex\endcsname{\the\catcode`@}
\else \message{Hey!  Apparently you were trying to
\string\input{alice2jlem.tex}  twice.   This does not make sense.}
\errmessage{Please edit your file (probably \jobname.tex) and remove
any duplicate ``\string\input'' lines}\endinput\fi

\expandafter\ifx\csname bib4plain.tex\endcsname\relax
  \expandafter\gdef\csname bib4plain.tex\endcsname{}
\else \message{Hey!  Apparently you were trying to \string\input
  bib4plain.tex twice.   This does not make sense.}
\errmessage{Please edit your file (probably \jobname.tex) and remove
any duplicate ``\string\input'' lines}\endinput\fi

\def\renewcommand{\newcommand}	       
\edef\cite{\the\catcode`@}%
\catcode`@ = 11
\let\@oldatcatcode = \cite
\chardef\@letter = 11
\chardef\@other = 12
%
%
%
%
\def\@innerdef#1#2{\edef#1{\expandafter\noexpand\csname #2\endcsname}}%
%
%
\@innerdef\@innernewcount{newcount}%
\@innerdef\@innernewdimen{newdimen}%
\@innerdef\@innernewif{newif}%
\@innerdef\@innernewwrite{newwrite}%
%
%
%
\def\@gobble#1{}%
%
%
%
\ifx\inputlineno\@undefined
   \let\@linenumber = \empty 
\else
   \def\@linenumber{\the\inputlineno:\space}%
\fi
%
%
%
\def\@futurenonspacelet#1{\def\cs{#1}%
   \afterassignment\@stepone\let\@nexttoken=
}%
\begingroup 
\def\\{\global\let\@stoken= }%
\\ 
\endgroup
\def\@stepone{\expandafter\futurelet\cs\@steptwo}%
\def\@steptwo{\expandafter\ifx\cs\@stoken\let\@@next=\@stepthree
   \else\let\@@next=\@nexttoken\fi \@@next}%
\def\@stepthree{\afterassignment\@stepone\let\@@next= }%
%
%
%
\def\@getoptionalarg#1{%
   \let\@optionaltemp = #1%
   \let\@optionalnext = \relax
   \@futurenonspacelet\@optionalnext\@bracketcheck
}%
%
%
\def\@bracketcheck{%
   \ifx [\@optionalnext
      \expandafter\@@getoptionalarg
   \else
      \let\@optionalarg = \empty
      \expandafter\@optionaltemp
   \fi
}%
\def\@@getoptionalarg[#1]{%
   \def\@optionalarg{#1}%
   \@optionaltemp
}%
%
%
%
\def\@nnil{\@nil}%
\def\@fornoop#1\@@#2#3{}%
\def\@for#1:=#2\do#3{%
   \edef\@fortmp{#2}%
   \ifx\@fortmp\empty \else
      \expandafter\@forloop#2,\@nil,\@nil\@@#1{#3}%
   \fi
}%
\def\@forloop#1,#2,#3\@@#4#5{\def#4{#1}\ifx #4\@nnil \else
       #5\def#4{#2}\ifx #4\@nnil \else#5\@iforloop #3\@@#4{#5}\fi\fi
}%
\def\@iforloop#1,#2\@@#3#4{\def#3{#1}\ifx #3\@nnil
       \let\@nextwhile=\@fornoop \else
      #4\relax\let\@nextwhile=\@iforloop\fi\@nextwhile#2\@@#3{#4}%
}%
%
%
%
\@innernewif\if@fileexists
\def\@testfileexistence{\@getoptionalarg\@finishtestfileexistence}%
\def\@finishtestfileexistence#1{%
   \begingroup
      \def\extension{#1}%
      \immediate\openin0 =
         \ifx\@optionalarg\empty\jobname\else\@optionalarg\fi
         \ifx\extension\empty \else .#1\fi
         \space
      \ifeof 0
         \global\@fileexistsfalse
      \else
         \global\@fileexiststrue
      \fi
      \immediate\closein0
   \endgroup
}%
%
%
%
%
\def\bibliographystyle#1{%
   \@readauxfile
   \@writeaux{\string\bibstyle{#1}}%
}%
\let\bibstyle = \@gobble
%
%
\let\bblfilebasename = \jobname
\def\bibliography#1{%
   \@readauxfile
   \@writeaux{\string\bibdata{#1}}%
   \@testfileexistence[\bblfilebasename]{bbl}%
   \if@fileexists
      \nobreak
      \@readbblfile
   \fi
}%
\let\bibdata = \@gobble
%
%
\def\nocite#1{%
   \@readauxfile
   \@writeaux{\string\citation{#1}}%
}%
\@innernewif\if@notfirstcitation
%
%
\def\cite{\@getoptionalarg\@cite}%
%
%
\def\@cite#1{%
   \let\@citenotetext = \@optionalarg
   \printcitestart
   \nocite{#1}%
   \@notfirstcitationfalse
   \@for \@citation :=#1\do
   {%
      \expandafter\@onecitation\@citation\@@
   }%
   \ifx\empty\@citenotetext\else
      \printcitenote{\@citenotetext}%
   \fi
   \printcitefinish
}%
\newif\ifweareinprivate
\weareinprivatetrue
\ifx\shlhetal\undefinedcontrolseq\weareinprivatefalse\fi
\ifx\shlhetal\relax\weareinprivatefalse\fi
\def\@onecitation#1\@@{%
   \if@notfirstcitation
      \printbetweencitations
   \fi
   \expandafter \ifx \csname\@citelabel{#1}\endcsname \relax
      \if@citewarning
         \message{\@linenumber Undefined citation `#1'.}%
      \fi
     \ifweareinprivate
      \expandafter\gdef\csname\@citelabel{#1}\endcsname{%
\strut 
\vadjust{\vskip-\dp\strutbox
\vbox to 0pt{\vss\parindent0cm \leftskip=\hsize 
\advance\leftskip3mm
\advance\hsize 4cm\strut\openup-4pt 
\rightskip 0cm plus 1cm minus 0.5cm ?  #1 ?\strut}}
         {\tt
            \escapechar = -1
            \nobreak\hskip0pt\pfeilsw
            \expandafter\string\csname#1\endcsname
\pfeilso
            \nobreak\hskip0pt
         }%
      }%
     \else  
      \expandafter\gdef\csname\@citelabel{#1}\endcsname{%
            {\tt\expandafter\string\csname#1\endcsname}
      }%
     \fi  
   \fi
   \csname\@citelabel{#1}\endcsname
   \@notfirstcitationtrue
}%
%
%
\def\@citelabel#1{b@#1}%
%
%
\def\@citedef#1#2{\expandafter\gdef\csname\@citelabel{#1}\endcsname{#2}}%
%
%
%
\def\@readbblfile{%
   \ifx\@itemnum\@undefined
      \@innernewcount\@itemnum
   \fi
   \begingroup
      \def\begin##1##2{%
         \setbox0 = \hbox{\biblabelcontents{##2}}%
         \biblabelwidth = \wd0
      }%
      \def\end##1{}
      %
      %
      \@itemnum = 0
      \def\bibitem{\@getoptionalarg\@bibitem}%
      \def\@bibitem{%
         \ifx\@optionalarg\empty
            \expandafter\@numberedbibitem
         \else
            \expandafter\@alphabibitem
         \fi
      }%
      \def\@alphabibitem##1{%
         \expandafter \xdef\csname\@citelabel{##1}\endcsname {\@optionalarg}%
         \ifx\biblabelprecontents\@undefined
            \let\biblabelprecontents = \relax
         \fi
         \ifx\biblabelpostcontents\@undefined
            \let\biblabelpostcontents = \hss
         \fi
         \@finishbibitem{##1}%
      }%
      \def\@numberedbibitem##1{%
         \advance\@itemnum by 1
         \expandafter \xdef\csname\@citelabel{##1}\endcsname{\number\@itemnum}%
         \ifx\biblabelprecontents\@undefined
            \let\biblabelprecontents = \hss
         \fi
         \ifx\biblabelpostcontents\@undefined
            \let\biblabelpostcontents = \relax
         \fi
         \@finishbibitem{##1}%
      }%
      \def\@finishbibitem##1{%
         \biblabelprint{\csname\@citelabel{##1}\endcsname}%
         \@writeaux{\string\@citedef{##1}{\csname\@citelabel{##1}\endcsname}}%
         \ignorespaces
      }%
      %
      %
      \let\em = \bblem
      \let\newblock = \bblnewblock
      \let\sc = \bblsc
      \frenchspacing
      \clubpenalty = 4000 \widowpenalty = 4000
      \tolerance = 10000 \hfuzz = .5pt
      \everypar = {\hangindent = \biblabelwidth
                      \advance\hangindent by \biblabelextraspace}%
      \bblrm
      \parskip = 1.5ex plus .5ex minus .5ex
      \biblabelextraspace = .5em
      \bblhook
      \input \bblfilebasename.bbl
   \endgroup
}%
%
%
\@innernewdimen\biblabelwidth
\@innernewdimen\biblabelextraspace
%
%
%
\def\biblabelprint#1{%
   \noindent
   \hbox to \biblabelwidth{%
      \biblabelprecontents
      \biblabelcontents{#1}%
      \biblabelpostcontents
   }%
   \kern\biblabelextraspace
}%
%
%
%
\def\biblabelcontents#1{{\bblrm [#1]}}%
%
%
\def\bblrm{\rm}%
%
%
\def\bblem{\it}%
%
%
\def\bblsc{\ifx\@scfont\@undefined
              \font\@scfont = cmcsc10
           \fi
           \@scfont
}%
%
%
\def\bblnewblock{\hskip .11em plus .33em minus .07em }%
%
%
\let\bblhook = \empty
%
%
%
\def\printcitestart{[}
\def\printcitefinish{]}
\def\printbetweencitations{, }
\def\printcitenote#1{, #1}
%
%
%
\let\citation = \@gobble
%
%
%
\@innernewcount\@numparams
%
%
\def\newcommand#1{%
   \def\@commandname{#1}%
   \@getoptionalarg\@continuenewcommand
}%
%
%
\def\@continuenewcommand{%
   \@numparams = \ifx\@optionalarg\empty 0\else\@optionalarg \fi \relax
   \@newcommand
}%
%
%
\def\@newcommand#1{%
   \def\@startdef{\expandafter\edef\@commandname}%
   \ifnum\@numparams=0
      \let\@paramdef = \empty
   \else
      \ifnum\@numparams>9
         \errmessage{\the\@numparams\space is too many parameters}%
      \else
         \ifnum\@numparams<0
            \errmessage{\the\@numparams\space is too few parameters}%
         \else
            \edef\@paramdef{%
               \ifcase\@numparams
                  \empty  No arguments.
               \or ####1%
               \or ####1####2%
               \or ####1####2####3%
               \or ####1####2####3####4%
               \or ####1####2####3####4####5%
               \or ####1####2####3####4####5####6%
               \or ####1####2####3####4####5####6####7%
               \or ####1####2####3####4####5####6####7####8%
               \or ####1####2####3####4####5####6####7####8####9%
               \fi
            }%
         \fi
      \fi
   \fi
   \expandafter\@startdef\@paramdef{#1}%
}%
%
%
%
%
\def\@readauxfile{%
   \if@auxfiledone \else 
      \global\@auxfiledonetrue
      \@testfileexistence{aux}%
      \if@fileexists
         \begingroup
            \endlinechar = -1
            \catcode`@ = 11
            \input \jobname.aux
         \endgroup
      \else
         \message{\@undefinedmessage}%
         \global\@citewarningfalse
      \fi
      \immediate\openout\@auxfile = \jobname.aux
   \fi
}%
%
%
\newif\if@auxfiledone
\ifx\noauxfile\@undefined \else \@auxfiledonetrue\fi
%
%
%
%
\@innernewwrite\@auxfile
\def\@writeaux#1{\ifx\noauxfile\@undefined \write\@auxfile{#1}\fi}%
%
%
%
\ifx\@undefinedmessage\@undefined
   \def\@undefinedmessage{No .aux file; I won't give you warnings about
                          undefined citations.}%
\fi
%
%
\@innernewif\if@citewarning
\ifx\noauxfile\@undefined \@citewarningtrue\fi
%
%
%
\catcode`@ = \@oldatcatcode

\def\pfeilso{\leavevmode
            \vrule width 1pt height9pt depth 0pt\relax
           \vrule width 1pt height8.7pt depth 0pt\relax
           \vrule width 1pt height8.3pt depth 0pt\relax
           \vrule width 1pt height8.0pt depth 0pt\relax
           \vrule width 1pt height7.7pt depth 0pt\relax
            \vrule width 1pt height7.3pt depth 0pt\relax
            \vrule width 1pt height7.0pt depth 0pt\relax
            \vrule width 1pt height6.7pt depth 0pt\relax
            \vrule width 1pt height6.3pt depth 0pt\relax
            \vrule width 1pt height6.0pt depth 0pt\relax
            \vrule width 1pt height5.7pt depth 0pt\relax
            \vrule width 1pt height5.3pt depth 0pt\relax
            \vrule width 1pt height5.0pt depth 0pt\relax
            \vrule width 1pt height4.7pt depth 0pt\relax
            \vrule width 1pt height4.3pt depth 0pt\relax
            \vrule width 1pt height4.0pt depth 0pt\relax
            \vrule width 1pt height3.7pt depth 0pt\relax
            \vrule width 1pt height3.3pt depth 0pt\relax
            \vrule width 1pt height3.0pt depth 0pt\relax
            \vrule width 1pt height2.7pt depth 0pt\relax
            \vrule width 1pt height2.3pt depth 0pt\relax
            \vrule width 1pt height2.0pt depth 0pt\relax
            \vrule width 1pt height1.7pt depth 0pt\relax
            \vrule width 1pt height1.3pt depth 0pt\relax
            \vrule width 1pt height1.0pt depth 0pt\relax
            \vrule width 1pt height0.7pt depth 0pt\relax
            \vrule width 1pt height0.3pt depth 0pt\relax}

\def\pfeilsw{ \leavevmode 
            \vrule width 1pt height0.3pt depth 0pt\relax
            \vrule width 1pt height0.7pt depth 0pt\relax
            \vrule width 1pt height1.0pt depth 0pt\relax
            \vrule width 1pt height1.3pt depth 0pt\relax
            \vrule width 1pt height1.7pt depth 0pt\relax
            \vrule width 1pt height2.0pt depth 0pt\relax
            \vrule width 1pt height2.3pt depth 0pt\relax
            \vrule width 1pt height2.7pt depth 0pt\relax
            \vrule width 1pt height3.0pt depth 0pt\relax
            \vrule width 1pt height3.3pt depth 0pt\relax
            \vrule width 1pt height3.7pt depth 0pt\relax
            \vrule width 1pt height4.0pt depth 0pt\relax
            \vrule width 1pt height4.3pt depth 0pt\relax
            \vrule width 1pt height4.7pt depth 0pt\relax
            \vrule width 1pt height5.0pt depth 0pt\relax
            \vrule width 1pt height5.3pt depth 0pt\relax
            \vrule width 1pt height5.7pt depth 0pt\relax
            \vrule width 1pt height6.0pt depth 0pt\relax
            \vrule width 1pt height6.3pt depth 0pt\relax
            \vrule width 1pt height6.7pt depth 0pt\relax
            \vrule width 1pt height7.0pt depth 0pt\relax
            \vrule width 1pt height7.3pt depth 0pt\relax
            \vrule width 1pt height7.7pt depth 0pt\relax
            \vrule width 1pt height8.0pt depth 0pt\relax
            \vrule width 1pt height8.3pt depth 0pt\relax
            \vrule width 1pt height8.7pt depth 0pt\relax
            \vrule width 1pt height9pt depth 0pt\relax
      }


\def\widestnumber#1#2{}

\def\citewarning#1{\ifx\shlhetal\relax 
    \else
    \par{#1}\par
    \fi
}

\def\rm{\fam0 \tenrm}

\def\fakesubhead#1\endsubhead{\bigskip\noindent{\bf#1}\par}



%
%
%

%

\font\textrsfs=rsfs10
\font\scriptrsfs=rsfs7
\font\scriptscriptrsfs=rsfs5

\newfam\rsfsfam
\textfont\rsfsfam=\textrsfs
\scriptfont\rsfsfam=\scriptrsfs
\scriptscriptfont\rsfsfam=\scriptscriptrsfs

\edef\oldcatcodeofat{\the\catcode`\@}
\catcode`\@11

\def\Cal@@#1{\noaccents@ \fam \rsfsfam #1}

\catcode`\@\oldcatcodeofat


\expandafter\ifx \csname margininit\endcsname \relax\else\margininit\fi

\pageheight{8.5truein}
\topmatter
\title{Non Cohen Oracle c.c.c.  \\
Sh669} \endtitle
\author {Saharon Shelah \thanks {\null\newline I would like to thank 
Alice Leonhardt for the beautiful typing. \null\newline
 First Typed - 1/9/98 \null\newline
 Latest revision  - 02/Nov/5} \endthanks} \endauthor 
\affil{Institute of Mathematics\\
 The Hebrew University\\
 Jerusalem, Israel
 \medskip
 Rutgers University\\
 Mathematics Department\\
 New Brunswick, NJ  USA} \endaffil

\abstract   The oracle c.c.c. is closely related to Cohen forcing.
During an iteration we can ``omit a type"; i.e. preserve ``the
intersection of a given family of Borel sets of reals is empty"
provided that Cohen forcing satisfies it.  We generalize this to other
cases.  In \S1 we replace Cohen by ``nicely" definable c.c.c., do the
parallel of the oracle c.c.c. and end with a criterion for extracting
a subforcing (not a complete subforcing, $\lessdot$!) of a given 
nicely one and satisfying the oracle.  \endabstract
\endtopmatter
\document  
 
\newpage

\head {\S0 Introduction} \endhead  \resetall \sectno=0
\bigskip

This answers a question fom \cite[Ch.IV]{Sh:b} (the chapter dealing with the
oracle c.c.c.) asking to replace Cohen by e.g. random.  Later we will
deal with the parallel for oracle proper and for the case $\bar
\varphi_\alpha$ is a (definition of a) nep forcing.  An application
will appear in a work with T. Bartoszynski. 

How do we use this framework?  We start with a universe satisfying
$\diamondsuit_{\aleph_1}$ and probably $2^{\aleph_1} = \aleph_2$ and
choose $\langle S^*_i:i < \omega_2 \rangle,S^*_i \subseteq S^* 
\subseteq \omega_1$ such
that $S^*_i/{\Cal D}_{\omega_1}$ is strictly increasing and
$\diamondsuit_{S^*_{i+1} \backslash S^*_i}$ and for simplicity $S^*_i
\subseteq S^*_{i+1}$ where ${\Cal D}_{\omega_1}$ is the club filter on
$\omega_1$.  We choose by induction on $i < \omega_2$, a
c.c.c. forcing $\Bbb P_i$ of cardinality $\aleph_1$ a sequence
$\bar M^i = \langle M^i_\alpha:\alpha \in S^*_i \rangle$ and a 1-commitment.

They are increasing in the relevant sense and the work at limit stages
is done by the general claims here.  For $i=j+1$ we have some freedom
in choosing $\Bbb P_{i+1}$, usually $\Bbb P_{i+1} = \Bbb P_i *
{\underset\sim {}\to {\Bbb Q}_i}$.  So, working in $\bold V^{{\Bbb
P}_i}$, $\Bbb Q_i$ has to satisfy a 0-commitment on $S^*_i$, and we
like it to satisfy some task, possibly connected with some $X_i
\subseteq \Bbb R^{V[{\Bbb P}_i]}$, say $X_i = {\underset\sim {}\to
X_i}[G_{{\Bbb P}_i}]$ toward whatever task we have.  We essentially
have to choose $\bar M^i$ such that $\bar M^i \restriction S^*_j 
= \bar M^j$ but we have
freedom to choose $\langle \bar M^i_\alpha:\alpha \in S^*_i \backslash
S^*_j \rangle$ and a 0-commitment on $S^*_i \backslash S^*_j$.  The
reals generic for the chosen forcing notion as well as $M^i_\alpha$
for $\alpha \in S^*_i \backslash S^*_j$ can be chosen considering
$X_i$.  E.g. $M^i_\alpha$ can be the Mostowski Collapse of some
$M \prec ({\Cal H}(\aleph_2),\in)$ to which $\Bbb P_j,\bar M^j,x_j$
and ${\underset\sim {}\to X_j}$ belongs.

So really this corresponds to the omitting type as in \cite[XI]{Sh:e}.
\newpage

\head {\S1 Non-Cohen oracle c.c.c.} \endhead  \resetall \sectno=1
\bigskip

\demo{\stag{12.1} Hypothesis}  
\roster
\item "{$(a)$}"  we assume CH, moreover $\diamondsuit_{S^*}$ where
$S^* \subseteq \{\delta < \omega_1:\delta \text{ limit}\}$ is stationary.
\endroster
\enddemo
\bigskip

\definition{\stag{1.1a} Definition/Notation}  1)  $\bar M$ denotes a
sequence of the form $\langle M_\delta:\delta \in S \rangle, M_\delta$
a transitive countable model of $ZFC^-_*,\delta \subseteq M_\delta,
M_\delta \models ``\delta$ is countable", $S \subseteq S^*$ 
stationary satisfying: for every 
$X \subseteq \omega_1$, the set $\{\delta \in S^*:
X \cap \delta \in M_\delta\}$ is stationary. \nl
2)  ${\Cal D}$ denotes a normal filter on $\omega_1$ usually
extending ${\Cal D}_{\bar M}$ which is defined in \scite{12.2}(1) below (of
course, the default value is ${\Cal D}_{\bar M}$, see \scite{12.3}(1)).
\enddefinition
\bn
We first give the old definitions from \cite[IV]{Sh:f}
\definition{\stag{12.2} Definition}  1) ${\Cal D}_{\bar M}$ is \newline
$\{X \subseteq \omega_1:\text{for some } Y \subseteq \omega_1 
\text{ we have}:
\delta \in S_{\bar M} \cap X 
\Rightarrow Y \cap \delta \in M_\delta\}$. \newline
2) A forcing notion $\Bbb P$ of cardinality $\le \aleph_1$ satisfies the
$(\bar M,{\Cal D})$-c.c. \ub{if} for 
some (equivalently any) one to one $f:\Bbb P \rightarrow \omega_1$ the set:

$$
\align
\bigl\{ \delta \in S_{\bar M}:&\text{if } X \in M_\delta \text{ and }
\{y \in \Bbb P:f(y) < \delta \text{ and } f(y) \in X\} \\
  &\text{is predense in } \Bbb P \restriction \{y \in \Bbb P:f(y) < \delta\}
\text{ \ub{then} $X$ is predense in } \Bbb P \bigr\}
\endalign
$$
\mn
belongs to ${\Cal D}$ and $\Bbb P$ has minimal element
$\emptyset_{\Bbb P}$. \nl
3) If ${\Cal D} = {\Cal D}_{\bar M}$ we may write ``$\bar M$-c.c.".
Recall that ${\Cal D}^+ = \{A \subseteq \omega_1:\omega_1 \backslash A
\notin {\Cal D}\}$. \nl
4) Let $\bar M^1 \le \bar M^2$ if $\bar M^\ell = \langle
M^\ell_\delta:\delta \in S_\ell \rangle$ and $\{\delta:\delta \in S_1 
\backslash S_2$ or $\delta \in S_1 \cap S_2 \and M^1_\delta \ne
M^2_\delta\}$ is not stationary. \nl
5) A forcing notion $\Bbb P$ satisfies the $(\bar M,{\Cal
D})$-c.c. \ub{if}: $|\Bbb P| \le \aleph_0$ or for every $X \subseteq \Bbb P$
of cardinality $\le \aleph_1$ there is $\Bbb P_1 \lessdot \Bbb P$ of
cardinality $\aleph_1$ which includes $X$ and satisfies the $(\bar
M,{\Cal D})$-c.c. 
\enddefinition
\bigskip

\demo{\stag{12.3} Fact}  1) ${\Cal D}_{\bar M}$ 
is a normal filter on $\omega_1$. \newline
2) The $\bar M$-c.c. implies the c.c.c., and if $D_{\bar M} \subseteq
{\Cal D}$ (or just there is a normal filter $D' \supseteq D_{\bar M}
\cup {\Cal D}$) then the $(\bar M,{\Cal D})$-c.c.c. implies the c.c.c.
and if ${\Cal D}_2 \supseteq {\Cal D}_1 \supseteq {\Cal D}_{\bar M}$ 
are normal filters,
\u{then} the $(\bar M,{\Cal D}_1)$-c.c. implies the $(\bar M,{\Cal D}_2)$-c.c.
\nl
3) We can find $\langle S^*_\zeta:\zeta < \omega_2 \rangle$ such that
$S^*_\zeta \subseteq S^*,\zeta < \xi \Rightarrow S^*_\zeta \subseteq S^*_\xi
\text{ mod } {\Cal D}_{\bar M}$, \newline
$S^*_\zeta \subseteq S^*_{\zeta +1}$ and
$S^*_{\zeta +1} \backslash S^*_\zeta \in {\Cal D}^+_{\bar M}$. \nl
4) If $\bar M^1 \le \bar M^2$ and the forcing notion $\Bbb P_2$
satisfies the $(\bar M^2,{\Cal D})$-c.c. and $\Bbb P_1 \lessdot \Bbb
P_2$, \ub{then} $\Bbb P_1$ satisfies the $(\bar M^1,{\Cal
D})$-c.c. \nl
5) If $\bar M^1 \le \bar M^2$ and $P_1 \lessdot P_2$ and $\Bbb P_2$
satisfies the $\bar M^2$-c.c. \ub{then} $\Bbb P_1$ satisfies the $\bar
M^1$-c.c. 
\enddemo
\bigskip

\demo{Proof}  See \cite[IV]{Sh:f}, but for the reader we prove one. \nl
1) Without loss of generality $\Bbb P_2$ has set of elements
$\omega_1$.  As $\Bbb P_1 \lessdot \Bbb P_2$ there is a function
$f:\Bbb P_2 \rightarrow \Bbb P_1$ such that
\mr
\item "{$(*)_1$}"  $g \in \Bbb P_2 \wedge f(g) \le_{{\Bbb P}_1} p \in
\Bbb P_1 \Rightarrow p,g$ are compatible in $\Bbb P_2$.
\ermn
Let $g:\Bbb P_2 \times \Bbb P_2 \rightarrow \Bbb P_2$ be such that
\mr
\item "{$(*)_2$}"  if $p,g \in \Bbb P_2$ are compatible then $g(p,g)$
is a common upper bound and $p,g \in \Bbb P_1 \Rightarrow g(p,g) \in
\Bbb P_1$.
\ermn
So there is a club $E$ of $\omega_2$ which is closed under $f,g$ so
\mr
\item "{$(*)_3$}"   if $\delta \in S,{\Cal I} \subseteq \Bbb P_1 \cap
\delta$ is predense in $\Bbb P_1 \restriction \delta$ \ub{then} ${\Cal
I}$ is predense in $\Bbb P_2 \restriction \delta$. \nl
[Why?  If $g \in \Bbb P_2 \cap \delta$ then $f(g) \in \Bbb P_1 \cap
\delta$ so by the assumption on ${\Cal I},f(g)$ is compatible with
some $r_1 \in {\Cal I} in \Bbb P_1 \cap \delta$, so there is $r_2 \in
\Bbb P_1 \restriction \delta$ above $g$ and $r_2$.  By the definition
of $f,r_2,g$ are compatible in $\Bbb P_2$ hence $g(r_2,g)$ is a common
upper bound of them in $\Bbb P_2 \restriction \delta$.]
\endroster
\enddemo

\definition{\stag{12.4} Definition}  1)  We say ${\Cal Y} = (S,\Phi,
{\underset\sim {}\to {\bar \eta}},\bar \nu) =
(S^{\Cal Y},\Phi^{\Cal Y},{\underset\sim {}\to {\bar \eta}^{\Cal Y}},
\bar \nu^{\Cal Y})$ is a $0$-\ub{commitment} for $\bar M$ \ub{if} for some 
$E \in {\Cal D}_{\bar M}$:
\mr
\item "{$(a)$}"  $S \subseteq S^*,S \in {\Cal D}^+_{\bar M}$
\sn
\item "{$(b)$}"  $\bar \eta = \langle {\underset\sim {}\to \eta_\alpha}:
\alpha \in S \rangle,\Phi = \langle {\bar \varphi}_\alpha:\alpha \in S 
\rangle$ and if $\alpha \in S \cap E$ then 
${\bar \varphi}_\alpha \in M_\alpha$ and \nl
$M_\alpha \models ``{\bar \varphi}_\alpha$ is an absolute definition of a
c.c.c. forcing notion called $\Bbb Q_\alpha$ with generic
real ${\underset\sim {}\to \eta_\alpha}"$; note, absolute here means that
forcing extensions of $M_\alpha$, preserve predensity of countable sets (in the
sense of $M_\alpha$) order and incompatibility
\sn
\item "{$(c)$}"  $\bar \nu = \langle \nu_\alpha:\alpha \in S \rangle$
where $\nu_\alpha \in 
{}^\omega \omega$ and for every $\alpha \in S \cap E$ the
real $\nu_\alpha$ is $(\Bbb Q_\alpha,
{\underset\sim {}\to \eta_\alpha})$-generic over $M_\alpha$.
\ermn
We ignore $\bar M$ if clear from the context.  We can replace $\bar M$
by $(\bar M,{\Cal D})$ if above $E \in {\Cal D},S \in {\Cal D}^+$.   \nl
1A) A forcing notion $\Bbb P$ of cardinality $\le \aleph_1$ satisfies the 
$0$-commitment ${\Cal Y}$ \ub{if}: $\Bbb P$ is an $\bar M$-c.c. forcing notion and
for any one-to-one mapping $h:\Bbb P \rightarrow \omega_1$ for some 
$E \in {\Cal D}_{\bar M}$ we have
\mr
\item "{$(d)$}"   if 
$\alpha \in S \cap E$ then $\Vdash_{\Bbb P}$ ``the real $\nu_\alpha$ is a 
$(\Bbb Q_\alpha,{\underset\sim {}\to \eta_\alpha})$-generic real 
over $M_\alpha[\alpha \cap h''{\underset\sim {}\to G_{\Bbb P}}]$".
\ermn
2) Let $\Bbb P \in {\Cal H}(\aleph_2)$ be an $\bar M$-c.c. forcing notion.  We
say that ${\Cal Y} = (S,{\underset\sim {}\to {\bar \Phi}},
{\underset\sim {}\to {\bar \eta}},{\underset\sim {}\to {\bar \nu}}) 
= (S^{\Cal Y},{\underset\sim {}\to \Phi^{\Cal Y}},
{\underset\sim {}\to {\bar \eta}^{\Cal Y}},\bar \nu)$
is a $1$-commitment on $\Bbb P$ if: for any $\bar N$ satisfying $(*)_1$ below, the 
clauses (a)-(d) of $(*)_2$ below hold
\mr
\item "{$(*)_1$}"  $\bar N = \langle N_\alpha:\alpha < \omega_1 \rangle$ is
increasing continuous, $N_\alpha \prec ({\Cal H}(\aleph_2),\in)$ is
countable, $\bar N \restriction (\alpha +1) \in N_{\alpha +1}$ and
$\{\bar M,\Bbb P\} \subseteq \dsize \bigcup_{\alpha < \omega_1} N_\alpha$
{\roster
\itemitem{ $(*)_2$(a) }  $S \subseteq \text{ Dom}(\bar M) \subseteq
S^*,S \in {\Cal D}^+_{\bar M}$
\sn
\itemitem{ (b) }  $\bar \eta = \langle {\underset\sim {}\to \eta_\alpha}:
\alpha \in S \rangle,\Phi = 
\langle {\underset\sim {}\to {\bar \varphi}_\alpha}:\alpha \in S \rangle$
so $({\underset\sim {}\to {\bar \varphi}_\alpha},
{\underset\sim {}\to \eta_\alpha})$ is a $\Bbb P$-name of a pair as in
\scite{12.4}(1)(a), both are hereditarily countable over $\Bbb P$ 
\sn
\itemitem{ (c) }  $\bar \nu = \langle {\underset\sim {}\to \nu_\alpha}:
\alpha \in S \rangle$ and ${\underset\sim {}\to \nu_\alpha}$ a 
$\Bbb P$-name of a 
real given by countably many conditions
\sn
\itemitem{ (d) }  the set of the $\alpha \in S$ satisfying the following
belongs to \nl
${\Cal D}_{\bar M} + S:{\underset\sim {}\to {\bar \varphi}_\alpha} \in
M_\alpha$, Mos Col$_{N_\alpha}(N_\alpha) \in M_\alpha$, and letting \nl
$\Bbb P'_\alpha =$ Mos Col$_{N_\alpha}(N_\alpha^{\Bbb P}) = 
\{\text{Mos Col}_{N_\alpha}(x):x \in \Bbb P \cap N_\alpha\} \in
M_\alpha$ we have $M_\alpha \models ``
{\underset\sim {}\to {\bar \varphi}_\alpha}$ is a $\Bbb P'_\alpha$-name of an
absolute definition of a c.c.c. forcing with generic real 
${\underset\sim {}\to \eta_\alpha}$" and $\Vdash_{\Bbb P}$ ``the real
${\underset\sim {}\to \nu_\alpha}$ is a 
$(\Bbb Q_{{\underset\sim {}\to {\bar \varphi}^\alpha}},
{\underset\sim {}\to \eta_\alpha})$-generic real over $M_\alpha
[{\underset\sim {}\to G_{\Bbb P}}]"$.
\endroster}
\ermn
For simplicity the reader may concentrate on the case $\langle
(\bar \varphi_\alpha,{\underset\sim {}\to \eta_\alpha}):\alpha \in S
\rangle \in \bold V$. \nl
3) Let $IS = \{(\Bbb P,{\Cal Y},\bar M):
\Bbb P \in {\Cal H}(\aleph_2)$ is an $\bar M$-c.c.
forcing notion and ${\Cal Y}$ \newline

$\qquad \qquad \qquad \qquad$ is a 1-commitment on $\Bbb P\}$. \newline
We shall omit $\bar M$ if clear from the context.
We can replace $\bar M$ by $(\bar M,{\Cal D})$-naturally and write 
IS$_{\Cal D}$, but the claims are the same.
 \nl
4) For $(\Bbb P^\ell,{\Cal Y}^\ell,\bar M^\ell) \in IS \, (\ell =
1,2)$ let $(\Bbb P^1,{\Cal Y}^1,\bar M^1)
\le^* (\Bbb P^2,{\Cal Y}^2,\bar M^2)$ means $\bar M^1 \le \bar M^2,
\Bbb P^1 \lessdot \Bbb P^2$ and for some $E \in {\Cal D}_{{\bar M}^1}$
we have \nl
$S^{{\Cal Y}_1} \cap E \subseteq S^{{\Cal Y}_2} \cap E,
{\underset\sim {}\to {\bar \Phi}^{{\Cal Y}^1}} \restriction (S^{{\Cal Y}_1}
\cap E) = {\underset\sim {}\to {\bar \Phi}^{{\Cal Y}^2}} \restriction 
(S^{\Cal Y} \cap E)$, \newline
${\underset\sim {}\to {\bar \eta}^{{\Cal Y}^1}} \restriction 
(S^{{\Cal Y}^1} \cap E) = {\underset\sim {}\to {\bar \eta}^{{\Cal Y}^2}} 
\restriction (S^{{\Cal Y}^1} \cap E)$ and
${\underset\sim {}\to {\bar \nu}^{{\Cal Y}^1}} \restriction (S^{{\Cal Y}^1}
\cap E) = {\underset\sim {}\to {\bar \nu}^{{\Cal Y}^2}} \restriction 
(S^{{\Cal Y}^1} \cap E)$.  We call $E$ a witness to 
$(\Bbb P^1,{\Cal Y}^1,\bar M^1) \le^*
(\Bbb P^2,{\Cal Y}^2,\bar M^2)$.   
\enddefinition
\bn

We point out the connection between $0$-commitment and $1$-commitment. \nl
\margintag{12.4A}\ub{\stag{12.4A} Fact}:  1) If ${\Cal Y}$ is a 1-commitment on 
$\Bbb P$ and $\Bbb P$ an
$\bar M$-c.c. forcing notion of cardinality $\le \aleph_1$, \ub{then}
$\Vdash_{\Bbb P} ``{\Cal Y}[{\underset\sim {}\to G_{\Bbb P}}] = (S^y,
{\underset\sim {}\to \Phi^{\Cal Y}}[{\underset\sim {}\to G_{\Bbb P}}],
{\underset\sim {}\to {\bar \eta}^{\Cal Y}}
[{\underset\sim {}\to G_{\Bbb P}}],{\underset\sim {}\to {\bar \nu}}
[{\underset\sim {}\to G_{\Bbb P}}])"$ is a 0-commitment so we call it
${\Cal Y}[{\underset\sim {}\to G_{\Bbb P}}]$.  Note 
$\underset\sim {}\to \eta[{\underset\sim {}\to G_{\Bbb P}}]$ 
is still a name.  \nl
2) If $\Bbb P = \{ \emptyset\}$ (the trivial forcing) then: ${\Cal Y}$ is a
1-commitment on $\Bbb P$ \ub{iff} ${\Cal Y}$ is a 0-commitment. \nl
3) If $\langle \bar M^i:i < \zeta \rangle$ is $\le$-increasing, $\zeta <
\omega_2$ and Dom$(\bar M^i) \backslash S$ is not stationary,
\ub{then} there is $\bar M$, Dom$(\bar M) = S$ such that $i < \zeta
\Rightarrow \bar M^i \le \bar M$. \nl
4) Increasing $\bar M$ preserves everything.
\bn
As a warm-up (see \cite{Sh:630} for more)
\proclaim{\stag{12.4.xad} Claim}  1) Assume
\mr
\item "{$(a)$}"  $M$ is a countable transitive model of ZFC$^-,M
\models ``\Bbb P_1$ is a countable forcing notion"
\sn
\item "{$(b)$}"  $M \models ``\varphi$ is an absolute definition of
c.c.c. forcing notion $\Bbb Q^\varphi$ with generic $\underset\sim {}\to
\eta:\alpha_1 \rightarrow \alpha_2$"
\sn
\item "{$(c)$}"  $\nu$ is $(M,\Bbb Q^\varphi)$-generic sequence, i.e.,
there is $G \subseteq (\Bbb Q^\varphi)^M$ generic over $M$ such that $\nu
= \underset\sim {}\to \eta[G]$.
\ermn
\ub{Then} we can find a countable $\Bbb P_2$ such that
\mr
\item "{$(\alpha)$}"  $\Bbb P_1 \subseteq_{ic} \Bbb P_2$ and every
${\Cal J} \in M$ which is predense in $\Bbb P_1$ is predense in $\Bbb
P_2$
\sn
\item "{$(\beta)$}"  $\Vdash_{{\Bbb P}_2} ``\nu$ is $(M',\Bbb
Q^\varphi)$-generic sequence where $M' = M[G_{{\Bbb P}_2} \cap \Bbb
P_1]"$.
\ermn
2) Similarly for $\varphi$ defining a nep forcing.
\endproclaim
\bigskip

\demo{Proof}  1) In $M$ we can define $\Bbb P^+ = \Bbb P_1 * 
({\underset\sim {}\to {\Bbb Q}^\varphi})
^{M[{\underset\sim {}\to G_{{\Bbb P}_1}}]}$, 
now as $Q^\varphi$ is absolutely c.c.c., we know
that $q \mapsto (\emptyset,q)$ is a complete embedding of $(\Bbb
Q^\varphi)^M$ into $\Bbb P^+$.  So if $G_* \subseteq (\Bbb
Q^\varphi)^M$ is generic over $M$ such that $\nu = \nu[G]$ then let
$\Bbb P^*_2 = \{(p,\underset\sim {}\to q) \in \Bbb P_1 *
({\underset\sim {}\to {\Bbb Q}^\varphi})^{M[{\underset\sim {}\to
G_{{\Bbb P}_1}}]}:(p,\underset\sim {}\to q)$ is compatible with
$(\emptyset,q')$ for every $q' \in G_*\}$.  Now check. \nl
2) See \cite{Sh:630}.  \hfill$\square_{\scite{12.4.xad}}$\margincite{12.4.xad}
\enddemo
\bigskip

\proclaim{\stag{12.5} Crucial Claim}  In $IS$, any $\le^*$-increasing
$\omega$-chain has an upper bound.
\endproclaim
\bigskip

\demo{Proof}  So assume $(\Bbb P^n,{\Cal Y}^n,\bar M^n) \in IS$ 
and $(\Bbb P^n,{\Cal Y}^n,\bar M^n)
\le^* (\Bbb P^{n+1},{\Cal Y}^{n+1},\bar M^{n+1})$ for $n < \omega$,
let $\bar M$ be such that 
$\bar M \ge \bar M^n$ for each $n$; so let $E_n \in {\Cal D}_{\bar M}$ 
witness it.  For simplicity assume that above any $p \in \Bbb P^n$ there 
are two incompatible elements, and $0 \in \Bbb P^0$ is minimal in all
$\Bbb P^n$, 
i.e. is $\emptyset_{{\Bbb P}_n}$.  Without loss of
generality the set of elements of $\Bbb P^n$ is $\subseteq \omega_1$ and
$\omega_1 \backslash \dsize \bigcup_{n < \omega} \Bbb P^n$ has cardinality
$\aleph_1$ and let $X^*$ be such that $\dsize \bigcup_{n < \omega}
\Bbb P_n \subseteq X^* \subseteq \omega_1$ and
$|X^* \backslash \dbcu_{n < \omega} \Bbb P^n| = 
\aleph_1$; notation helps in a future use, also there we replace $\omega$ 
by a (countable) ordinal of cofinality $\aleph_0$.  We can define functions
$F_n,f_{n,m+1},F_{n,m,\ell}$ 
(when $n < m < \omega,\ell < \omega)$ such that
\mr
\item "{$(a)_n$}"  if $p,q \in \Bbb P_n$ are compatible then $F_n(p,q)
\in \Bbb P_n$ is a common upper bound
\sn
\item "{$(b)_{n,m}$}"  if $n \le m$ and
$p \in \Bbb P_m$, then $\langle F_{n,m,\ell}(p):\ell
< \omega \rangle$ is a maximal antichain of $\Bbb P_n$, such that for each
$\ell$: \newline
\ub{either} $p,F_{n,m,\ell}(p)$ are incompatible (in $\Bbb P_m$) \newline
\ub{or} $p$ is compatible in $\Bbb P_m$ with every $q \in \Bbb P_n$
which  is above $F_{n,m,\ell}(p)$
\sn
\item "{$(c)_{n,m}$}"  if $p \in \Bbb P_m,g \in \Bbb P_n$ then 
$F_{n,m}(p,g) \in \Bbb P_n$ and if there is $r,g \le_{{\Bbb P}_n} r$
and $r,p$ is incompatible in $\Bbb P_m$ then $p,F_{n,m}(p,g)$ are 
incompatible in $\Bbb P_m$.
\ermn
Let $E$ be a club of $\omega_1$, such that $\delta \in E 
\Rightarrow \delta$ is closed under $F_n,F_{n,m,\ell}$ and \newline
otp$(X^* \cap \delta 
\backslash \dsize \bigcup_{n < \omega} \Bbb P^n) = \delta$. 
\newline
We would like to define a forcing notion $\Bbb P^\omega$ with universe $X^*$,
and $1$-commitment ${\Cal Y}^\omega$, and functions $F_\omega,F_{n,\omega,
\ell}$ satisfying the natural requirements.  First, let \newline
$E_\omega = \dsize \bigcap_{n < \omega} E_n \cap E,S^{{\Cal Y}^\omega} 
= \dsize \bigcup_{n < \omega}
S^{{\Cal Y}^n} \cap E_\omega$, and for $\alpha \in S^{{\Cal Y}^\omega}$ 
the triple (${\underset\sim {}\to {\bar \varphi}^{{\Cal Y}^\omega}_\alpha},  
{\underset\sim {}\to \eta^{{\Cal Y}^\omega}_\alpha},
{\underset\sim {}\to \nu^{{\Cal Y}^\omega}_\alpha})$ is 
$({\underset\sim {}\to {\bar \varphi}^{{\Cal Y}^{n(\alpha)}}_\alpha},
{\underset\sim {}\to \eta^{{\Cal Y}^{n(\alpha)}}_\alpha},
{\underset\sim {}\to \nu^{{\Cal Y}^{n(\alpha)}}_\alpha})$ where $n(\alpha) =
\text{ Min}\{n:\alpha \in S^{{\Cal Y}^n}\}$.
\mn
Defining $\Bbb P_\omega,
F_\omega,F_{n,\omega,\ell}$ is harder, so we first define
$AP$, a set of approximations to it.  A member $t$ of $AP$ has the form
$(\delta^t,\Bbb P^t,F^t_\omega,
F^t_{n,\omega,\ell},\Gamma^t)_{\ell < \omega}$ 
satisfying
\mr
\item "{$(\alpha)$}"  $\delta^t \in E$
\sn
\item "{$(\beta)$}"  $\Bbb P^t$ is a forcing notion with set of elements
$\subseteq X^* \cap \delta^t$ and $\supseteq \delta^t \cap \dbcu_n
\Bbb P_n$ and $0 \le_{{\Bbb P}^t} p$ for every $p \in \Bbb P^t$
\sn
\item "{$(\gamma)$}" $\Bbb P^t \restriction (\Bbb P^n \cap \delta^t) =
\Bbb P^n \restriction (\Bbb P^n \cap \delta^t)$
\sn     
\item "{$(\delta)$}"  if $p,q \in \Bbb P^t$ are compatible in $\Bbb P^t$ then
$F^t_\omega(p,q)$ is such an upper bound
\sn
\item "{$(\varepsilon)$}"  if $p \in \Bbb P^t,n < \omega$ then $\langle
F^t_{n,\omega,\ell}(p):\ell < \omega \rangle$ is a maximal antichain of
$\Bbb P_n$, the members are $< \delta^t$, and for each $\ell$, \ub{either}
$p,F^t_{n,\omega,\ell}(p)$ are incompatible in $\Bbb P^t$ \ub{or}
$(\forall q \in 
\Bbb P_n \cap \delta)(\Bbb P_n \models ``F_{n,\omega,\ell}(p) \le
q" \Rightarrow p,q$ are compatible in $\Bbb P^t)$ and for at least one
$\ell$ the second case occurs
\sn 
\item "{$(\zeta)$}"  if $p \in \Bbb P^t \cap \Bbb P_m \backslash
\dbcu_{\ell < m} \Bbb P_\ell$ then $F^t_{n,\omega,\ell}(p) = F_{n,m,\ell}(p)$
\sn
\item "{$(\eta)$}"  $\Gamma^t$ is a sequence $\langle \bar p^t_\zeta:\zeta
< \zeta^t \rangle,\zeta^t < \omega_1$ and $\bar p^t_\zeta$ is a sequence of
length $\le \omega$ of members of $\Bbb P^t$ which form a maximal antichain
(of $\Bbb P^t$)
\sn
\item "{$(\theta)$}"   if $p \in \Bbb P^t$ and $n < m < \omega$ and $r
\in \Bbb P_n \cap \delta^t$ and $[r \le r' \in \Bbb P_n \cap \delta^t
\Rightarrow r',p$ are compatible in $\Bbb P^t]$, \ub{then} the 
set $\{F^t_{m,\omega,\ell}(p):\ell < \omega$ and $p$ is compatible
with $F^t_{m,\omega,\ell}(p)$ in $\Bbb P^t\}$ satisfies: if $r \le g
\in \Bbb P_n$ then in $\Bbb P_m,q$ is compatible with some member of
this set
\sn
\item "{$(\iota)$}"     if $\zeta < \zeta^t$ and 
$n < \omega$ \ub{then}: \newline
$\{F^t_{n,\omega,\ell}(p^t_{\zeta,k}):k < \omega,\ell < \omega$ and
$p^t_{\zeta,k},F^t_{n,\omega,\ell}(p^t_{\zeta,k})$ are compatible in 
$\Bbb P^t\}$ is a predense subset of $\Bbb P_n$.  Note that trivial,
this subset is predense in $\Bbb P_n \cap \delta$; similarly in clause
$(\kappa)$ 
\ermn
Moreover, 
\mr
\item "{$(\kappa)$}"  if $p^* \in \Bbb P^t$ and $n < \omega$ and
$\zeta < \zeta^t$  then 
$$
\align
{\Cal I}^t_{\zeta,n,p^*} =: \bigl\{r' \in \Bbb P_n \cap \delta:&(i) 
\quad r',p^* \text{ incompatible in } \Bbb P^t \text{ \ub{or}} \\
  &(ii) \quad \text{for some } k < \omega 
\text{ and } p' \text{ we have} \\
  &\quad \quad (\forall r'')[r' \le r'' \in \Bbb P_n \cap 
\delta \rightarrow \{r'',p'\} \\
  &\quad \quad \text{has an upper bound in } \Bbb P^t] \text{ and }
    p^* \le_{{\Bbb P}^t} p',p^t_{\zeta,k} \le_{{\Bbb P}^t} p' \bigr\}.
\endalign
$$
is predense in $\Bbb P_n$. 
\endroster
\enddemo
\bn
\centerline {$* \qquad * \qquad *$}
\bn
We define the (natural) partial order 
$\le^*$ on $AP$: for $t,s, \in AP$ as follows; we let $t \le^* s$ \ub{iff}:
\mr
\widestnumber\item{$(iii)$}
\item "{$(i)$}"  $\delta^t \le \delta^s$
\sn
\item "{$(ii)$}"  $\Bbb P^t \subseteq \Bbb P^s$
\sn
\item "{$(iii)$}"  $F^t_\omega \subseteq F^s_\omega$
\sn
\item "{$(iv)$}"  $F^t_{n,\omega,\ell} \subseteq F^s_{n,\omega,\ell}$
\sn
\item "{$(v)$}"   $\Gamma^t$ is an initial segment of $\Gamma^s$.
\ermn
\ub{Fact A}:  $AP \ne \emptyset$.
\bigskip

\demo{Proof}  Easy: choose $\delta \in E$, let $\Bbb P^t = \bigl(
\dsize \bigcup_{n < \omega} \Bbb P^n \bigr) \restriction \delta,
F_\omega(p,q) =
F_{n(p,q)}(p,q)$ where $n(p,q) = \text{ Min}\{n:p \in \Bbb P_n$ and 
$q \in \Bbb P_n\}$.
\newline
For $n < \omega,p \in \Bbb P^t \cap \delta$ let $\langle F^t_{n,\omega,\ell}
(p):\ell < \omega \rangle$ be $\langle F_{n,m,\ell}(p):\ell < \omega \rangle$
for the first $m \ge n$ such that $p \in \Bbb P_m$. \newline
Lastly, $\Gamma =$ empty sequence.
\enddemo
\bn
\ub{Fact B}: If $t \in AP$ and $\delta^t < \delta \in E$, \ub{then} there is
$s$ satisfying $t \le^* s \in AP$ with 
$\delta^s \ge \delta,\zeta^s = \zeta^t$.
\bigskip

\demo{Proof}  Without loss of generality 
$t,\langle \Bbb P^n \restriction \delta:
n < \omega \rangle,X^* \cap \delta$ belongs to $M_\delta$ and $\delta \in E
\cap \dsize \bigcap_{n < \omega} E_n$ and $X^* \cap \delta \backslash
\dbcu_{n < \omega} \Bbb P_n \backslash \delta^t$ is infinite.  
[Why?  As $\emptyset \notin
{\Cal D}_{\bar M}$ and we can increase $\delta$.]   So (for the last
phrase see the proof of \scite{12.3}(4))
\mr
\item "{$(*)$}"  any ${\Cal J} \in M_\delta$ which is a predense
subset of $\Bbb P_n \restriction \delta$ is a predense subset of $\Bbb
P_n$ and $n < m \Rightarrow \Bbb P_n \restriction \delta \lessdot
\Bbb P_m \restriction \delta$.
\ermn
Let $A = \Bbb P^t \subseteq 
X^* \cap \delta^t,B = \dsize \bigcup_{n < \omega} \Bbb P^n \cap
\delta$.  We define a forcing notion $\Bbb Q$, with set of elements $\subseteq
A \times B$ identifying $(p,0)$ with $p$ and $(0,q)$ with $q$.  Now $(p,q) \in
A \times B$ belongs to $\Bbb Q$ \ub{iff}: $p=0$ \ub{or} $q=0$ \ub{or} there
are $r \in A \cap B$ and $n=n(p,q)$ such that:  $\Bbb P_n \models ``r \le q"$, 
and $(\forall r')[r \le r' \in \Bbb P_n \cap \delta^t \rightarrow r',p
\text{ compatible in } \Bbb P^t]$; we call such $r$ a witness and $n$ a possible
value for $n(p,q)$.  The order on $\Bbb Q$ is $(p,q) \le (p',q')
\Leftrightarrow p \le_{{\Bbb P}^t} p' \and \dsize \bigvee_n 
q \le_{{\Bbb P}_n} q'$. \newline
Now note
\mr
\item "{$(\alpha)$}"  $\Bbb Q \restriction A = \Bbb P^t$
\sn
\item "{$(\beta)$}"  $\Bbb Q \restriction B = \dsize \bigcup_{n <
\omega} \Bbb P_n \restriction \delta$
\sn
\item "{$(\gamma)$}"  If $(p,q) \in \Bbb Q,m = n(p,q)$ and 
$q \in \Bbb P_m \restriction \delta$ and 
$\Bbb P_m \restriction \delta \models 
``q \le q'\,"$ and $\Bbb P^t \models ``p' \le
p"$, \ub{then} $(p',q') \in \Bbb Q$ and $\Bbb Q \models ``(p,q) \le (p,q')"$
\sn
\item "{$(\delta)$}"  if $(p,q) \in \Bbb Q$ and $n=n(p,q) \le m < \omega$,
\ub{then} for some $q_1$ we have: $(p,q) \le (p,q_1) \in \Bbb Q$ and
$n(p,q_1) = m$, or at least $m$ is a possible value for $n(p,q_1)$. \nl
[Why?  Let $r \in \Bbb P_{n(p,q)}$ be a witness in particular $r$ is
compatible with $p$ in $\Bbb P_t$.
By clause $(\theta)$ of the Definition of AP the set ${\Cal J} =
\{F^t_{m,\omega,\ell}(p):\ell < \omega$ and $p$ is compatible with
$F^t_{m,\omega,\ell}(p)$ in $\Bbb P^t\}$ is predense above $r$ in
$\Bbb P_m$.  $\Bbb P_n \models r \le g$ hence $\Bbb P_m \models r \le
g$ so for some $\ell,F^t_{m,\omega,\ell}(p) \in {\Cal J}$ is
compatible with $q$ in $\Bbb P_m$ so
there is $q_1 \in \Bbb P_m \cap \delta$ such that $\Bbb P_m \models q
\le q_1 \wedge F^t_{m,\omega,\ell}(p) \le q_1$.  So $(p,q_1) \in \Bbb
Q$ as witnessed by $m$ and $r' = F^t_{m,\omega,\ell}(p)$, is as required.]
\sn
\item "{$(\varepsilon)$}"  $\Bbb P_n \restriction 
\delta \lessdot \Bbb Q$ \newline
[Why?  Let $(p^0,q^0) \in \Bbb Q$, of 
course, we can replace this pair by any larger
one, so by clause $(\delta)$ above \wilog \, some $m \in [n,\omega)$, is a
possible value for $n(p^0,q^0)$ so we have
$q^0 \in \Bbb P_m \restriction \delta$, hence recalling that $\Bbb P_n
\restriction \delta \lessdot \Bbb P_m \restriction \delta$ there is $q^1 \in 
\Bbb P_n \restriction \delta$ such that:
$$
(\forall r \in \Bbb P_n)(\Bbb P_n \restriction \delta \models q^1 \le r \Rightarrow r,
q^0 \text{ compatible in } \Bbb P_m \restriction \delta).
$$
Assume $q^1 \le r \in \Bbb P_n \restriction \delta$.  So $r,q^0$ are compatible
in $\Bbb P_m \restriction \delta$ hence has a common upper bound $q^2
\in \Bbb P_m \restriction \delta$. \newline
In particular $q^0 \le q^2 \in \Bbb P_m \restriction \delta$ so by clause
$(\gamma)$ we have $(p^0,q^2) \in \Bbb Q$ and 
$(p^0,q^0) \le^Q (p^0,q^2)$; also $r = (0,r) \le
(p^0,q^2)$ as $r \le q^2$ together $r,(p^0,q^0)$ are compatible in
$\Bbb Q$, so $[q^1 \le r \in \Bbb P_n \restriction n \Rightarrow (p^0,q^0),r = (0,r)$ are
compatible in $\Bbb Q]$.  As $(p^0,q^0) \in \Bbb Q$ was arbitrary we are done.]
\sn
\item "{$(\zeta)$}"  if $p_1,p_2 \in \Bbb P^t$ are incompatible in
$\Bbb P^t$
then they are incompatible in $\Bbb Q$. \newline
[Why?  Look at the order of $\Bbb Q$].
\sn
\item "{$(\eta)$}"  if $\zeta < \zeta^t$ then $\bar p^t_\zeta$ is a maximal
antichain in $\Bbb Q$. \newline
[Why?  If not some $(p^*,q^*) \in \Bbb Q$ is incompatible in $\Bbb Q$ with every 
$(p^t_{\zeta,k},0)$ for $k < \omega$.  Let $n < \omega$ be a possible value
of $n(p^*,q^*)$ so $q^* \in \Bbb P_n \restriction \delta$ and there is a witness
$r^* \le q^*,r^* \in \Bbb P_n \restriction \delta^t$ for $(p^*,q^*)
\in \Bbb Q$.
\newline
By clause $(\kappa )$ in the definition of $t \in AP$ we know that for some
$r \in \Bbb P_n \cap \delta^t$ we have:
{\roster
\itemitem{ $(i)$ }  $r \in {\Cal I}^t_{\zeta,n,p^*}$
\sn
\itemitem{ $(ii)$ }  $q^*,r$ are compatible in $\Bbb P_n$
\endroster}
As $q^*,r$ are compatible and $r^* \le q^*$ also $r^*,r$ are compatible in 
$\Bbb P_n$ hence in $\Bbb P_n \cap \delta^t$, 
so by the demand on $r^*$, we have: 
$r,p^*$ are compatible in $\Bbb P^t$.  So in clause $(\iota)$ 
of definition of 
$AP$, in the definition of ${\Cal I}^t_{\zeta,n,p^*}$ for our $r$
subclause $(i)$ fails hence subclause $(ii)$ holds so there are $k,p'$ as in
subclause $(ii)$ there.  Also let $q^1 \in \Bbb P_n \restriction \delta$ be 
a common upper bound of $q^*,r$. So $r$ witness that $(p',q^1) \in
\Bbb Q$ with $n$ a possible value of $n(p',q^1)$.  
Clearly it is above $(p^*,q^*)$ and above
$p^t_{\zeta,k}$ so we are done.]
\ermn
Let $\delta^s = \delta$.  Clearly $\Bbb Q \in M_\delta$ and recall
that $M_\delta \models ``\delta$ is countable, $\Bbb P^t$ is
countable" so there is $f \in M_\delta$ such that 
$f:\Bbb Q \rightarrow X^* \cap \delta$ is a one to one (into or even onto), 
extending id$_A \cup \text{ id}_B$, and define $\Bbb P^s$ such that
$f$ is an isomorphism from $\Bbb Q$ onto $\Bbb P^s$.  
We can define $F^s_\omega,F^s_{n,\omega,\ell} \,
(n,\ell < \omega)$ extending $F^t_\omega,F^s_{n,\omega,\ell}$ as
required, e.g., $F^s_{n,\omega,\ell}((p,q)) = F^{\frak
s}_{n,m,\ell}(g)$ for some $m > n$ such that $g \in \Bbb P_m$ except
when $g=0$ then $F^s_{n,\omega,\ell}((p,0)) = F^t_{n,\omega,\ell}(p)$.
Now it is easy to check clause $(\theta)$ of the definition of
$s \in AP$, recalling $(*)$ above and clauses $(i),(\kappa)$ holds
since the construction is made in $M_\delta$.  
Lastly, let $\Gamma^s = \Gamma^t$.
\enddemo
\bn
\ub{Fact C}:  If $t^n \in AP$ and 
$t^n \le^* t^{n+1}$ for $n < \omega$ \ub{then}
there is $t$ such that \nl
$n < \omega \Rightarrow t^n \le^* t \in AP$ and
$\delta^t = \dsize \bigcup_{n < \omega} \delta^{t^n}$ and $\zeta^t =
\dsize \bigcup_{n < \omega} \zeta^{t^n}$. \newline
[Why?  Just let $\delta^t,\zeta^t$ be as above, $\Bbb P^t = \dsize 
\bigcup_{n < \omega} \Bbb P^{t^n},F^t_\omega = \dsize \bigcup_{n < \omega} 
F^{t^n}_n,F^t_{m,\omega,\ell} = \dsize \bigcup_{n < \omega} 
F^{t^n}_{m,\omega,\ell}$
and $p^t_{\zeta,k} = p^{t^n}_{\zeta,k}$ for every $n$ large enough.
Now check.]
\bn
\ub{Main Fact D}:  Assume $t \in AP,\delta^t \in E \cap \dsize
\bigcap_{n < \omega} E_n,t \in M_\delta$ and $\delta^t \in \dbcu_{n < \omega}
S^{{\Cal Y}_n}$.  \ub{Then} there is $s \in AP$
such that $t \le^* s$ and ${\underset\sim {}\to \nu_\delta}$ is actually a
$\Bbb P^s$-name (i.e. all the countably many conditions appearing in its
definition belongs to $\dsize \bigcup_{n < \omega} 
\Bbb P_m \cap \delta^t \subseteq \Bbb P^t$) and: 
\mr
\item "{$(*)$}"  \ub{if} $\Bbb P^t \subseteq_{ic} \Bbb Q$, and for 
each $\zeta < \zeta^s$ the sequence
$\bar p^t_\zeta$ is a maximal antichain of $\Bbb Q$, \ub{then} \newline
$\Vdash_{\Bbb Q} ``$there is 
$G' \subseteq \Bbb Q^{M_\delta[\underset\sim {}\to G]}
_{\underset\sim {}\to {\bar \varphi}^t_\delta}$ generic over $M_\delta[G]$
such that $(\underset\sim {}\to \eta[\underset\sim {}\to G])[G'] =
\underset\sim {}\to \nu_\delta"$. 
\ermn
[Why?  Chase arrows.]
\mn
So we can choose $t_\varepsilon \in AP$ 
by induction on $\varepsilon < \omega_1$ such
that $t^\varepsilon$ is $\le^*$-increasing continuous, 
$\delta^{t^{\varepsilon +1}} > \delta^{t^\varepsilon}$, and if $t^\varepsilon
\in M_{(\delta^{t^\varepsilon})},\delta^{t^\varepsilon} \in \dsize
\bigcap_{n < \omega} E_n \cap E \cap \dbcu_{n < \omega} S^{{\Cal Y}_n}$ 
then $t^{\varepsilon +1}$ is gotten by
Fact D.  No problem to carry this ($\varepsilon = 0$ by Fact A,
$\varepsilon = \varepsilon_1 + 1$ by Fact D if possible and by Fact B if not;
lastly, if $\varepsilon$ is a limit ordinal, use Fact C).

Now let $\Bbb P^\omega = \dsize \bigcup_{\varepsilon < \omega_1} 
\Bbb P^{t^\varepsilon}$ and it should be clear how to define ${\Cal
Y}^\omega$; now check the requirements. \hfill$\square_{\scite{12.5}}$\margincite{12.5}
\bn
\definition{\stag{12.6} Definition}  Let $\bar C^* = \langle C^*_\delta:
\delta < \omega_2$ a limit ordinal$\rangle$ (and $C^*_\alpha = \emptyset$ 
otherwise) be a square and $\bar X^* = \langle X^*_i:i < \omega_1 \rangle$ be
an increasing sequence of subsets of $\omega_1,|X^*_i \backslash 
\dsize \bigcup_{j < i} X^*_j| = \aleph_1,X^*_{\omega_1} =
\dsize \bigcup_{i < \omega_1} X^*_i$.

We say that $\langle (\Bbb P_i,{\Cal Y}_i,f_i,\bar M_i):
i < \alpha \rangle$ is a 
$(\bar C^*,\bar X^*)$-iteration (we omit $\bar M^i$ and write $(\bar
M,\bar C^*,\bar X^*)$-iteration if $i < \alpha \Rightarrow \bar M^i
\le M$ or an $\bar M$-iteration when 
$\bar C^*,\bar X$ are clear from context) \ub{if}:
\mr
\item "{$(a)$}"  $(\Bbb P_i,{\Cal Y}_i,\bar M^i) 
\in IS$ is $<^*$-increasing and Dom$(\bar M^i) = S^{{\Cal Y}_i}$
\sn
\item "{$(b)$}"  $f_i$ is a one to one function from $\Bbb P_i$ onto
$X^*_{\text{otp}(C^*_\alpha)}$, and let $(\Bbb P'_i,{\Cal Y}'_i)$ be such that 
$f_i$ maps $(\Bbb P_i,{\Cal Y}_i)$ to $(\Bbb P'_i,{\Cal Y}'_i)$
\sn
\item "{$(c)$}"  if $j \in \text{ acc}(C_i)$ then $f_j \subseteq f_i$
\sn
\item "{$(d)$}"  if cf$(i) = \aleph_0$ and $i = \text{ sup acc}(C^*_i)$ then
$(\Bbb P'_i,{\Cal Y}'_i)$ is gotten 
from $\langle (\Bbb P'_j,{\Cal Y}'_j):j \in 
\text{ acc}(C^*_\delta) \rangle$ as in the proof of \scite{12.5} 
(using $\langle X^*_j:j \in \text{ acc}(C^*_i)\rangle,
X^*_{\text{otp}(C^*_i)}$ instead of $\langle X_n:n < \omega \rangle,
X_\omega$ so acc$(C^*_i)$ replace $\omega$ and we generate \nl
$\langle t^i_\alpha:\alpha < \omega_1 \rangle$ and by it 
define $(\Bbb P'_i,{\Cal Y}'_i)$ hence $(\Bbb P_i,{\Cal Y}_i)$
\sn
\item "{$(e)$}"  in clause (d), assume $\delta = \text{ otp}(C^*_i),\langle
(\Bbb P'_j,{\Cal Y}'_j) \restriction \delta:j \in \text{ acc}(C^*_i) \rangle \in 
M_\delta$ and for $j_1 < j_2$ from acc$(C^*_i)$ the ordinal $\delta$ 
belongs to the club \nl
$\{\alpha < \omega_1:\alpha \text{ limit closed under the functions }
F^{j_1} \text{ and } F^{j_1,j_2}_\ell$ (see clause (f) below)$\}$ and
$\delta^{t^{j_1}_\delta} = \delta$.
Let $t^i_* \in AP$ be defined by $\delta^i_* = \text{ otp}(C^*_i),
\Bbb P^{t^i_*} = \cup \{\Bbb P'_j 
\restriction \delta:j \in \text{ acc}(C^*_i)\},
F^{t^i_*}_\omega = \cup\{F^{j_1} \restriction \delta_{i,j_1} \in
\text{ acc}(C^*_i)\},F^{t^i_*}_{j_1,\omega,\ell} = \cup\{F^{j_1,j_2}_\ell:
j_2 \in \text{ acc}(C^*_i) \backslash j_1\}$ and let
$\Gamma^{t^i_*}$ be empty.  If $t^i_* \in M_\delta$ then let $t^i_0$ be
gotten from $t^i_*$ as in Fact (D).
\sn
\item "{$(f)$}"  $F^j$ is a (partial) two-place function from
$X^*_{\text{otp}(C^*_j)}$ to itself such that $F^j(p,q)$ is the 
$<$-first common upper bound of $p$ and $q$ in $\Bbb P'_j$ and if $j_1 \in 
\text{ acc}(C^*_{j_2})$ then we have $\langle F^{j_1,j_2}_n(p):n < \omega 
\rangle$ is a maximal antichain of $\Bbb P'_{j_1}$ satisfying: for
each $n$, either $F^{j_1,j_2}_n(p)$ is incompatible with $p$ in $\Bbb
P'_{j_2}$ \ub{or} $p$ is compatible with $r$ in $\Bbb P'_{j_2}$
wherever $\Bbb P'_{j_1} \models F^{j_1,j_2}_n(p) \le r$.
\endroster
\enddefinition
\bn
\proclaim{\stag{12.7} Claim}  (iteration at limit)  1) Assume
$\langle (\Bbb P_i,{\Cal Y}_i,f_i):i < \zeta \rangle$ is a 
$(\bar M,\bar C^*,\bar X^*)$-iteration 
where $\zeta < \omega_2$ is a limit ordinal.  \ub{Then}
\mr
\item "{$(a)$}"  we can find 
$(\Bbb P_\zeta,{\Cal Y}_\zeta,f_\zeta)$ such that
$\langle (\Bbb P_i,{\Cal Y}_i,f_i):i < \zeta + 1 \rangle$ 
is an $\bar M$-iteration 
\sn
\item "{$(b)$}"  if $S \subseteq S^*,i < \zeta \Rightarrow S^{{\Cal Y}_i} 
\subseteq S$ mod ${\Cal D}_{\bar M}$, \ub{then} we can demand 
$S^{{\Cal Y}_\zeta} = S$.
\endroster
\endproclaim
\bn
\demo{Proof}  If cf$(\zeta) = \aleph_0$ we use \scite{12.5} but taking
care of clause (e), this just dictates to us how to start the
induction there.  If cf$(\zeta) = \aleph_1$, then by the square bookkeeping (see clause (e) in
Definition \scite{12.6}) our work is done
(using $f_\zeta = \cup \{f_\xi:\xi \in \text{ acc}(C_\zeta)\}$).
\hfill$\square_{\scite{12.7}}$\margincite{12.7}
\enddemo
\bigskip

\proclaim{\stag{12.7A} Claim}  1) Assume
\mr
\item "{$(a)$}"  ${\Cal Y} = (S,\bar \Phi,{\underset\sim {}\to {\bar \eta}},
{\underset\sim {}\to {\bar \nu}})$ is a 1-commitment on the forcing notion
$\Bbb P \in {\Cal H}(\aleph_2)$ for $\bar M$
\sn
\item "{$(b)$}"  $G_{\Bbb P} \subseteq \Bbb P$ is generic 
over $\bold V,{\bar \nu}^0 =
\langle \nu^0_\alpha:\alpha \in S \rangle$ where $\nu^0_\alpha = 
{\underset\sim {}\to \nu_\alpha}[G_{\Bbb P}]$, \nl
$\bar M^1 = \bar M[G_{\Bbb P}] 
= \langle M_\delta[f'',G_{\Bbb P}]:\delta \in S^* \rangle$ for some
one to one function $f$ from $\Bbb P$ into $\omega_1$
\sn
\item "{$(c)$}"  in $\bold V[G_{\Bbb P}],{\Cal Y}^1 = (S^1,\bar \Phi^1,
{\underset\sim {}\to {\bar \eta}^1},\bar \nu^1)$ is a 0-commitment, \nl
$S \subseteq S^1$ {\text {\rm mod\/}} 
${\Cal D}_{\bar M[G_{\Bbb P}]},\bar \Phi^1 \restriction
(S \cap S^1) = \bar \Phi \restriction (S \cap S_1)$, \nl
${\underset\sim {}\to {\bar \eta}^1} \restriction (S \cap S_1) =
{\underset\sim {}\to {\bar \eta}} \restriction (S \cap S_1),\bar \nu^1
\restriction (S \cap S_1) = \bar \nu^0 \restriction (S \cap S_1)$ and \nl
$(S^1,\bar \Phi^1,{\underset\sim {}\to {\bar \eta}^1}) \in V$
\sn
\item "{$(d)$}"  in $\bold V[G_{\Bbb P}],\Bbb Q$ 
is a forcing notion satisfying the $0$-commitment ${\Cal Y}^1$ for
$\bar M^1$.
\ermn
\ub{Then} for some $\Bbb P$-name ${\underset\sim {}\to {\Bbb Q}}$ 
and 1-commitment ${\Cal Y}^2$ we have:
\mr
\item "{$(a)$}"  $(\Bbb P,{\Cal Y}) \le^* (\Bbb P * {\underset\sim
{}\to {\Bbb Q}},{\Cal Y}^2)$
\sn
\item "{$(b)$}"  $S^{{\Cal Y}^2} = S^1,\Phi^{{\Cal Y}^2} = \bar \Phi^1,
{\underset\sim {}\to {\bar \eta}^{{\Cal Y}_2}} = 
{\underset\sim {}\to {\bar \eta}^1},{\underset\sim {}\to {\bar \nu}}[G_P] =
\bar \nu^1$
\sn
\item "{$(c)$}"  ${\underset\sim {}\to {\Bbb Q}}[G_{\Bbb P}] = \Bbb Q$.
\ermn
2) If for every $G_{\Bbb P} \subseteq \Bbb P$ generic over $\bold V$ 
there are $\Bbb Q$ satisfying some $\psi_1$ and \nl
$(S^1,\bar \Phi^1,{\underset\sim {}\to {\bar \eta}^1},
\bar \nu) \in \bold V[G_{\Bbb P}]$ 
as above satisfying some $\psi_2$, \ub{then} we can
demand
\mr
\item "{$(d)$}"  $\Vdash_{\Bbb P} ``{\underset\sim {}\to {\Bbb Q}}
[{\underset\sim {}\to G_{\Bbb P}}],{\Cal Y}^2$ 
as above satisfies $\psi_1,\psi_2$ respectively".
\ermn
3) We may allow $\langle (\bar \varphi_\alpha,
{\underset\sim {}\to {\bar \eta}
_\alpha}):\alpha \in S^1 \rangle$ be a sequence of $\Bbb P$-names 
and even $(\Bbb P * {\underset\sim {}\to {\Bbb Q}})$-names.
\endproclaim
\bigskip

\demo{Proof}  Straight.
\enddemo
\bigskip

\proclaim{\stag{12.8} Claim}  (iteration in successor case: increase
the commitment).

Assume $\langle (\Bbb P_i,{\Cal Y}_i,f_i):i < \zeta \rangle$ is an 
$\bar M$-iteration and
$\zeta = \xi + 1,S^{{\Cal Y}_\xi} \subseteq S \subseteq S^*,S
\subseteq {\text {\rm Dom\/}}(\bar M)$ and
$\langle ({\underset\sim {}\to \varphi_\alpha},
{\underset\sim {}\to \eta_\alpha}):\alpha \in S 
\backslash S^{{\Cal Y}_\xi} \rangle$
is as required in Definition \scite{12.4}.  Lastly 
${\underset\sim {}\to Z_\alpha} \subseteq {}^\omega 2$ is a 
$\Bbb P_\xi$-name of 
a positive set for $({\underset\sim {}\to {\bar \varphi}_\alpha},
{\underset\sim {}\to \eta_\alpha})$ for every such $\alpha$.

\ub{Then} we can find $(\Bbb P_\zeta,{\Cal Y}_\zeta,f_\zeta)$ such that
\mr
\item "{$(i)$}"  $\langle (\Bbb P_i,{\Cal Y}_i,f_i):
i < \zeta +1 \rangle$ is an $\bar M$-iteration
\sn
\item "{$(ii)$}"  $\Bbb P_\zeta = \Bbb P_\xi,S^{{\Cal Y}_\xi} = S,
({\underset\sim {}\to {\bar \varphi}^{{\Cal Y}_\xi}_\alpha},
{\underset\sim {}\to \eta^{{\Cal Y}_\xi}_\alpha}) =
({\underset\sim {}\to {\bar \varphi}_\alpha},
{\underset\sim {}\to \eta_\alpha})$ if $\alpha \in S \backslash 
S^{{\Cal Y}_\xi}$.
\endroster
\endproclaim
\bigskip

\demo{Proof}  Straight.
\enddemo
\bigskip

\proclaim{\stag{12.9} Claim}  (iteration at successor: increasing the forcing)
\newline
Suppose
\mr
\item "{$(a)$}"  $(\Bbb P,{\Cal Y}) \in IS$ and the set of 
elements of $\Bbb P$ is $X_i$ (the $X^*_j$'s as in \scite{12.6})
\sn
\item "{$(b)$}"  ${\underset\sim {}\to {\Bbb Q}}$ is 
a $\Bbb P$-name satisfying, for every
$G \subseteq \Bbb P$ generic over $\bold V$, the following:
{\roster
\itemitem{ (i) }  ${\underset\sim {}\to {\Bbb Q}}[G]$ is a forcing notion with set of
elements beings $X_{i+1} \backslash X_i$
\itemitem{ (ii) }  $\bigl\{\delta < \omega_1$: if $\Bbb P \restriction \delta \in
M_\delta$ and $G \cap \delta$ is a generic subset of \nl

$\qquad \qquad \qquad \Bbb P \restriction \delta,{\underset\sim {}\to
{\Bbb Q}}[G] 
\restriction \delta \in 
M_\delta[G \cap \delta]$ and ${\underset\sim {}\to \nu_\delta}[G]$ is 
forced to be \nl

$\qquad \qquad \qquad$ generic for 
$\bigl( (\Bbb Q^{{\underset\sim {}\to \varphi_\delta}[G]})^{M_\delta[G]},
{\underset\sim {}\to \eta_\delta}[G] \bigr) \bigr\} \in {\Cal D}_{\bar M[G]}$
\endroster}
\ermn
\ub{Then} we can find $(\Bbb P^+,{\Cal Y}^+)$ 
such that $(\Bbb P,{\Cal Y}) \le^* 
(\Bbb P^+,{\Cal Y}^+) \in IS$ and the $\Bbb P$-name 
$\Bbb P^+/{\underset\sim {}\to G_{\Bbb P}}$ is equivalent to
${\underset\sim {}\to {\Bbb Q}}[{\underset\sim {}\to G_{\Bbb P}}]$.
\endproclaim
\bn
\demo{Proof}  Straight.
\enddemo
\bn
\demo{\stag{12.10} Conclusion}  Assume $(\bar C^*,\bar X^*)$ is as in
\scite{12.5}.  
Let $\Phi$ be a set of definitions of forcing
notions with some real parameters, and
$\langle S^*_i:i < \omega_2 \rangle$ is as in \scite{12.3}
for ${\Cal D}_{\bar M}$.

We can find $\langle (\Bbb P_i,{\Cal Y}_i,f_i,\bar M^i):
i < \omega_2 \rangle$ such that
\mr
\item "{$(a)$}"  it is an $(\bar C^*,\bar X^*)$-iteration
\sn
\item "{$(b)$}"  $\Bbb P = \dsize \bigcup_{i < \omega_2} \Bbb P_i$ 
is a c.c.c. forcing notion of cardinality $\aleph_2$ (so in 
$\bold V^{\Bbb P},2^{\aleph_0} \le \aleph_2)$ and
except in degenerated cases equality holds)
\sn
\item "{$(c)$}"  $S^{{\Cal Y}_i} = S^*_i$ from \scite{12.3}(3)
\sn
\item "{$(d)$}"  if in $\bold V^{{\Bbb P}_i}$ we have 
$(\bar \varphi,\underset\sim {}\to \eta)$
is a case of $\Phi$ as in \scite{12.4}, moreover \newline 
$\Vdash_{{\Bbb P}_i} ``\{\delta \in S^*_{i+1} \backslash S^*_i:
M^{i+1}_\delta[f''_i,{\underset\sim {}\to G_{{\Bbb P}_i}}]"
\models ``(\bar \varphi,\underset\sim {}\to \eta)$ as required in
\scite{12.4}$"\} \in {\Cal D}^+_{{\bar M}^{i+1}}$ \newline
(even less with more bookkeeping) and $Z \subseteq 
({}^\omega 2)^{\bold V^{\Bbb P}}$ is
positive for $(\bar \varphi,\eta)$, \ub{then}
\sn
{\roster
\itemitem{ $(\alpha)$ }  $\{\delta \in S^{{\Cal Y}_{i+1}} \backslash
S^{{\Cal Y}_i}:(\bar \varphi,
{\underset\sim {}\to \eta_\delta})/G_{{\Bbb P}_i} = (\bar \varphi,\underset\sim {}\to
\eta)$ and ${\underset\sim {}\to \nu_\delta}[G_{{\Bbb P}_i}] \in Z\} \in
{\Cal D}^+_{\bar M}$, in fact the set is forced to include such 
old set (from $\bold V$) by this we can get
\sn
\itemitem{ $(\beta)$ }  for some $j > i,\delta \in S^{{\Cal Y}_{j+1}} 
\backslash S^{{\Cal Y}_0}_j \Rightarrow (\bar \varphi,
{\underset\sim {}\to \eta_\delta})/G_{{\Bbb P}_i} = (\bar \varphi,
\underset\sim {}\to \eta),{\underset\sim {}\to \nu_\delta}
[G_{{\Bbb P}_i}] \in Z$
\endroster}
\sn
\item "{$(e)$}"  if $H$ is a pregiven function such that for every $i <
\omega_2$ and $(\Bbb P,{\Cal Y},\bar M)$ satisfying \newline
$(\Bbb P_i,{\Cal Y}_i) \le^* (\Bbb P,{\Cal Y}) \in IS$ such that
$S^{\Cal Y} = {\Cal S}_i$ we have
$(\Bbb P,{\Cal Y}) \le^* H(\Bbb P,{\Cal Y}) \in IS$ 
such that $H(\Bbb P,{\Cal Y},\bar M)$ satisfies the demands from (a) +
(c) on $(\Bbb P_{i+1},{\Cal Y}_{i+1},\bar M^{i+1})$, \ub{then} 
we can demand 
$(\exists^{\aleph_2}j)
[(\Bbb P_{j+1},{\Cal Y}_{j+1}) = H(\Bbb P_j,{\Cal Y}_j)]$; moreover,
if $S^* \subseteq \omega_2$ is stationary we can demand $\{j \in
S:(\Bbb P_{j+1},{\Cal Y}_{j+1}) = H(\Bbb P_j,{\Cal Y}_n)\}$ is
stationary. \nl
(Of course, we can promise this for $\aleph_2$ such functions).
\endroster
\enddemo
\bn
\demo{Proof}  Put together the previous claims.  (Concerning clause
(e) \wilog \, $\{i < \omega_1:\text{otp}(C^*_i)=0\}$ is stationary) so
in those stages we have no influence of clause (e) of \scite{12.6};
anyhow the influence of \scite{12.6}(e) is minor.
\enddemo
\bn
\demo{\stag{12.11} Discussion}  We discuss here some possible extensions.
\enddemo
\bigskip

\proclaim{\stag{12.10a} Claim}  Assume $\langle S_i:i < \omega_2
\rangle$ is a sequence of pairwise almost disjoint stationary subsets
of $\omega_1$, each with diamond and $i < j \Rightarrow S_i \subseteq
S^+_i \, {\text {\rm mod\/}} \, {\Cal D}_{\omega_1}$, 
so $S^+_i \subseteq \omega_1$ and $S_i \cap S^+_i = \emptyset$.

\ub{Then} in the following game the between the bookkeeper and the forcer,
the bookkeeper has a winning strategy. \nl
A Play last $\omega_2$ moves, before the $\alpha$-th move a sequence
$\langle(\Bbb P_i,\Bbb Q_i,
{\underset\sim {}\to {\bar M}_i},{\Cal Y}):
i < \alpha \rangle$ is defined such that
\mr
\item "{$(a)$}"  $\Bbb P_i$ a c.c.c. forcing notion of cardinality
$\aleph_1$, say $\subseteq {\Cal H}_{< \aleph_1}(\aleph_2)$
\sn
\item "{$(b)$}"  $\Bbb Q_i$ is a $\Bbb P_i$-name of a forcing notion
of cardinality $\le \aleph_1$, say $\subseteq \omega_1$
\sn
\item "{$(c)$}"  $\Bbb P_i$ is $\lessdot$-increasing
\sn
\item "{$(d)$}"  $\Bbb P_{i+1},\Bbb P_i * {\underset\sim {}\to {\Bbb Q}_i}$
are isomorphic over $\Bbb P_i$
\sn
\item "{$(e)$}"  ${\underset\sim {}\to {\bar M}_i}$ is a $\Bbb
P_i$-name of an $S_i$-oracle
\sn
\item "{$(f)$}"  ${\underset\sim {}\to Y_i}$ is a $\Bbb P_i$-name of
a $S_i$-commitment.
\ermn
In the $i$-th move:

\block
the bookkeeper chooses $\Bbb P_i$ and a $\Bbb P_i$-name
$({\underset\sim {}\to M^+_i},{\underset\sim {}\to {\Cal Y}^+})$ of
an $S^+_i$-oracle and $0$-commitment the forces choose
${\underset\sim {}\to {\Bbb Q}_i}$ and $({\underset\sim {}\to
{\bar M}_i},{\underset\sim {}\to {\Cal Y}_i}),\Bbb P_i$-names such
that ${\underset\sim {}\to {\Bbb Q}_i}$ satisfies
$({\underset\sim {}\to M^+_i},{\underset\sim {}\to {\Cal Y}^+_i})$
and $({\underset\sim {}\to {\bar M}_i},{\underset\sim {}\to {\Cal
Y}_i})$.
\endblock
\mn
In the end the bookkeeper wins if

$$
i < j < \omega_2 \Rightarrow \Bbb P_j/\Bbb P_i \text{ satisfies }
({\underset\sim {}\to {\bar M}_i},{\underset\sim {}\to {\Cal
Y}_i}).
$$
\endproclaim
\bigskip

\demo{Proof}  Similar to earlier proofs.
\enddemo
\bn
We give an easy criterion for existence. \nl
The following uses more from \cite{Sh:630}.
\proclaim{\stag{12.12} Claim}  Assume
\mr
\item "{$(a)$}"  $(\Bbb P,\le,\le_n)_{n < \omega}$ is a definition of
a forcing notion satisfying condition $A$ of Baumgartner with $\le_n$ as
witness and $ZFC^-_*$ says this, in a way preserved by suitable forcing
\sn
\item "{$(b)$}"  ${\Cal Y}$ is a 0-commitment,
\sn
\item "{$(c)$}"  $\Bbb P$ is absolutely nep such that for each $\alpha \in
S^{\Cal Y}$ it is $\le_n$-purely \nl
$I_{{\Bbb Q}^{\bar \varphi_\alpha}}$-preserving, i.e.
{\roster
\itemitem{ $(*)$ }  if $M$ is a $\Bbb P$-candidate and a 
$\Bbb Q^{{\bar \varphi}_\alpha}$-candidate, $p \in \Bbb P^M,n < \omega$
and $q \in (\Bbb Q^{{\bar \varphi}_\alpha})^M$
\ub{then} for some $p',\eta,\nu$ we have $p \le_n p' \in \Bbb P,p'$ is
$\langle M,\Bbb P \rangle$-generic and $\nu$ is
$(\Bbb Q^{\bar \varphi_\alpha},{\underset\sim {}\to \eta_\alpha})$-generic over
$M$ satisfying $q$ (check def) and $p' \Vdash_{\Bbb P} ``\nu$ is
$(\Bbb Q^{\bar \varphi_\alpha},{\underset\sim {}\to \eta_\alpha})$-generic over
$M \langle {\underset\sim {}\to G_{\Bbb P}} \cap P^M \rangle"$.
\endroster}
\ermn
\ub{Then} there is a c.c.c. forcing notion $\Bbb P' \subseteq \Bbb P$ (not necessarily
$\Bbb P' \lessdot \Bbb P$) satisfying the 0-commitment ${\Cal Y}$ and \nl
$\Vdash_{\Bbb P}$ ``for a club of $\delta < \omega_1,\varphi(
\underset\sim {}\to \nu,{\underset\sim {}\to \eta^*_\delta})$".
\endproclaim
\bigskip

\remark{Remark}  1) Why the $\varphi_\delta$'s?  We hope it helps, for example in the
following; suppose we are given
$f:\Bbb R \rightarrow \Bbb R$, we like to force 
$A \subseteq \Bbb R$ which is not
in $I_{{\Bbb Q}^{\bar \varphi_\alpha}}$ and on which the function $f$ is continuous; 
i.e. to force a
continuous ${\underset\sim {}\to f^*}$ such that $\{\eta \in {}^\omega 2:
{\underset\sim {}\to f^*}(\eta) = f(\eta)\} \in (I^{ex}_{\Bbb Q})^+$.  So not only
do we like to find $q \Vdash ``\eta_\delta$ is 
$({\underset\sim {}\to {\Bbb Q}_\delta},{\underset\sim {}\to \eta_\delta})$-generic 
over $M_\delta[{\underset\sim {}\to G_{\Bbb P}}]"$ but also 
$q \Vdash_{{\Bbb P}'} ``\underset\sim {}\to f (\eta_\delta) = f(\eta_\delta)"$.  
This is what $\bar \varphi$ says.
\endremark
\bigskip

\demo{Proof}  We choose by induction on $\alpha < \omega_1$, a pair
$(\Bbb P_\alpha,\Gamma_\alpha)$ such that:
\mr
\item "{$(\alpha)$}"  $\Bbb P_\alpha \subseteq \Bbb P$ is countable
\sn
\item "{$(\beta)$}"  $\Gamma_\alpha$ is a countable
family of predense subsets of $\Bbb P_\alpha$
\sn
\item "{$(\gamma)$}"  if ${\Cal I} \in \Gamma_\alpha$ and $p \in \Bbb P_\alpha$
and $n < \omega$ then for some $q$ we have $p \le_n q \in \Bbb P_\alpha$ and
${\Cal I}$ is predense above $q$ in $\Bbb P$
\sn
\item "{$(\delta)$}"  $\Bbb P_\alpha$ is increasing continuous in $\alpha$
\sn
\item "{$(\varepsilon)$}"  $\Gamma_\alpha$ is increasing continuous in
$\alpha$.
\endroster
\enddemo
\bn
\ub{Case 1}:  $\alpha = 0$.  Trivial.
\bn
\ub{Case 2}:  $\alpha = \beta +1,\beta$ non-limit or $(\Bbb P_\beta,\Gamma_\beta)
\notin M_\beta$.

Let $(\Bbb P_\alpha,\Gamma_\alpha) = (\Bbb P_\beta,\Gamma_\beta)$.
\bn
\ub{Case 3}:  $\alpha$ limit.

Let $(\Bbb P_\alpha,\Gamma_\alpha) = (\dbcu_{\beta < \alpha} \Bbb P_\beta,
\dbcu_{\beta < \alpha} \Gamma_\beta$.
\bn
\ub{Case 4}:  $\alpha = \delta+1$ where $\delta$ is a limit ordinal
and $(\Bbb P_\delta,\Gamma_\delta) \in M_\delta$.

We can find $g \subseteq \text{ Levy}(\aleph_0,|\Bbb P|)^{M_\delta}$, generic
over $M_\delta$ such that $\eta^*_\delta$ is still \nl
$\Bbb Q_\delta$-generic over $M_\delta[g]$ (see \cite{Sh:630},\S6).

In $M_\delta[g]$ we define $\Bbb P^+_\delta = \{p:M_\delta[g] \models
p \in \Bbb P \text{ and } {\Cal I} \in \Gamma_\alpha \Rightarrow {\Cal I} 
\text{ predense above } p'\}$
using the induction hypothesis, as in $M_\delta[g]$ the sets $\Gamma_\delta$
is countable:
\mr
\item "{$(*)$}"  for every $p \in \Bbb P_\delta$ and $n < \omega$ there is
$p' \in \Bbb P^+_\delta$ such that $\Bbb P \models p \le_n p'$.
\ermn
Again by \cite{Sh:630},\S6 for every
$n < \omega$ and $p \in \Bbb P^+_\delta$, there is $q_{p,n} \in \Bbb P$ such that
$p \le_n q_{p,n} \in \Bbb P,q_{p,n}$ is $(M_\delta[g],\Bbb Q)$-generic and
$q_{p,n} \Vdash_{\Bbb P} ``\eta_\delta$ is a $(\Bbb Q_\delta,
{\underset\sim {}\to \eta_\delta})$-generic real over $M_\delta[g]
[{\underset\sim {}\to G_{\Bbb P}}]$ 
and $\varphi_\delta(\underset\sim {}\to \nu,
\eta_\delta)"$.

Let $\Bbb P_{\delta +1} = \Bbb P_\delta \cup \{q_{p,n}:p \in \Bbb P^+_\delta \text{ and }
n < \omega\}$ and $\Gamma_{\delta + 1} = \Gamma_\delta \cup 
\{{\Cal I}_\delta\}$ \nl
where ${\Cal I}_\delta = 
\{q_{p,n}:p \in \Bbb P^+_\delta \text{ and } n < \omega\}$.
\hfill$\square_{\scite{12.12}}$\margincite{12.12}
\newpage

     \shlhetal 

\newpage
    
REFERENCES.  
\bibliographystyle{lit-plain}
\bibliography{lista,listb,listx,listf,liste}

\enddocument